\documentclass[12pt]{amsart}
\usepackage[osf,sc]{mathpazo}
\usepackage{amssymb}

\usepackage{geometry}
\usepackage{bm}
\usepackage{graphicx}
\usepackage{multirow}
\usepackage{hyperref}
\hypersetup{
	colorlinks=true, 
	linktoc=all,     
	linkcolor=blue,
	citecolor=red,
	filecolor=black,
	urlcolor=blue	
}
\usepackage{enumerate}
\usepackage[inline]{enumitem}
\makeatletter
\newcommand{\inlineitem}[1][]{%
	\ifnum\enit@type=\tw@
	{\descriptionlabel{#1}}
	\hspace{\labelsep}%
	\else
	\ifnum\enit@type=\z@
	\refstepcounter{\@listctr}\fi
	\quad\@itemlabel\hspace{\labelsep}%
	\fi} \makeatother
\newcommand{\ga}{\alpha}

\newcommand{\gd}{\delta}

\newcommand{\gz}{\zeta}

\newcommand{\gl}{\lambda}
\newcommand{\gm}{\mu}
\newcommand{\gn}{\nu}
\newcommand{\gx}{\xi}

\newcommand{\gp}{\pi}

\newcommand{\gr}{\rho}

\newcommand{\gf}{\phi}

\newcommand{\gc}{\psi}

\newcommand{\Gd}{\Delta}

\newcommand{\Gl}{\Lambda}

\newcommand{\Gs}{\Sigma}

\newcommand{\Gom}{\Omega}

\newcommand{\ugl}{\ul{\gl} = (\gl_{1}^{\gr _1}>\gl_{2}^{\gr _2}>\gl_{3}^{\gr _3}>\ldots>\gl_{k}^{\gr _k})}

\newcommand{\grpp}{\mcl{A}_{\ul{\gl}}}
\newcommand{\grppp}{\mcl{A}_{\ul{\gl}^{'}}}
\newcommand{\grpppp}{\mcl{A}_{\ul{\gl}^{''}}}
\newcommand{\grpppI}{\mcl{A}_{\ul{\gl^{'}/I}}}
\newcommand{\autgp}{\mcl{G}_{\ul{\gl}}}
\newcommand{\autgpp}{\mcl{G}_{\ul{\gl}^{'}}}
\newcommand{\autgppp}{\mcl{G}_{\ul{\gl}^{''}}}

\newcommand{\hgrpp}{\mcl{H}^0_{\ul{\gl}}}

\newcommand{\grmpp}{\mcl{A}_{\ul{\gm}}}
\newcommand{\grmppp}{\mcl{A}_{\ul{\gm}^{'}}}
\newcommand{\grmpppp}{\mcl{A}_{\ul{\gm}^{''}}}
\newcommand{\grmpppI}{\mcl{A}_{\ul{\gm^{'}/I}}}
\newcommand{\autmgp}{\mcl{G}_{\ul{\gm}}}

\newcommand{\hgrmpp}{\mcl{H}^0_{\ul{\gm}}}

\newcommand{\grnpp}{\mcl{A}_{\ul{\gn}}}
\newcommand{\grnppp}{\mcl{A}_{\ul{\gn}^{'}}}
\newcommand{\grnpppp}{\mcl{A}_{\ul{\gn}^{''}}}
\newcommand{\grnpppI}{\mcl{A}_{\ul{\gn^{'}/I_1}}}
\newcommand{\autngp}{\mcl{G}_{\ul{\gn}}}

\newcommand{\hgrnpp}{\mcl{H}^0_{\ul{\gn}}}

\newcommand{\subs}{\subset}

\newcommand{\sbnq}{\subsetneq}

\newcommand{\bs}{\backslash}
\newcommand{\fs}{/}

\newcommand{\nin}{\notin}

\newcommand{\ti}{\tilde}

\newcommand{\mbb}{\mathbb}

\newcommand{\mcl}{\mathcal}

\newcommand{\ul}{\underline}
\newcommand{\ol}{\overline}

\newcommand{\us}{\underset}
\newcommand{\os}{\overset}

\newcommand{\lra}{\longrightarrow}
\newcommand{\llra}{\longleftrightarrow}

\newcommand{\N}{\mbb N}

\newcommand{\R}{\mcl R}

\newcommand{\RR}[1]{\R/\pi^{#1}\R}

\newcommand{\ra}{\rightarrow}
\newcommand{\Ra}{\Rightarrow}
\newcommand{\Llra}{\Longleftrightarrow}

\newcommand{\es}{\emptyset}

\newcommand{\fo}[3]{
	\begingroup
	{\fontsize{#1}{#2}\selectfont {#3}}
	\endgroup
}
\newcommand{\equ}[1]{%
	\begin{equation*}
		#1
	\end{equation*}
}
\newcommand{\equa}[1]{%
	\begin{equation*}
		\begin{aligned}
			#1
		\end{aligned}
	\end{equation*}
}
\newcommand{\equan}[2]{%
	\begin{equation}
		\label{Eq:#1}
		\begin{aligned}
			#2
		\end{aligned}
	\end{equation}
}

\DeclareMathOperator{\Hom}{Hom}
\DeclareMathOperator{\End}{End}
\DeclareMathOperator{\Aut}{Aut}

\newcommand{\mattwo}[4]{%
	\begin{pmatrix}
		#1 & #2\\ #3 & #4
	\end{pmatrix}
}

\theoremstyle{plain}
\newtheorem{theorem}{Theorem}[section]

\newtheorem{prop}[theorem]{Proposition}

\newtheorem{cor}[theorem]{Corollary}
\newtheorem{conj}[theorem]{Conjecture}

\newtheorem{claim}[theorem]{Claim}

\makeatletter
\def\namedlabel#1#2{\begingroup
	\def\@currentlabel{#2}%
	\label{#1}\endgroup
}
\makeatother

\newtheorem*{thmDelta}{\bf{Theorem} $\bm{\Gd}$}
\theoremstyle{definition}

\theoremstyle{remark}
\newtheorem{remark}[theorem]{Remark}
\newtheorem{example}[theorem]{Example}
\numberwithin{equation}{section}

\begin{document}

	\title[Positivity Conjecture]{On the Positivity Conjecture for Finite Abelian \MakeLowercase{p}-Groups}
	\author{C P Anil Kumar}
	\address{School of Mathematics, Harish-Chandra Research Institute, HBNI, Chhatnag Road, Jhunsi, Prayagraj (Allahabad), 211 019,  India. \,\, email: {\tt akcp1728@gmail.com}}
	\thanks{This work is done while the author is a Post Doctoral Fellow at Harish-Chandra Research Institute, Prayagraj(Allahabad).}
	\date{\sc \today}
	\vspace*{-3cm}
	\begin{abstract}
		For a partition $\ugl$ and its associated finite $\R$-module $\grpp=\us{i=1}{\os{k}{\oplus}} (\RR {\gl_i})^{\gr_i}$, where $\R$ is a discrete valuation ring, with maximal ideal generated by a uniformizing element $\gp$, having finite residue field ${\bf k}=\RR {} \cong \mbb{F}_q$, the number of orbits of pairs 
		$n_{\ul{\gl}}(q)= \mid \autgp\bs \big(\grpp\times \grpp\big)\mid$ for the diagonal action of the automorphism group $\autgp=\Aut(\grpp)$, is a polynomial in $q$ with integer coefficients. Positivity conjecture states that these coefficients are in fact non-negative. In this article, we prove this conjecture.
	\end{abstract}
\maketitle
\section{\bf{Introduction}}
The automorphism orbits in finite abelian groups have been understood quite well for over a hundred years (See G.~A.~Miller~\cite{MR1505955}, G.~Birkhoff~\cite{Birkhoff}). The combinatorics of automorphism orbits in abelian groups is also studied by K.~Dutta, A.~Prasad~\cite{MR2793603} and C.~P.~Anil Kumar~\cite{MR3261812},~\cite{MR3647154}. Some other authors who have worked on this subject are  S.~Delsarte~\cite{MR0025463}, M.~Schwachh\"{o}fer, M.~Stroppel~\cite{MR1656579}, B.~L.~Kerby and E.~Rode~\cite{MR2782608}.  In this article we mention a conjecture for finite abelian $p$-groups called the positivity conjecture about the total number of orbits of pairs. This Conjecture~\ref{conj:PositivityConjecture} is stated below in the context of finite modules over discrete valuation rings where a finite abelian $p$-group is naturally a finite module over the ring of $p$-adic integers. This was initially jointly conjectured by  C.~P.~Anil Kumar, Amritanshu Prasad in~\cite{MR3261812}.

Let $\Gl_0$ denote the set of all sequences of the form 
\equan{Partition}{\ul{\gl}=(\gl_1^{\gr_1},\gl_2^{\gr_2},\ldots,\gl_k^{\gr_k})} 
where $\gl_1>\gl_2>\ldots>\gl_k$ is a strictly decreasing sequence of positive integers and $\gr_1,\gr_2,\ldots,\gr_k$ are positive integers. We allow the case where $k=0$, resulting in the empty sequence, which we denote by $\es$. Let $\R$ be a discrete valuation ring with maximal ideal generated by a uniformizing element $\gp$ having finite residue field ${\bf k}=\RR {}$. Every finite $\R$-module, up to an isomorphism is of the form 
\equan{FiniteModule}{\grpp=\us{i=1}{\os{k}{\oplus}} (\RR {\gl_i})^{\gr_i}} for a unique $\ul{\gl}\in \Gl_0$.

Now we mention the positivity conjecture which was intially stated in~\cite{MR3261812}. 

\begin{conj}
	\label{conj:PositivityConjecture}
	Let $\ugl\in \Gl_0$ be a partition and $\R$ be a discrete valuation ring with maximal ideal generated by a uniformizing element $\gp$ having finite residue field ${\bf k}=\RR {} \cong \mbb{F}_q$. Let $\grpp=\us{i=1}{\os{k}{\oplus}} (\RR {\gl_i})^{\gr_i}$ be its associated finite $\R$-module and $\autgp$ be its automorphism group. Let $n_{\ul{\gl}}(q)$ be the number of orbits of pairs in $\grpp \times \grpp$ for the diagonal action of $\autgp$ on $\grpp \times \grpp$. Then $n_{\ul{\gl}}(q)$ is a polynomial of degree $\gl_1$ with non-negative integer coefficients.
\end{conj}

\subsection{Results Known about the Positivity Conjecture}
We already know the following.
\begin{enumerate}
	\item The number $n_{\ul{\gl}}(q)$ is a polynomial in $q$ with integer coefficients (C.~P.~Anil Kumar, Amritanshu Prasad~\cite{MR3261812}, Theorem 5.11, Page 17).
	\item The degree of $n_{\ul{\gl}}(q)$ is $\gl_1$ the highest part of $\ul{\gl}$
	(C.~P.~Anil Kumar, Amritanshu Prasad~\cite{MR3261812}, Theorem 5.11, Page 17).
	\item The polynomial $n_{\ul{\gl}}(q)=n_{\ul{\gl}^2}(q)$ where \equ{\ul{\gl}^2=(\gl_1^{\min(\gr_1,2)}>\gl_2^{\min(\gr_2,2)}>\ldots>\gl_k^{\min(\gr_k,2)})}
	(C.~P.~Anil Kumar, Amritanshu Prasad~\cite{MR3261812}, Corollary 4.5, Page 12).
\end{enumerate} 
So, for the computation of the polynomial $n_{\ul{\gl}}(q)$, any part of the partition which repeats more than twice can be reduced to two. Here in this article we prove the positivity conjecture.  

\section{\bf{Preliminaries}}
In this section we mention some preliminaries about finite $\R$-modules.
\subsection{Orbits of Elements}
For the present purposes, the combinatorial description of orbits due to K.~Dutta and A.~Prasad~\cite{MR2793603} is useful. We recapitulate the notation and some of the results in~\cite{MR2793603} here for the case of finite $\R$-modules. 

It turns out in the case of finite $\R$-modules, that, for any module $\grpp$, the $\autgp$-orbits in $\grpp$ are in bijective correspondence with certain classes of ideals in a poset $P$, which we call the fundamental poset. As a set 
\equan{FP}{P=\{(v,l)\mid v\in \N\cup \{0\}, l\in \N, 0\leq v< l\}.}
The partial order on $P$ is defined by setting 
\equ{(v,l)\leq (v',l') \text{ if and only if }v\geq v' \text{ and }l-v\leq l'-v'.}

Let $\mcl{J}(P)$ denote the lattice of order ideals in $P$. A typical element of $\grpp$ is a column vector of the form $x=(x_{\gl_i,r_i})$ where $i$ runs over the set $\{1,\ldots,k\}$, and, for each $i,r_i$ runs over the set $\{1,\ldots,\gr_i\}$. To $x\in \grpp$, we associate the order ideal $I(x)\in \mcl{J}(P)$ generated by the elements $(v(x_{\gl_i,r_i}),\gl_i)$ for all pairs $(i,r_i)$ such that $x_{\gl_i,r_i}\neq 0$ in $\RR {\gl_i}$. Here for any $x\in \grpp, v(x)$ denotes the largest $l$ for which $x\in \gp^l\grpp$ (in particular, $v(0)=\infty$).  

A key observation is the following theorem.
\begin{theorem}
	Let $\grpp,\mcl{A}_{\ul{\gm}}$ be two finite $\R$-modules. An element $y\in \mcl{A}_{\ul{\gm}}$ is a homomorphic image of $x\in \grpp$ if and only if $I(y)\subseteq I(x)$. Moreover for $x,y\in \grpp$, they lie in the same $\autgp$-orbit if and only if $I(x)=I(y)$.
\end{theorem}
Note that the orbit of $0$ corresponds to the empty ideal.

For each $\ul{\gl}\in \Gl_0$ Let $\mcl{J}(P_{\ul{\gl}})$ denote the sublattice of $\mcl{J}(P)$ consisting of ideals such that $\max I$ is contained in the set 
\equan{FPLambda}{P_{\ul{\gl}}=\{(v,l)\mid 0\leq v< l,l=\gl_i,1\leq i\leq k\}.}
Then the $\autgp$-orbits in $\grpp$ are in bijective correspondence with this set $\mcl{J}(P_{\ul{\gl}})$ of order ideals. The lattice $\mcl{J}(P_{\ul{\gl}})$ is  isomorphic to the lattice $\mcl{J}(P_{\ul{\gl}})$ of order ideals in the induced subposet $P_{\ul{\gl}}$. For each order ideal $I\in \mcl{J}(P_{\ul{\gl}})$, we use the notation 
\equ{(\grpp)^*_I=\{x\in \grpp\mid I(x)=I\}} for the orbit corresponding to $I$.

A convenient way to think about the ideals in $P$ is in terms of what we call boundaries: for each positive integer $l$ define the boundary valuation of $I$ at $l$ to be \equ{\partial_lI=\min\{v\mid \min (v,l)\in I\}.}
We denote the sequence $\{\partial_lI\}$ of boundary valuations $\partial I$ and call it the boundary of $I$.

For each order ideal $I\subset P$ let $\max I$ denote the set of maximal elements in $I$. The ideal $I$ is completely determined by $\max I$: in fact, taking $I$ to $\max I$ gives a bijection from the lattice  $\mcl{J}(P_{\ul{\gl}})$ to the set of antichains in $P_{\ul{\gl}}$.

An alternative description for $(\grpp)^*_{I}$ is 
\equa{(\grpp)^*_{I}&=\{x=(x_{\gl_i,r_i})\in \grpp\mid v(x_{\gl_i,r_i})\geq \partial_{\gl_i}I \text{ for all }\gl_i \text{ and }r_i\text{ and for }\\&(\partial_{\gl_i}I,\gl_i)\in \max I \text{ we must have }\us{r_i}{\max} \{v(x_{\gl_i,r_i})\mid \text{ for all }r_i\}=\partial_{\gl_i}I\}.}

For each ideal $I\in \mcl{J}(P_{\ul{\gl}})$ with \equ{\max I=\{(v_1,l_1),\ldots(v_s,l_s)\}} define an element $e_I$ of $\grpp$ a column vector whose co-ordinates are given by 
\equan{CanonicalForm}{x_{\gl_i,r_i}=\begin{cases}
		\ol{\gp}^{v_j}, \text{ if }\gl_i=l_j\text{ and }r_j=1,\\
		0,\text{ otherwise.}\\		
	\end{cases}	
}
For the $\R$-module $\grpp$ the functions $x\lra I(x)$ and $I\lra e_I$ induce mutually inverse bijections between the set of $\autgp$-orbits in $\grpp$ and the set of order ideals in $\mcl{J}(P_{\ul{\gl}})$.

For an ideal $I\in \mcl{J}(P_{\ul{\gl}})$, define 
\equ{(\grpp)_I=\us{J\subseteq I}{\bigsqcup}(\grpp)^*_J}
This is a submodule of $\grpp$ which is $\autgp$-invariant. The description of $(\grpp)_I$ in terms of valuations of co-ordinates and boundary valuations is as follows.
\equ{(\grpp)_I=\{x=(x_{\gl_i,r_i})\mid v(x_{\gl_i,r_i})\geq \partial_{\gl_i}I\}.} 
Such a module $(\grpp)_I$ given by an ideal $I\in \mcl{J}(P_{\ul{\gl}})$ is called a characteristic submodule. The $\autgp$ orbits in $\grpp$ are parametrized by the finite distributive lattice $\mcl{J}(P_{\ul{\gl}})$. Moreover, each order ideal $I\in \mcl{J}(P_{\ul{\gl}})$ gives rise to a characteristic submodule $(\grpp)_I$ of $\grpp$. The lattice structure of $\mcl{J}(P_{\ul{\gl}})$ gets reflected in the poset structure of the characteristice submodules $(\grpp)_I$ when they are partially ordered by inclusion. For ideals $I,J\in \mcl{J}(P_{\ul{\gl}})$, we have 
\equ{(\grpp)_{I\cup J}=(\grpp)_I+(\grpp)_J\text{ and }(\grpp)_{I\cap J}=(\grpp)_I\cap(\grpp)_J.}
The map $I\ra (\grpp)_I$ gives an isomorphism from $\mcl{J}(P_{\ul{\gl}})$ to the poset of $\autgp$-invariant characteristic submodules. In fact when the residue field ${\bf k}$ of $\R$ has at least three elements then every $\autgp$-invariant submodule is a characteristic submodule. Therefore $\mcl{J}(P_{\ul{\gl}})$ is isomorphic to the lattice of $\autgp$-invariant submodules (B.~L.~Kerby and E.~Rode~\cite{MR2782608}).

\subsection{Stabilizer of $e_I$}
By the description of $\autgp$-orbits in $\grpp$, every $\autgp$-orbit of pairs of elements $(x_1,x_2)\in \grpp\times \grpp$ contains a pair of the form $(e_I,x)$ for some $I\in \mcl{J}(P_{\ul{\gl}})$ and $x\in \grpp$. Then the $\autgp$-orbits of pairs in $\grpp\times \grpp$ which contain an element of the form $(e_I,x)$ for some $x\in \grpp$ are in bijective correspondence with $(\autgp)_I$-orbits in $\grpp$. Here we give a description of $(\autgp)_I$ which facilitates the classification of $(\autgp)_I$-orbits in $\grpp$.

The main idea here is to decompose $\grpp$ into a direct sum of two $\R$-modules (this decomposition depends on $I$):
\equan{Decomp}{\grpp=\grppp\oplus \grpppp,}
where $\grppp$ consists of those cyclic summands in the decomposition~\ref{Eq:FiniteModule} of $\grpp$ where $e_I$ has non-zero co-ordinates, and $\grpppp$ consists of remaining cyclic summands. With respect to this decomposition we have 
\equ{e_I=(e_I',0).}
The reason for introducing this decomposition is that the description of stabilizer of $e_I'$ in the automorphism group $\autgpp$ of $\grppp$ is quite nice (see Lemma $4.2$ on Page $11$ in~\cite{MR3261812} ). The stabilizer of $e_I'$ in $\autgpp$ is 
\equ{(\autgpp)_I=\{id_{\grppp}+u\mid u\in \Hom(\grppp,\grppp)\text{ satisfies }u(e_I')=0\}.}	
Every endomorphism of $\grpp$ can be written as a matrix of the form $\mattwo {g_{11}}{g_{12}}{g_{21}}{g_{22}}$ where $g_{11}:\grppp\ra \grppp,g_{22}:\grpppp\ra \grpppp,g_{12}:\grpppp\ra \grppp,g_{21}:\grppp\ra \grpppp$. The stabilizer of $e_I$ in $\autgp$ consists of matrices of the form 
\equ{\mattwo {id_{\grppp}+u}{g_{12}}{g_{21}}{g_{22}}}
where $u\in \Hom(\grppp,\grppp)$ satisfies $u(e_I')=0,g_{12}\in \Hom(\grpppp,\grppp)$ is arbitrary, $g_{21}\in \Hom(\grppp,\grpppp)$ satisfies $g_{21}(e_I')=0$ and $g_{22}\in \autgppp=\Aut(\grpppp)$ is invertible. For this see Theorem $4.4$ on page $11$ in~\cite{MR3261812}.

\subsection{The Stabilizer Orbit of an Element}

Let $(\autgp)_I$ denote the stabilizer of $e_I\in \grpp$. Write each element $x\in \grpp$ as $x=(x',x'')$ with respect to the decomposition~\ref{Eq:Decomp} of $\grpp$. Also for any $x'\in \grppp$ let $\ol{x}'$ denote the image of $x'$ in $\grpppI=\grppp\bs \R e_I'$. The partition $\ul{\gl^{'}/I}$ is completely determined by the partition $\ul{\gl}'$ and the ideal $I\cap P_{\ul{\gl}'}\in J(\mcl{P})_{\ul{\gl}'}$ (see Lemma $6.2$ on page $19$ in~\cite{MR3261812}). 

Now we describe the stabilizer $(\autgp)_I$-orbit of $x\in \grpp$. 
\begin{theorem}
	\label{theorem:StabilizerOrbit}
	Given $x,y\in \grpp, y=(y',y'')$ lies in the $(\autgp)_I$-orbit of $x\in \grpp$ if and only if the following conditions hold:
	\equan{FirstCond}{y'\in x'+(\grppp)_{I(\ol{x}')\cup I(x'')}} and 
	\equan{SecondCond}{y''\in (\grpppp)^*_{I(x'')}+(\grpppp)_{I(\ol{x}')}} 
	The size of the orbit is $\mid(\grppp)_{I(\ol{x}')\cup I(x'')}\mid\mid (\grpppp)^*_{I(x'')}+(\grpppp)_{I(\ol{x}')} \mid$.
\end{theorem}
For the proof of this result see Theorem $5.1$ on page $13$ in~\cite{MR3261812}.

\section{\bf{Statement and Proof of the Theorem in the case of Partitions with a Repeated Part}}
\label{sec:RepeatedPartCase}
We state the theorem.
\begin{thmDelta}[The Repeated Part Case]
	\namedlabel{theorem:RepeatedPartCase}{$\Gd$}
	Let $\ugl\in \Gl_0$ be such that $\gr_{i_0}= 2$ for some $1\leq i_0\leq k$. Let $\ul{\gm} =\big(\gl_1^{\gr_1}>\gl_2^{\gr_2}>\ldots>\gl_{i_0-1}^{\gr_{i_0-1}}>\gl_{i_0}^{\gr_{i_0}-1}=\gl_{i_0}^1>\gl_{i_0+1}^{\gr_{i_0+1}}>\ldots>\gl_{k-1}^{\gr_{k-1}}>\gl_k^{\gr_k}\big)$.
	Let $\ul{\gn}=\big((\gl_1-2)^{\gr_1}>(\gl_2-2)^{\gr_2}>\ldots>(\gl_{i_0-1}-2)^{\gr_{i_0-1}}\geq (\gl_{i_0}-1)^{\gr_{i_0}}\geq \gl_{i_0+1}^{\gr_{i_0+1}}>\ldots>\gl_{k-1}^{\gr_{k-1}}>\gl_k^{\gr_k}\big)$. Let  $n_{\ul{\gl}}(q)= \mid \autgp\bs \grpp\times \grpp\mid,n_{\ul{\gm}}(q)=\mid \autmgp\bs \grmpp\times \grmpp\mid,n_{\ul{\gn}}(q)= \mid \autngp\bs \grnpp\times \grnpp\mid$. Then we have 
	\equ{n_{\ul{\gl}}(q)=n_{\ul{\gm}}(q)+n_{\ul{\gn}}(q).}
\end{thmDelta}
\begin{example}
\begin{itemize}
\item If $\ul{\gl}=(\gl_1^2)=(1^2)$ then $i_0=1,\ul{\gm}=(1),\ul{\gn}=\es$.
\item If $\ul{\gl}=(\gl_1=2>\gl_2^2=1^2)$ then $i_0=2,\ul{\gm}=(2>1),\ul{\gn}=\es$.
\item If $\ul{\gl}=(\gl_1^2=2^2>\gl_2=1)$ then $i_0=1,\ul{\gm}=(2,1),\ul{\gn}=(1)$.
\item If $\ul{\gl}=(\gl_1^2=2^2>\gl_2^2=1^2)$ and $i_0=1$ then $\ul{\gm}=(2>1^2),\ul{\gn}=(1^2)$.
\item If $\ul{\gl}=(\gl_1^2=2^2>\gl_2^2=1^2)$ and $i_0=2$ then $\ul{\gm}=(2^2>1),\ul{\gn}=\es$.
\item If $\ul{\gl}=(\gl_1^2=4^2>\gl_2^2=3^2>\gl_3^2=2^2)$ and $i_0=2$ then $\ul{\gm}=(4^2>3>2^2),\ul{\gn}=(2^2\geq 2^2\geq 2^2)=(2^6)$. Alternatively, we can choose $\ul{\gn}^2=(2^{\min(6,2)})=(2^2)$ for the computation of the polynomial $n_{\ul{\gn}}(q)=n_{\ul{\gn}^2}(q)$, as any part of the partition which repeats more than twice can be reduced to two.

\end{itemize}
\end{example}
As a consequence of Theorem~\ref{theorem:RepeatedPartCase} we have the following theorem. 
\begin{theorem}
	\label{theorem:RepeatedPartsToDistinctParts}
	Let $\ugl\in \Gl_0$ be a partition. Let  $n_{\ul{\gl}}(q)= \mid \autgp\bs \grpp\times \grpp\mid$. Then we have that,
	$n_{\ul{\gl}}(q)$ is expressible as the sum of polynomials $n_{\ul{\gm}}(q)$ corresponding to partitions $\ul{\gm}$ which have all its parts distinct.
	As a consequence if $n_{\ul{\gm}}(q)$ is a polynomial with nonnegative integer coefficients for any partition $\ul{\gm}$ which has all its parts distinct then $n_{\ul{\gl}}(q)$ is a polynomial with nonnegative integer coefficients for any partition $\ul{\gl}\in \Gl_0$ in general.
\end{theorem}
\begin{proof}
	This theorem follows by applying repeatedly Theorem~\ref{theorem:RepeatedPartCase} and reducing the number of repetitions of any part in any partition, which occurs more than twice to two (refer C.~P.~Anil Kumar, Amritanshu Prasad~\cite{MR3261812}, Corollary 4.5, Page 12).	
\end{proof}
\begin{example}
Let us illustrate Theorem~\ref{theorem:RepeatedPartsToDistinctParts} with an example where $\ul{\gl}=(3^2>2^2>1^2)$.
First we have
\begin{itemize}
\item $n_{(3>2^2>1^2)}(q)=n_{(3>2>1^2)}(q)+n_{(1^2)}(q)$.
\item $n_{(3>2>1^2)}(q)=n_{(3>2>1)}(q)+n_{(1)}(q)$.
\item $n_{(2^2>1^2)}(q)=n_{(2>1^2)}(q)+n_{(1^2)}(q)$.
\item $n_{(2>1^2)}(q)=n_{(2>1)}(q)+n_{\es}(q)$.
\item $n_{(1^2)}(q)=n_{(1)}(q)+n_{\es}(q)$.
\end{itemize}
So \equa{n_{(3^2>2^2>1^2)}(q)&=n_{(3>2^2>1^2)}(q)+n_{(2^2>1^2)}(q)\\&=n_{(3>2>1^2)}(q)+n_{(1^2)}(q)+n_{(2>1^2)}(q)+n_{(1^2)}(q)\\&=
n_{(3>2>1^2)}(q)+n_{(2>1^2)}(q)+2n_{(1^2)}(q)\\&=n_{(3>2>1)}(q)+n_{(1)}(q)+n_{(2>1)}(q)+n_{\es}(q)+2(n_{(1)}(q)+n_{\es}(q))\\&=n_{(3>2>1)}(q)+n_{(2>1)}(q)+3n_{(1)}(q)+3n_{\es}(q).}
\end{example}

We prove another important theorem and state a few other required theorems before proving Theorem~\ref{theorem:RepeatedPartCase}.

\begin{theorem}
	\label{theorem:LatticeIso}
	Let $\ugl\in \Gl_0$ and $\gl_{i_0}$ be a part of $\ul{\gl}$ for some $1\leq i_0\leq k$. Let $\ul{\gm}=\big(\gl_1^{\gr_1}>\gl_2^{\gr_2}>\ldots>\gl_{i_0-1}^{\gr_{i_0-1}}>\gl_{i_0+1}^{\gr_{i_0+1}}>\ldots>\gl_{k-1}^{\gr_{k-1}}>\gl_k^{\gr_k}\big)$ and let $\ul{\gn}=\big((\gl_1-2)^{\gr_1}>(\gl_2-2)^{\gr_2}>\ldots>(\gl_{i_0-1}-2)^{\gr_{i_0-1}}\geq (\gl_{i_0}-1)^{\gr_{i_0}}\geq \gl_{i_0+1}^{\gr_{i_0+1}}>\ldots>\gl_{k-1}^{\gr_{k-1}}>\gl_k^{\gr_k}\big)$.
	\begin{enumerate}
		\item Consider the subset $\mcl{L} \subs \mcl{J}(P_{\ul{\gl}})$ of ideals given by 
		\equ{\mcl{L}=\{I\in \mcl{J}(P_{\ul{\gl}})\mid \max(I)\cap P_{(\gl_{i_0})}=\es\}.}
		Then $\mcl{L}$ is a sublattice. It is isomorphic to the lattice $\mcl{J}(P_{\ul{\gm}})$ with the lattice isomorphism being \equ{I\in \mcl{J}(P_{\ul{\gl}}) \lra I\cap P_{\ul{\gm}}\in \mcl{J}(P_{\ul{\gm}}).}   
		\item Consider the subset $\mcl{M} \subs \mcl{J}(P_{\ul{\gl}})$ of ideals given by 
		\equ{\mcl{M}=\{I\in \mcl{J}(P_{\ul{\gl}})\mid \max(I)\cap P_{(\gl_{i_0})}\neq\es\}=\mcl{J}(P_{\ul{\gl}})\bs \mcl{L}.}
		Then $\mcl{M}$ is a sublattice. It is isomorphic to the lattice $\mcl{J}(P_{\ul{\gn}})$ with the lattice isomorphism being
		\equ{\gc: I\in \mcl{M}\lra J=\langle \max(J)\rangle\in \mcl{J}(P_{\ul{\gn}})}
		where $\max(J)$ is given as follows. Let $K=\{(v-1,\gl-2)\mid (v,\gl)\in \max(I),\gl>\gl_{i_0}\}\cup \{(v,\gl)\mid (v,\gl)\in \max(I),\gl<\gl_{i_0}\}\subs P_{\ul{\gn}}$. (Note that $v-1<\gl-2$ for $(v,\gl)\in \max(I),\gl>\gl_{i_0}$.) Let $\max(I)\cap P_{(\gl_{i_0})}=\{(u,\gl_{i_0})\}$. Then 
		\begin{center}
			$\max(J)=
			\begin{cases}
				K \cup \{(u,\gl_{i_0}-1)\} \text{ if }K\cup\{(u,\gl_{i_0}-1)\} \text{ is an antichain,}\\
				K \text{ if }u=\gl_{i_0}-1 \text{ or if }K\cup\{(u,\gl_{i_0}-1)\} \text{ is not an antichain.}\\ 
			\end{cases}$
		\end{center}
	\end{enumerate}
\end{theorem}
\begin{proof}
	Assertion $(1)$ is clear and immediate. We prove assertion $(2)$. We first show that $K$ is an antichain.
	Let $(v,\gl),(w,\gm)\in \max(I),\gl>\gl_{i_0},\gm<\gl_{i_0}$. We know that $(u,\gl_{i_0})\in \max(I)$. So we have $v>u>w,\gl-v>\gl_{i_0}-u>\gm-w\Ra v-1>w,(\gl-2)-(v-1)=\gl-v-1>\gm-w$. Hence $(v-1,\gl-2),(w,\gm)$ are not comparable with each other. Therefore $K$ is an antichain.
	
	Given an ideal $I\in \mcl{M}$, we obtain $J=\gc(I)$ directly as follows. 
	\equa{J=&\{(v-1,\gl-2)\mid (v,\gl)\in I,\gl>\gl_{i_0},v+1<\gl\}\\&\cup\{(v,\gl_{i_0}-1)\mid (v,\gl_{i_0})\in I,v\neq \gl_{i_0}-1\}\\&\cup\{(v,\gl)\mid (v,\gl)\in I,\gl<\gl_{i_0}\}.}
	For this direct definition of $J=\gc(I)$, the set $\max(J)$ is precisely given as in the statement of the theorem.
	This is because of the following. Let $\max([I]_{(\gl_i)})=\{(a_i,\gl_i)\}$ if it is nonempty and define $a_i=\gl_i$ if $\max([I]_{(\gl_i)})$ is empty for $1\leq i\leq k$. We have $a_{i_0}=u$ and it follows that since $(u,\gl_{i_0})\in \max(I)$, $\max([I]_{(\gl_{i_0}+1)})=\{(u+1,\gl_{i_0}+1)\},\max([I]_{(\gl_{i_0}-1)})=\{(u,\gl_{i_0}-1)\}$ and $a_{i_0-1}-1\geq a_{i_0}= u,\gl_{i_0}-a_{i_0}-1\geq \gl_{i_0+1}-a_{i_0+1}$.
	Using this and from the definition of $J$, we observe that $\max([J]_{(\gl_i-2)})=\{(a_i-1,\gl_i-2)\}$ for all $\gl_i>\gl_{i_0}, \max([J]_{(\gl_i)})=\{(a_i,\gl_i)\}$ for all $\gl_i<\gl_{i_0}$ and $\max([J]_{(\gl_{i_0}-1)})=\{(a_{i_0},\gl_{i_0}-1)\}$.
	Since $I$ is an ideal in $\mcl{J}(P_{\ul{\gl}})$ we have 
	\equa{a_1\geq a_2\geq \ldots\geq a_{i_0-1}\geq &a_{i_0}\geq a_{i_0+1}\geq \ldots \geq a_{k-1}\geq a_k\\
		\gl_1-a_1\geq \ldots\geq \gl_{i_0-1}-a_{i_0-1}\geq &\gl_{i_0}-a_{i_0}\geq \gl_{i_0+1}-a_{i_0+1}\geq \ldots \geq  \gl_k-a_k.}
	Now we observe that
	\fo{10}{10}{ 
		\equa{a_1-1\geq a_2-1\geq \ldots\geq a_{i_0-1}-1\geq &a_{i_0}\geq a_{i_0+1}\geq \ldots \geq a_{k-1}\geq a_k\\
			\gl_1-a_1-1\geq \ldots\geq \gl_{i_0-1}-a_{i_0-1}-1\geq &\gl_{i_0}-a_{i_0}-1\geq \gl_{i_0+1}-a_{i_0+1}\geq \ldots \geq \gl_k-a_k.}
	}
	So $J$ is indeed an ideal in $\mcl{J}(P_{\ul{\gn}})$. Now it is easy to check that $\max(J)$ is precisely given as in the statement of the theorem. Moreover we observe that $(a_{i_0},\gl_{i_0}-1)=(u,\gl_{i_0}-1)\in \max(J)$ if and only if $a_{i_0-1}-1>a_{i_0}=u,\gl_{i_0}-a_{i_0}-1>\gl_{i_0+1}-a_{i_0+1}$.
	
	Now the bijection $\gc:\mcl{M}\lra \mcl{J}(P_{\ul{\gn}})$ is clear as we can define $J$ from $I$ and define $I$ from $J$. It is also clear that $\gc$ is a lattice isomorphism. This proves the theorem.
\end{proof}
\begin{example}
	We give three examples here to Theorem~\ref{theorem:LatticeIso}(2).
	Let $\ul{\gl}=(\gl^1_1=4^1>\gl^1_2=\gl^1_k=2^1)$. Then $\ul{\gn}=((\gl_1-2)^1=2^1>(\gl_k-1)^1=1^1)$. Then the bijection is given as \equa{\max(I),I\in \mcl{J}(P_{\ul{\gl}}),\max(I)\cap P_{(\gl_k)}\neq\es &\os{\gc}{\llra} \max(J),J\in \mcl{J}(P_{\ul{\gn}})\\
		\{(1,4),(0,2)\} &\llra \{(0,2)\}\\
		\{(0,2)\} &\llra \{(0,1)\}\\
		\{(2,4),(1,2)\}&\llra \{(1,2)\}\\
		\{(1,2)\}&\llra \es.}
	Let $\ul{\gl}=(\gl^1_1=5^1>\gl^1_2=\gl^1_k=2^1)$. Then $\ul{\gn}=((\gl_1-2)^1=3^1>(\gl_k-1)^1=1^1)$. Then the bijection is 
	given as \equa{\max(I),I\in \mcl{J}(P_{\ul{\gl}}),\max(I)\cap P_{(\gl_k)}\neq\es &\os{\gc}{\llra} \max(J),J\in \mcl{J}(P_{\ul{\gn}})\\
		\{(1,5),(0,2)\} &\llra \{(0,3)\}\\
		\{(2,5),(0,2)\} &\llra \{(1,3),(0,1)\}\\
		\{(0,2)\}&\llra \{(0,1)\}\\
		\{(2,5),(1,2)\}&\llra \{(1,3)\}\\
		\{(3,5),(1,2)\}&\llra \{(2,3)\}\\
		\{(1,2)\}&\llra \es.}
	Let $\ul{\gl}=(\gl^1_1=3^1>\gl^1_2=\gl^1_k=2^1)$. Then $\ul{\gn}=((\gl_1-2)^1=1^1\geq (\gl_k-1)^1=1^1)=(1^2)$. Then the bijection is 
	given as \equa{\max(I),I\in \mcl{J}(P_{\ul{\gl}}),\max(I)\cap P_{(\gl_k)}\neq\es &\os{\gc}{\llra} \max(J),J\in \mcl{J}(P_{\ul{\gn}})\\
		\{(0,2)\}&\llra \{(0,1)\}\\
		\{(1,2)\}&\llra \es.}
	
\end{example}
\begin{theorem}
\label{theorem:LambdaMuOne}
Let $\ugl\in \Gl_0$ be such that $\gr_{i_0}= 2$ for some $1\leq i_0\leq k$. Let $\ul{\gm} =\big(\gl_1^{\gr_1}>\gl_2^{\gr_2}>\ldots>\gl_{i_0-1}^{\gr_{i_0-1}}>\gl_{i_0}^{\gr_{i_0}-1}=\gl_{i_0}^1>\gl_{i_0+1}^{\gr_{i_0+1}}>\ldots>\gl_{k-1}^{\gr_{k-1}}>\gl_k^{\gr_k}\big)$. Let $I\in \mcl{J}(P_{\ul{\gl}}),J\in \mcl{J}(P_{\ul{\gl'/I}}),K\in \mcl{J}(P_{\ul{\gl}''})$.
Define $X^{\ul{\gl}}_{I,J,K}=\{(x',x'')\in \grpp\mid I(\ol{x}')=J,I(x'')=K\}$  and
$\ga^{\ul{\gl}}_{I,J,K}=\mid(\grppp)_{J\cup K}\oplus\big((\grpppp)^*_K+(\grpppp)_J\big)\mid$. Let $I_1\in \mcl{J}(P_{\ul{\gm}}),J_1\in \mcl{J}(P_{\ul{\gm'/I_1}}),K_1\in \mcl{J}(P_{\ul{\gm}''})$.
Define $X^{\ul{\gm}}_{I_1,J_1,K_1}=\{(x',x'')\in \grmpp\mid I(\ol{x}')=J_1,I(x'')=K_1\}$  and
$\ga^{\ul{\gm}}_{I_1,J_1,K_1}=\mid(\grmppp)_{J_1\cup K_1}\oplus\big((\grmpppp)^*_{K_1}+(\grmpppp)_{J_1}\big)\mid$.

Then 
\equan{TwoPointFive}{\us{I\in \mcl{J}(P_{\ul{\gl}}),\max(I)\cap P_{(\gl_{i_0})}=\es }{\sum}&\bigg(\us{J\in \mcl{J}(P_{\ul{\gl^{'}/I}}),K\in \mcl{J}(P_{\ul{\gl}''})}{\sum}\frac{\mid X^{\ul{\gl}}_{I,J,K}\mid}{\ga^{\ul{\gl}}_{I,J,K}}\bigg)=\\ &\us{I_1\in \mcl{J}(P_{\ul{\gm}}), \max(I_1)\cap P_{(\gl_{i_0})}=\es }{\sum}\bigg(\us{J_1\in \mcl{J}(P_{\ul{\gm^{'}/I_1}}),K_1\in \mcl{J}(P_{\ul{\gm}''})}{\sum}\frac{\mid X^{\ul{\gm}}_{I_1,J_1,K_1}\mid}{\ga^{\ul{\gm}}_{I_1,J_1,K_1}}\bigg).}
\end{theorem}

\begin{theorem}
	\label{theorem:LambdaMuTwo}
	Let $\ugl\in \Gl_0$ be such that $\gr_{i_0}= 2$ for some $1\leq i_0\leq k$. Let $\ul{\gm} =\big(\gl_1^{\gr_1}>\gl_2^{\gr_2}>\ldots>\gl_{i_0-1}^{\gr_{i_0-1}}>\gl_{i_0}^{\gr_{i_0}-1}=\gl_{i_0}^1>\gl_{i_0+1}^{\gr_{i_0+1}}>\ldots>\gl_{k-1}^{\gr_{k-1}}>\gl_k^{\gr_k}\big)$. Let $I\in \mcl{J}(P_{\ul{\gl}}),J\in \mcl{J}(P_{\ul{\gl'/I}}),K\in \mcl{J}(P_{\ul{\gl}''})$.
	Define $X^{\ul{\gl}}_{I,J,K}=\{(x',x'')\in \grpp\mid I(\ol{x}')=J,I(x'')=K\}$  and
	$\ga^{\ul{\gl}}_{I,J,K}=\mid(\grppp)_{J\cup K}\oplus\big((\grpppp)^*_K+(\grpppp)_J\big)\mid$. Let $I_1\in \mcl{J}(P_{\ul{\gm}}),J_1\in \mcl{J}(P_{\ul{\gm'/I_1}}),K_1\in \mcl{J}(P_{\ul{\gm}''})$.
	Define $X^{\ul{\gm}}_{I_1,J_1,K_1}=\{(x',x'')\in \grmpp\mid I(\ol{x}')=J_1,I(x'')=K_1\}$  and
	$\ga^{\ul{\gm}}_{I_1,J_1,K_1}=\mid(\grmppp)_{J_1\cup K_1}\oplus\big((\grmpppp)^*_{K_1}+(\grmpppp)_{J_1}\big)\mid$.
	
	Then 
\equan{FourPointFive}{\us{I\in \mcl{J}(P_{\ul{\gl}}), \max(I)\cap P_{(\gl_{i_0})}\neq \es }{\sum}&\bigg(\us{\us{\us{\max([J]_{(\gl_{i_0})})<\max([K]_{(\gl_{i_0})}),\max(K)\cap P_{(\gl_{i_0})}=\es}{\max([J]_{(\gl_{i_0})})\geq \max([K]_{(\gl_{i_0})}) \text{ or }}}{J\in \mcl{J}(P_{\ul{\gl^{'}/I}}),K\in \mcl{J}(P_{\ul{\gl}''})}}{\sum}\frac{\mid X^{\ul{\gl}}_{I,J,K}\mid}{\ga^{\ul{\gl}}_{I,J,K}}\bigg)\\ &=\us{\os{I_1\in \mcl{J}(P_{\ul{\gm}})}{\max(I_1)\cap P_{(\gl_{i_0})}\neq \es} }{\sum}\bigg(\us{J_1\in \mcl{J}(P_{\ul{\gm^{'}/I_1}}),K_1\in \mcl{J}(P_{\ul{\gm}''})}{\sum}\frac{\mid X^{\ul{\gm}}_{I_1,J_1,K_1}\mid}{\ga^{\ul{\gm}}_{I_1,J_1,K_1}}\bigg).}
\end{theorem}
\begin{theorem}
	\label{theorem:LambdaNu}	
	Let $\ugl\in \Gl_0$ be such that $\gr_{i_0}= 2$ for some $1\leq i_0\leq k$ and $\ul{\gn}=\big((\gl_1-2)^{\gr_1}>(\gl_2-2)^{\gr_2}>\ldots>(\gl_{i_0-1}-2)^{\gr_{i_0-1}}\geq (\gl_{i_0}-1)^{\gr_{i_0}}\geq \gl_{i_0+1}^{\gr_{i_0+1}}>\ldots>\gl_{k-1}^{\gr_{k-1}}>\gl_k^{\gr_k}\big)$. Let $I\in \mcl{J}(P_{\ul{\gl}}),J\in \mcl{J}(P_{\ul{\gl'/I}}),K\in \mcl{J}(P_{\ul{\gl}''})$.
	Define $X^{\ul{\gl}}_{I,J,K}=\{(x',x'')\in \grpp\mid I(\ol{x}')=J,I(x'')=K\}$  and
	$\ga^{\ul{\gl}}_{I,J,K}=\mid(\grppp)_{J\cup K}\oplus\big((\grpppp)^*_K+(\grpppp)_J\big)\mid$. Let $I_1\in \mcl{J}(P_{\ul{\gn}}),J_1\in \mcl{J}(P_{\ul{\gn'/I_1}}),K_1\in \mcl{J}(P_{\ul{\gn}''})$.
	Define $X^{\ul{\gn}}_{I_1,J_1,K_1}=\{(x',x'')\in \grnpp\mid I(\ol{x}')=J_1,I(x'')=K_1\}$  and
	$\ga^{\ul{\gn}}_{I_1,J_1,K_1}=\mid(\grnppp)_{J_1\cup K_1}\oplus\big((\grnpppp)^*_{K_1}+(\grnpppp)_{J_1}\big)\mid$. Then we have 
	\equan{FiveFive}{\us{I\in \mcl{J}(P_{\ul{\gl}}), \max(I)\cap P_{(\gl_{i_0})}\neq \es }{\sum}&\bigg(\us{\us{\max([J]_{(\gl_{i_0})})<\max([K]_{(\gl_{i_0})}),\max(K)\cap P_{(\gl_{i_0})}\neq\es}{J\in \mcl{J}(P_{\ul{\gl^{'}/I}}),K\in \mcl{J}(P_{\ul{\gl}''})}}{\sum}\frac{\mid X^{\ul{\gl}}_{I,J,K}\mid}{\ga^{\ul{\gl}}_{I,J,K}}\bigg)\\&=\us{I_1\in \mcl{J}(P_{\ul{\gn}})}{\sum}\bigg(\us{J_1\in \mcl{J}(P_{\ul{\gn^{'}/I_1}}),K_1\in \mcl{J}(P_{\ul{\gn}''})}{\sum}\frac{\mid X^{\ul{\gn}}_{I_1,J_1,K_1}\mid}{\ga^{\ul{\gn}}_{I_1,J_1,K_1}}\bigg)}
\end{theorem}
We prove Theorem~\ref{theorem:RepeatedPartCase} first and then prove Theorems~\ref{theorem:LambdaMuOne},~\ref{theorem:LambdaMuTwo},~\ref{theorem:LambdaNu} sequentially.
\begin{proof}[Proof of Theorem~\ref{theorem:RepeatedPartCase}]
	First we observe that for an ideal $I\in\mcl{J}(P_{\ul{\gl}})=\mcl{J}(P_{\ul{\gm}})$ we have the following data. The partitions $\ul{\gl}',\ul{\gl}''$, the canonical form $e_I=(e'_I,0)$ with $e'_I\in \grppp$, the partition $\ul{\gl^{'}/I}$ associated to $\grpppI=\grppp\fs \R e_I'$, the stabilizer subgroup $(\autgp)_I\subseteq \autgp$ of the canonical form $e_I$.

Let $x=(x',x'')\in \grppp\oplus \grpppp=\grpp$. For $\ol{x}'\in \grpppI$, let $I(\ol{x}')=J\in \mcl{J}(P_{\ul{\gl^{'}/I}})$ and $K=I(x'')\in \mcl{J}(P_{\ul{\gl}''})$. Let $\ga^{\ul{\gl}}_{I,J,K}$ be the size of the $(\autgp)_I$ orbit of $x=(x',x'')$ which is the cardinality of the valuative set \equ{(\grppp)_{J\cup K}\oplus \big((\grpppp)^*_K+(\grpppp)_J\big)\subseteq \grppp\oplus \grpppp=\grpp.}
Let $X^{\ul{\gl}}_{I,J,K}=\{(x',x'')\in \grpp\mid I(\ol{x}')=J,I(x'')=K\}$. Then
the number of $(\autgp)_I$ orbits in $\grpp$ is a polynomial in $q$ with integer coefficients and is given by 
\equ{n_{\ul{\gl},I}(q)=\mid (\autgp)_I\bs \grpp\mid=\us{J\in \mcl{J}(P_{\ul{\gl^{'}/I}})}{\sum}\us{K\in \mcl{J}(P_{\ul{\gl}''})}{\sum}\frac{\mid X^{\ul{\gl}}_{I,J,K}\mid}{\ga^{\ul{\gl}}_{I,J,K}}.}
We also have 
\equ{n_{\ul{\gl}}(q)=\us{I\in \mcl{J}(P_{\ul{\gl}})}{\sum}\bigg(\us{J\in \mcl{J}(P_{\ul{\gl^{'}/I}}),K\in \mcl{J}(P_{\ul{\gl}''})}{\sum}\frac{\mid X^{\ul{\gl}}_{I,J,K}\mid}{\ga^{\ul{\gl}}_{I,J,K}}\bigg).}

Similarly we have the number of $(\autmgp)_I$ orbits in $\grmpp$ is a polynomial in $q$ with integer coefficients and is given by  
\equ{n_{\ul{\gm},I}(q)=\mid (\autmgp)_I\bs \grmpp\mid=\us{J\in \mcl{J}(P_{\ul{\gm^{'}/I}})}{\sum}\us{K\in \mcl{J}(P_{\ul{\gm}''})}{\sum}\frac{\mid X^{\ul{\gm}}_{I,J,K}\mid}{\ga^{\ul{\gm}}_{I,J,K}}.}
We also have 
\equ{n_{\ul{\gm}}(q)=\us{I\in \mcl{J}(P_{\ul{\gm}})}{\sum}\bigg(\us{J\in \mcl{J}(P_{\ul{\gm^{'}/I}}),K\in \mcl{J}(P_{\ul{\gm}''})}{\sum}\frac{\mid X^{\ul{\gm}}_{I,J,K}\mid}{\ga^{\ul{\gm}}_{I,J,K}}\bigg).}
For an ideal $I\in \mcl{J}(P_{\ul{\gn}})$, we have the number of $(\autngp)_I$ orbits in $\grnpp$ is a polynomial in $q$ with integer coefficients and is given by  
\equ{n_{\ul{\gn},I}(q)=\mid (\autngp)_I\bs \grnpp\mid=\us{J\in \mcl{J}(P_{\ul{\gn^{'}/I}})}{\sum}\us{K\in \mcl{J}(P_{\ul{\gn}''})}{\sum}\frac{\mid X^{\ul{\gn}}_{I,J,K}\mid}{\ga^{\ul{\gn}}_{I,J,K}}.}
We also have 
\equ{n_{\ul{\gn}}(q)=\us{I\in \mcl{J}(P_{\ul{\gn}})}{\sum}\bigg(\us{J\in \mcl{J}(P_{\ul{\gn^{'}/I}}),K\in \mcl{J}(P_{\ul{\gn}''})}{\sum}\frac{\mid X^{\ul{\gn}}_{I,J,K}\mid}{\ga^{\ul{\gn}}_{I,J,K}}\bigg).}

Theorem~\ref{theorem:RepeatedPartCase} follows if we establish the identity in Equation~\ref{Eq:TwoPointFive} and the following identity.

\equan{FourPointFivePlusFiveFive}{&\us{I\in \mcl{J}(P_{\ul{\gl}}), \max(I)\cap P_{(\gl_{i_0})}\neq \es }{\sum}\bigg(\us{J\in \mcl{J}(P_{\ul{\gl^{'}/I}}),K\in \mcl{J}(P_{\ul{\gl}''})}{\sum}\frac{\mid X^{\ul{\gl}}_{I,J,K}\mid}{\ga^{\ul{\gl}}_{I,J,K}}\bigg)=\\ &\us{\os{I\in \mcl{J}(P_{\ul{\gm}})}{\max(I)\cap P_{(\gl_{i_0})}\neq \es} }{\sum}\bigg(\us{J\in \mcl{J}(P_{\ul{\gm^{'}/I}}),K\in \mcl{J}(P_{\ul{\gm}''})}{\sum}\frac{\mid X^{\ul{\gm}}_{I,J,K}\mid}{\ga^{\ul{\gm}}_{I,J,K}}\bigg)+\us{I\in \mcl{J}(P_{\ul{\gn}})}{\sum}\bigg(\us{J\in \mcl{J}(P_{\ul{\gn^{'}/I}}),K\in \mcl{J}(P_{\ul{\gn}''})}{\sum}\frac{\mid X^{\ul{\gn}}_{I,J,K}\mid}{\ga^{\ul{\gn}}_{I,J,K}}\bigg).}

Theorem~\ref{theorem:LambdaMuOne} proves the identity in Equation~\ref{Eq:TwoPointFive}. The identity in Equation~\ref{Eq:FourPointFivePlusFiveFive} is obtained by summing the identities in Equations~\ref{Eq:FourPointFive},\ref{Eq:FiveFive}.
Theorem~\ref{theorem:LambdaMuTwo} proves the identity in Equation~\ref{Eq:FourPointFive} and Theorem~\ref{theorem:LambdaNu} proves the identity in Equation~\ref{Eq:FiveFive}. Hence Theorem~\ref{theorem:RepeatedPartCase} follows.
\end{proof}

\begin{proof}[Proof of Theorem~\ref{theorem:LambdaMuOne}]
We consider the case of ideals $I\in \mcl{J}(P_{\ul{\gl}})=\mcl{J}(P_{\ul{\gm}})$ such that $\max(I)\cap P_{(\gl_{i_0})}=\es$. For such an ideal $I$, we have $\ul{\gl}'=\ul{\gm}',\ul{\gl}''$ differs from $\ul{\gm}''$ by an extra part $\gl_{i_0}$, that is, the isotypic part of $\ul{\gl}''$ corresponding to $\gl_{i_0}$ is $\gl_{i_0}^{\gr_{i_0}}=\gl_{i_0}^2$ whereas the isotypic part of $\ul{\gm}''$ corresponding to $\gl_{i_0}$ is $\gl_{i_0}^{\gr_{i_0}-1}=\gl_{i_0}^1$. We have $\mcl{J}(P_{\ul{\gl}''})=\mcl{J}(P_{\ul{\gm}''})$. We also have $\ul{\gl^{'}/I}=\ul{\gm^{'}/I},\grpppI=\grppp/\R e_I'=\grmpp/\R e_I'=\grmpppI$.
\begin{center}
	*****************************************************************
\end{center}
Let $J\in \mcl{J}(P_{\ul{\gl^{'}/I}})=\mcl{J}(P_{\ul{\gm^{'}/I}})$ and $K\in \mcl{J}(P_{\ul{\gl}''})=\mcl{J}(P_{\ul{\gm}''})$ be such that 
\equ{\max([J]_{(\gl_{i_0})})\geq \max([K]_{(\gl_{i_0})})}
that is, either  $\max([K]_{(\gl_{i_0})}=\es$ or if $\max([K]_{(\gl_{i_0})})= \{(u,\gl_{i_0})\}$ and $\max([J]_{(\gl_{i_0})})= \{(v,\gl_{i_0})\}$ then $v\leq u$. Choose $u=\gl_{i_0}$ if $\max([K]_{(\gl_{i_0})})=\es$ and $v=\gl_{i_0}$ if $\max([J]_{(\gl_{i_0})})=\max([K]_{(\gl_{i_0})})=\es$ . Then we have 

\equ{(\grppp)_{J\cup K}=(\grmppp)_{J\cup K}, \mid (\grpppp)^*_K+(\grpppp)_J\mid = q^{\gl_{i_0}-v}\mid (\grmpppp)^*_K+(\grmpppp)_J \mid.}
Hence \equ{\ga^{\ul{\gl}}_{I,J,K}=(q^{\gl_{i_0}-v})\ga^{\ul{\gm}}_{I,J,K}.}

Now, if $\max([J]_{(\gl_{i_0})})= \max([K]_{(\gl_{i_0})})$, that is, $v=u$ and if $\max(K)\cap P_{(\gl_{i_0})}=\es$, that is, $(u,\gl_{i_0})\nin \max(K)$ then the isotypic component of $(\grpppp)^*_K$ corresponding to $\gl_{i_0}$ is $\gp^u(\R/\gp^{\gl_{i_0}}\R)^{\gr_{i_0}}$ and the isotypic component of $(\grmpppp)^*_K$ corresponding to $\gl_{i_0}$ is $\gp^u(\R/\gp^{\gl_{i_0}}\R)^{\gr_{i_0}-1}$. Hence we have 
\equ{\mid X^{\ul{\gl}}_{I,J,K}\mid =(q^{\gl_{i_0}-u})\mid X^{\ul{\gm}}_{I,J,K}\mid =(q^{\gl_{i_0}-v})\mid X^{\ul{\gm}}_{I,J,K}\mid.}
In this case, we therefore have 
\equ{\frac{\mid X^{\ul{\gl}}_{I,J,K}\mid}{\ga^{\ul{\gl}}_{I,J,K}}=\frac{\mid X^{\ul{\gm}}_{I,J,K}\mid}{\ga^{\ul{\gm}}_{I,J,K}}.}
\begin{center}
	*****************************************************************
\end{center} 
If $\max([J]_{(\gl_{i_0})})=\max([K]_{(\gl_{i_0})})$, that is, $v=u$ and if $\max(K)\cap P_{(\gl_{i_0})}\neq\es$, that is, $(u,\gl_{i_0})\in \max(K)$ then the isotypic component of $(\grpppp)^*_K$ corresponding to $\gl_{i_0}$ is $\gp^u(\R/\gp^{\gl_{i_0}}\R)^{\gr_{i_0}}-\gp^{u+1}(\R/\gp^{\gl_{i_0}}\R)^{\gr_{i_0}}$ and the isotypic component of $(\grmpppp)^*_K$ corresponding to $\gl_{i_0}$ is $\gp^u(\R/\gp^{\gl_{i_0}}\R)^{\gr_{i_0}-1}-\gp^{u+1}(\R/\gp^{\gl_{i_0}}\R)^{\gr_{i_0}-1}$. Let $K_1\in \mcl{J}(P_{\ul{\gl}''})$ be such that \equ{\max(K_1)=\max(K)\bs\{(u,\gl_{i_0})\}.} Here, we consider ideals $L\in \mcl{J}(P_{\ul{\gl}''})$ such that 
\equ{K_1\subseteq L\subseteq K.}
Now $(\grpppp)^*_{K_1}\subseteq (\grpppp)^*_{L}\subseteq (\grpppp)^*_{K}$ and $\max(K)=\max(K_1)\sqcup\{(u,\gl_{i_0})\}$ implies that $\max(K_1)\subseteq \max(L)$. 

Conversely, given any $L\in \mcl{J}(P_{\ul{\gl}''})$ such that $\max([J]_{(\gl_{i_0})})> \max([L]_{(\gl_{i_0})})$, there is exactly one $K\in \mcl{J}(P_{\ul{\gl}''})$ such that $\{(v,\gl_{i_0})\}=\max([J]_{(\gl_{i_0})})=\max([K]_{(\gl_{i_0})})=\{(u,\gl_{i_0})\}$, and $K_1\subseteq L\subsetneq K$ where $K_1\in \mcl{J}(P_{\ul{\gl}''})$ is such that $\max(K_1)=\max(K)\bs\{(u,\gl_{i_0})\}$. The ideals $K$ and $K_1$ are obtained from $L$ in a unique manner as follows. The set $\max(K_1)$ is obtained from $\max(L)$ by excluding those elements from $\max(L)$ which are comparable with $(u,\gl_{i_0})$. The set $\max(K)=\max(K_1)\sqcup \{(u,\gl_{i_0})\}$ is an antichain and the ideal $K$ is generated by this antichain. Clearly we have $K_1\subseteq L \subsetneq K$.

Now, for any such $L$ with $K_1\subseteq L\subseteq K$ we have $J\cup L=J\cup K$ in the big fundamental poset $P$ as defined in~\ref{Eq:FP}, because $(u,\gl_{i_0})\in J,\max(K_1)\subseteq \max(K)\cap \max(L)$. Hence 
\equ{(\grppp)_{J\cup L}=(\grppp)_{J\cup K}.} 
We also have in this scenario that, $\max(K)\bs [J]_{\ul{\gl}''}=\max(L)\bs [J]_{\ul{\gl}''}$ implying
\equ{\mid (\grpppp)^*_L+(\grpppp)_J \mid=\mid(\grpppp)^*_K+(\grpppp)_J\mid. }
Therefore 
\equ{\ga^{\ul{\gl}}_{I,J,L}=\ga^{\ul{\gl}}_{I,J,K}\text{ for all }K_1\subseteq L\subseteq K.}
Similarly we have 
\equ{(\grmppp)_{J\cup L}=(\grmppp)_{J\cup K},\mid (\grmpppp)^*_L+(\grmpppp)_J \mid=\mid(\grmpppp)^*_K+(\grmpppp)_J\mid.}
Hence 
\equ{\ga^{\ul{\gm}}_{I,J,L}=\ga^{\ul{\gm}}_{I,J,K}\text{ for all }K_1\subseteq L\subseteq K.}
Moreover we have 
\equan{OneOne}{\ga^{\ul{\gl}}_{I,J,L}=\ga^{\ul{\gl}}_{I,J,K}=(q^{\gl_{i_0}-u})\ga^{\ul{\gm}}_{I,J,K}=(q^{\gl_{i_0}-u})\ga^{\ul{\gm}}_{I,J,L}.}	
Now we observe that 
\equ{\us{K_1\subseteq L\subseteq K}{\sqcup} (\grpppp)^*_{L}=\text{Product of two sets }=T_1\times T_2}
and 
\equ{\us{K_1\subseteq L\subseteq K}{\sqcup} (\grmpppp)^*_L=\text{Product of two sets }=S_1\times S_2}
where $S_1=T_1$ and $T_2,S_2$ are products of sets of the form $\gp^j(\R/\gp^{\gl}\R)^{\gr}$ such that the isotypic component of $T_2$ corresponding to $\gl_{i_0}$ is $\gp^u(\R/\gp^{\gl_{i_0}}\R)^{\gr_{i_0}}$, the isotypic component of $S_2$ corresponding to $\gl_{i_0}$ is $\gp^u(\R/\gp^{\gl_{i_0}}\R)^{\gr_{i_0}-1}$ and the remaining isotypic components of $T_2$ matches with the remaining isotypic components of $S_2$. Hence we have $\mid T_2\mid=(q^{\gl_{i_0}-u})\mid S_2\mid$.

So 
\equan{TwoTwo}{\mid \us{K_1\subseteq L\subseteq K}{\sqcup} X^{\ul{\gl}}_{I,J,L}\mid=(q^{\gl_{i_0}-u})\mid \us{K_1\subseteq L\subseteq K}{\sqcup} X^{\ul{\gm}}_{I,J,L}\mid.}

So we obtain by combining Equations~\ref{Eq:OneOne},~\ref{Eq:TwoTwo}, 
\equ{\us{K_1\subseteq L\subseteq K}{\sum}\frac{\mid X^{\ul{\gl}}_{I,J,L}\mid}{\ga^{\ul{\gl}}_{I,J,L}}=\frac{\mid \us{K_1\subseteq L\subseteq K}{\sqcup} X^{\ul{\gl}}_{I,J,L}\mid}{\ga^{\ul{\gl}}_{I,J,K}}=\frac{\mid \us{K_1\subseteq L\subseteq K}{\sqcup} X^{\ul{\gm}}_{I,J,L}\mid}{\ga^{\ul{\gm}}_{I,J,K}}=\us{K_1\subseteq L\subseteq K}{\sum}\frac{\mid X^{\ul{\gm}}_{I,J,L}\mid}{\ga^{\ul{\gm}}_{I,J,L}}.}
\begin{center}
	*****************************************************************
\end{center}
Now we consider the case  $J\in \mcl{J}(P_{\ul{\gl^{'}/I}})=\mcl{J}(P_{\ul{\gm^{'}/I}})$ and $K\in \mcl{J}(P_{\ul{\gl}''})=\mcl{J}(P_{\ul{\gm}''})$ such that 
\equ{\max([J]_{(\gl_{i_0})})< \max([K]_{(\gl_{i_0})}),}
that is, if $\max([J]_{(\gl_{i_0})})= \{(v,\gl_{i_0})\}$ and $\max([K]_{(\gl_{i_0})})= \{(u,\gl_{i_0})\}$ then $v>u$ or $\max([J]_{(\gl_{i_0})})=\es$ and $\max([K]_{(\gl_{i_0})})= \{(u,\gl_{i_0})\}\neq \es$. 
Here we have 
\equ{(\grppp)_{J\cup K}=(\grmppp)_{J\cup K}}
and if $(u,\gl_{i_0})\nin \max(K)$ then
\equ{\mid (\grpppp)^*_K+(\grpppp)_J \mid=q^{\gl_{i_0}-u}\mid(\grmpppp)^*_K+(\grmpppp)_J\mid, \mid (\grpppp)^*_K\mid = q^{\gl_{i_0}-u}\mid (\grmpppp)^*_K\mid.}
So we get that \equ{\ga^{\ul{\gl}}_{I,J,K}=(q^{\gl_{i_0}-u})\ga^{\ul{\gm}}_{I,J,K},\mid X^{\ul{\gl}}_{I,J,K}\mid =(q^{\gl_{i_0}-u})\mid X^{\ul{\gm}}_{I,J,K}\mid.}
Therefore 
\equ{\frac{\mid X^{\ul{\gl}}_{I,J,K}\mid}{\ga^{\ul{\gl}}_{I,J,K}}=\frac{\mid X^{\ul{\gm}}_{I,J,K}\mid}{\ga^{\ul{\gm}}_{I,J,K}}.}

Now, if $(u,\gl_{i_0})\in \max(K)$ then
\equa{\mid (\grpppp)^*_K+(\grpppp)_J \mid&=\frac{(q^{(\gl_{i_0}-u)\gr_{i_0}}-q^{(\gl_{i_0}-u-1)\gr_{i_0}})}{(q^{(\gl_{i_0}-u)(\gr_{i_0}-1)}-q^{(\gl_{i_0}-u-1)(\gr_{i_0}-1)})}\mid(\grmpppp)^*_K+(\grmpppp)_J\mid,\\
	\mid (\grpppp)^*_K\mid &= \frac{(q^{(\gl_{i_0}-u)\gr_{i_0}}-q^{(\gl_{i_0}-u-1)\gr_{i_0}})}{(q^{(\gl_{i_0}-u)(\gr_{i_0}-1)}-q^{(\gl_{i_0}-u-1)(\gr_{i_0}-1)})}\mid (\grmpppp)^*_K\mid.}
So we get that \equa{\ga^{\ul{\gl}}_{I,J,K}&=\frac{(q^{(\gl_{i_0}-u)\gr_{i_0}}-q^{(\gl_{i_0}-u-1)\gr_{i_0}})}{(q^{(\gl_{i_0}-u)(\gr_{i_0}-1)}-q^{(\gl_{i_0}-u-1)(\gr_{i_0}-1)})}\ga^{\ul{\gm}}_{I,J,K},\\
	\mid X^{\ul{\gl}}_{I,J,K}\mid &=\frac{(q^{(\gl_{i_0}-u)\gr_{i_0}}-q^{(\gl_{i_0}-u-1)\gr_{i_0}})}{(q^{(\gl_{i_0}-u)(\gr_{i_0}-1)}-q^{(\gl_{i_0}-u-1)(\gr_{i_0}-1)})}\mid X^{\ul{\gm}}_{I,J,K}\mid.}
Therefore again
\equ{\frac{\mid X^{\ul{\gl}}_{I,J,K}\mid}{\ga^{\ul{\gl}}_{I,J,K}}=\frac{\mid X^{\ul{\gm}}_{I,J,K}\mid}{\ga^{\ul{\gm}}_{I,J,K}}.}

\begin{center}
	*****************************************************************
\end{center}
Hence we have, for a fixed $I\in \mcl{J}(P_{\ul{\gl}})=\mcl{J}(P_{\ul{\gm}})$ such that $\max(I)\cap P_{(\gl_{i_0})}=\es$, 
\equ{\bigg(\us{J\in \mcl{J}(P_{\ul{\gl^{'}/I}}),K\in \mcl{J}(P_{\ul{\gl}''})}{\sum}\frac{\mid X^{\ul{\gl}}_{I,J,K}\mid}{\ga^{\ul{\gl}}_{I,J,K}}\bigg)=\bigg(\us{J\in \mcl{J}(P_{\ul{\gm^{'}/I}}),K\in \mcl{J}(P_{\ul{\gm}''})}{\sum}\frac{\mid X^{\ul{\gm}}_{I,J,K}\mid}{\ga^{\ul{\gm}}_{I,J,K}}\bigg).}

Upon summation over $I\in \mcl{J}(P_{\ul{\gl}})=\mcl{J}(P_{\ul{\gm}})$ such that $\max(I)\cap P_{(\gl_{i_0})}=\es$, the identity in Equation~\ref{Eq:TwoPointFive} follows.
\end{proof}

\begin{proof}[Proof of Theorem~\ref{theorem:LambdaMuTwo}]
	Now we consider the case when $\max(I)\cap P_{(\gl_{i_0})}\neq \es, I\in \mcl{J}(P_{\ul{\gl}})$. We have to prove the identity in Equation~\ref{Eq:FourPointFive}.

\begin{center}
	*****************************************************************
\end{center}
For an ideal $I\in \mcl{J}(P_{\ul{\gl}})=\mcl{J}(P_{\ul{\gm}}),\max(I)\cap P_{(\gl_{i_0})}\neq \es$ we have $\ul{\gl}'=\ul{\gm}'$ with the isotypic part corresponding to $\gl_{i_0}$ being $\gl_{i_0}^1$. We also have $\ul{\gl}''$ differs from $\ul{\gm}''$ by a part $\gl_{i_0}$, that is, the isotypic part of $\ul{\gl}''$ corresponding to $\gl_{i_0}$ being $\gl_{i_0}^{\gr_{i_0}-1}=\gl_{i_0}^1$ and there is no isotypic part of $\ul{\gm}''$ corresponding to $\gl_{i_0}$. We also have $\ul{\gl^{'}/I}=\ul{\gm^{'}/I}$ and $\grmpppI=\grmppp\fs \R e_I'=\grppp\fs \R e_I'=\grpppI$.

\begin{center}
	*****************************************************************
\end{center}
First we observe that ideals $K\in \mcl{J}(P_{\ul{\gl}''})$ such that $\max(K)\cap P_{(\gl_{i_0})}=\es$ are in bijection with ideals $L\in \mcl{J}(P_{\ul{\gm}''})$ using Theorem~\ref{theorem:LatticeIso}(1). The bijection being 
\equ{K\lra L=K\cap P_{\ul{\gm}''}\text{ with } \max(K)=\max(L).}

Let $J\in \mcl{J}(P_{\ul{\gl^{'}/I}})=\mcl{J}(P_{\ul{\gm^{'}/I}}), K\in \mcl{J}(P_{\ul{\gl}''})$. We define two ideals now, that arise from $K$. Let $K_1\in \mcl{J}(P_{\ul{\gl}''})$ be such that $\max(K_1)=\max(K)\bs P_{(\gl_{i_0})}\subs P_{\ul{\gm}''}$. Let $K_2=K\cap P_{\ul{\gm}''}\in \mcl{J}(P_{\ul{\gm}''})$. If $\max(K)\cap P_{(\gl_{i_0})}=\es$ then $\max(K_1)=\max(K_2)=\max(K)$ and $K_1=K,K_1\cap P_{\ul{\gm}''}=K_2$. If 
$\max([K]_{(\gl_{i_0})})=\{(u,\gl_{i_0})\}$ and $\max(K)\cap P_{(\gl_{i_0})}\neq\es$ then we have $\max(K)=\max(K_1)\sqcup\{(u,\gl_{i_0})\}$ and $\max(K_1)\subseteq \max(K_2)$. In general $\max(K_1)$ need not be equal to $\max(K_2)$.

Suppose \equ{\es\neq\{(v,\gl_{i_0})\}=\max([J]_{(\gl_{i_0})})\geq \max([K]_{(\gl_{i_0})})=\{(u,\gl_{i_0})\}\neq \es,}
that is, $v\leq u$. Choose $v=\gl_{i_0}$ if $\max([J]_{(\gl_{i_0})})=\es$ and $u=\gl_{i_0}$ if $\max([K]_{(\gl_{i_0})})=\es$. Then $v\leq u$ continues to hold. Let $K'\in \mcl{J}(P_{\ul{\gl}''}), L'\in \mcl{J}(P_{\ul{\gm}''})$ be such that $K_1\subseteq K'\subseteq K$ and $K_1\cap P_{\ul{\gm}''}=[K_1]_{\ul{\gm}''}\subseteq L'\subseteq K_2$. Then we observe that the elements in the sets $\max(K')\bs \max(K_1), \max(L')\bs \max(K_1)$ are comparable and less than or equal to $(u,\gl_{i_0})$ if the sets are nonempty. Hence these elements are less than or equal to $(v,\gl_{i_0})$. As a consequence we have $\max(K')\bs \max(K_1)\subseteq [J]_{\ul{\gl}''},\max(L')\bs \max(K_1)\subseteq [J]_{\ul{\gm}''}$ implying $\max(K')\bs [J]_{\ul{\gl}''}=\max(K_1)\bs [J]_{\ul{\gl}''}=\max(K)\bs [J]_{\ul{\gl}''}$ and $\max(L')\bs [J]_{\ul{\gm}''}=\max(K_1)\bs [J]_{\ul{\gm}''}=\max(K_2)\bs [J]_{\ul{\gl}''}$. So

\equa{(\grpppp)^*_{K'}+(\grpppp)_J &= (\grpppp)^*_K+(\grpppp)_J=(\grpppp)^*_{K_1}+(\grpppp)_J,\\ (\grmpppp)^*_{L'}+(\grmpppp)_J &= (\grmpppp)^*_{K_2}+(\grmpppp)_J=(\grmpppp)^*_{K_1}+(\grmpppp)_J.}

We also have that  
\equa{\mid (\grpppp)^*_{K'}+(\grpppp)_J\mid&=\mid (\grpppp)^*_K+(\grpppp)_J\mid\\&= q^{\gl_{i_0}-v}\mid (\grmpppp)^*_{K_2}+(\grmpppp)_J \mid=q^{\gl_{i_0}-v}\mid (\grmpppp)^*_{L'}+(\grmpppp)_J \mid.}

We observe that $J\cup K, J\cup K_1,J\cup K_2, J\cup K', J\cup L'$ all generate
the same ideal in the big fundamental poset $P$ defined in Equation~\ref{Eq:FP}. Hence we have 
\equ{(\grppp)_{J\cup K'}=(\grppp)_{J\cup K}=(\grppp)_{J\cup K_1}=(\grmppp)_{J\cup K_1}=(\grmppp)_{J\cup K_2}= (\grmppp)_{J\cup L'}.}
So
\equan{ThreeThree}{\ga^{\ul{\gl}}_{I,J,K'}=\ga^{\ul{\gl}}_{I,J,K}=\ga^{\ul{\gl}}_{I,J,K_1}=(q^{\gl_{i_0}-v})\ga^{\ul{\gm}}_{I,J,K_1}=(q^{\gl_{i_0}-v})\ga^{\ul{\gm}}_{I,J,K_2}=(q^{\gl_{i_0}-v})\ga^{\ul{\gm}}_{I,J,L'}.}
\begin{center}
	*****************************************************************
\end{center}

Now if $\max([J]_{(\gl_{i_0})})= \max([K]_{(\gl_{i_0})})$, that is, $v=u$ and if $\max(K)\cap P_{(\gl_{i_0})}=\es$, then the isotypic component of $(\grpppp)^*_K$ corresponding to $\gl_{i_0}$ is $\gp^u(\R/\gp^{\gl_{i_0}}\R)$ and $(\grpppp)^*_K=(\grmpppp)^*_{K_2}\oplus \gp^u(\R/\gp^{\gl_{i_0}}\R)$. Hence we have 
\equ{\mid X^{\ul{\gl}}_{I,J,K}\mid =(q^{\gl_{i_0}-u})\mid X^{\ul{\gm}}_{I,J,K_2}\mid =(q^{\gl_{i_0}-v})\mid X^{\ul{\gm}}_{I,J,K_2}\mid.}
In this case, we therefore have from Equation~\ref{Eq:ThreeThree} 
\equ{\frac{\mid X^{\ul{\gl}}_{I,J,K}\mid}{\ga^{\ul{\gl}}_{I,J,K}}=\frac{\mid X^{\ul{\gm}}_{I,J,K_2}\mid}{\ga^{\ul{\gm}}_{I,J,K_2}}.}
\begin{center}
	*****************************************************************
\end{center}
Now if $\max([J]_{(\gl_{i_0})})=\max([K]_{(\gl_{i_0})})$, that is, $v=u$ and if $\max(K)\cap P_{(\gl_{i_0})}\neq\es$, that is, $(u,\gl_{i_0})\in \max(K)$ then we consider ideals $K'\in \mcl{J}(P_{\ul{\gl}''})$ such that 
\equ{K_1\subseteq K'\subseteq K.}
Note that \equ{K_1\subseteq K'\subseteq K, K'\neq K \Llra K_1\subseteq K'\subseteq K, \max([J]_{(\gl_{i_0})})>\max([K']_{(\gl_{i_0})}).}
We also consider ideals $L'\in \mcl{J}(P_{\ul{\gm}''})$ such that 
\equ{K_1\cap P_{\ul{\gm}''}= [K_1]_{\ul{\gm}''}\subseteq L'\subseteq K_2=[K]_{\ul{\gm}''}=K\cap P_{\ul{\gm}''}.}
Note that for such ideals we have $\max([J]_{(\gl_{i_0})})>\max([L']_{(\gl_{i_0})})$. 

Given any $K'\in \mcl{J}(P_{\ul{\gl}''})$ such that $\max([J]_{(\gl_{i_0})})> \max([K']_{(\gl_{i_0})})$, there is exactly one $K\in \mcl{J}(P_{\ul{\gl}''})$ such that $\{(v,\gl_{i_0})\}=\max([J]_{(\gl_{i_0})})=\max([K]_{(\gl_{i_0})})=\{(u,\gl_{i_0})\}$, and $K_1\subseteq K'\subsetneq K$ where $K_1\in \mcl{J}(P_{\ul{\gl}''})$ is such that $\max(K_1)=\max(K)\bs\{(u,\gl_{i_0})\}$. The ideals $K$ and $K_1$ are obtained from $K'$ in a unique manner as follows. The set $\max(K_1)$ is obtained from $\max(K')$ by excluding those elements from $\max(K')$ which are comparable with $(u,\gl_{i_0})$. The set $\max(K)=\max(K_1)\sqcup \{(u,\gl_{i_0})\}$ is an antichain and the ideal $K$ is generated by this antichain. Clearly we have $K_1\subseteq K' \subsetneq K$ and $K,K_1$ are uniquely determined from $K'$.

Similarly given any $L'\in \mcl{J}(P_{\ul{\gm}''})$ such that $\max([J]_{(\gl_{i_0})})> \max([L']_{(\gl_{i_0})})$, there is exactly one $K\in \mcl{J}(P_{\ul{\gl}''})$ such that $\{(v,\gl_{i_0})\}=\max([J]_{(\gl_{i_0})})=\max([K]_{(\gl_{i_0})})=\{(u,\gl_{i_0})\}$ and $[K_1]_{\ul{\gm}''}\subseteq L'\subseteq K_2$ where $K_1\in \mcl{J}(P_{\ul{\gl}''})$ is such that $\max(K_1)=\max(K)\bs\{(u,\gl_{i_0})\}$ and $K_2=K\cap P_{\ul{\gm}''}$. The ideals $K$ and $K_1$ are obtained from $L'$ as follows. The set $\max(K_1)$ is obtained from $\max(L')$ by excluding those elements from $\max(L')$ which are comparable with $(u,\gl_{i_0})$. The set $\max(K)=\max(K_1)\sqcup \{(u,\gl_{i_0})\}$ is an antichain and the ideal $K$ is generated by this antichain. Clearly we have $[K_1]_{\ul{\gm}''}\subseteq L'\subseteq K_2=K\cap P_{\ul{\gm}''}$ and $K,K_1,K_2$ are uniquely determined from $L'$.

Now we observe that 
\equ{\us{K_1\subseteq K'\subseteq K}{\sqcup} (\grpppp)^*_{K'}=\text{Product of two sets }=T_1\times T_2}
and 
\equ{\us{[K_1]_{\ul{\gm}''}\subseteq L'\subseteq K_2}{\sqcup} (\grmpppp)^*_{L'}=\text{Product of two sets }=S_1\times S_2}
where $S_1=T_1$ and $T_2,S_2$ are products of sets of the form $\gp^j(\R/\gp^{\gl}\R)^{\gr}$ such that the isotypic component of $T_2$ corresponding to $\gl_{i_0}$ is $\gp^u(\R/\gp^{\gl_{i_0}}\R)$ and the other isotypic components of $T_2$ matches with the isotypic components of $S_2$. Hence we have $\mid T_2\mid=(q^{\gl_{i_0}-u})\mid S_2\mid=(q^{\gl_{i_0}-v})\mid S_2\mid$.

So 
\equan{FourFour}{\mid \us{K_1\subseteq K'\subseteq K}{\sqcup} X^{\ul{\gl}}_{I,J,K'}\mid=(q^{\gl_{i_0}-v})\mid \us{[K_1]_{\ul{\gm}''}\subseteq L'\subseteq K_2}{\sqcup} X^{\ul{\gm}}_{I,J,L'}\mid.}

So we obtain by combining Equations~\ref{Eq:ThreeThree},~\ref{Eq:FourFour}, 
\equ{\us{K_1\subseteq K'\subseteq K}{\sum}\frac{\mid X^{\ul{\gl}}_{I,J,K'}\mid}{\ga^{\ul{\gl}}_{I,J,K'}}=\frac{\mid \us{K_1\subseteq K'\subseteq K}{\sqcup} X^{\ul{\gl}}_{I,J,K'}\mid}{\ga^{\ul{\gl}}_{I,J,K}}=\frac{\mid \us{[K_1]_{\ul{\gm}''}\subseteq L'\subseteq K_2}{\sqcup} X^{\ul{\gm}}_{I,J,L'}\mid}{\ga^{\ul{\gm}}_{I,J,K_2}}=\us{[K_1]_{\ul{\gm}''}\subseteq L'\subseteq K_2}{\sum}\frac{\mid X^{\ul{\gm}}_{I,J,L'}\mid}{\ga^{\ul{\gm}}_{I,J,L'}}.}

\begin{center}
	*****************************************************************
\end{center} 

Suppose \equ{\es \text{ or }\{(v,\gl_{i_0})\}=\max([J]_{(\gl_{i_0})})< \max([K]_{(\gl_{i_0})})=\{(u,\gl_{i_0})\},\text{ that is, } v>u.} 
Choose $v=\gl_{i_0}$ if $\max([J]_{(\gl_{i_0})})=\es$.
Also let $K\in \mcl{J}(P_{\ul{\gl}''})$ be such that $\max(K)\cap P_{(\gl_{i_0})}=\es, L=K\cap P_{\ul{\gm}''}\in \mcl{J}(P_{\ul{\gm}''})$.
Then we have $\max([J]_{(\gl_{i_0})})< \max([L]_{(\gl_{i_0})})=\{(u,\gl_{i_0})\}$.
Also $J\cup K,J \cup L$ generate the same ideal in the big fundamental poset $P$ defined in Equation~\ref{Eq:FP}. Hence we have \equ{(\grppp)_{J\cup K}=(\grmppp)_{J\cup L}.}

Since $\max(K)\cap P_{(\gl_{i_0})}=\es$, that is, $(u,\gl_{i_0})\nin \max(K)$ we have $(\grpppp)^*_K=(\grmpppp)^*_L\oplus \gp^u(\R/\gp^{\gl_{i_0}}\R)$. So $\mid(\grpppp)^*_K\mid=(q^{\gl_{i_0}-u})\mid (\grmpppp)^*_L\mid$. Hence 
\equ{\mid X^{\ul{\gl}}_{I,J,K}\mid =q^{\gl_{i_0}-u}\mid X^{\ul{\gm}}_{I,J,L}\mid.} 
\equa{\mid\big((\grpppp)^*_K+(\grpppp)_J\big)\mid&=q^{\gl_{i_0}-u}\mid\big((\grmpppp)^*_L+(\grmpppp)_J\big)\mid\Ra\\& \ga^{\ul{\gl}}_{I,J,K}=q^{\gl_{i_0}-u}\ga^{\ul{\gm}}_{I,J,L}.}
Therefore 
\equ{\frac{\mid X^{\ul{\gl}}_{I,J,K}\mid}{\ga^{\ul{\gl}}_{I,J,K}}=\frac{\mid X^{\ul{\gm}}_{I,J,L}\mid}{\ga^{\ul{\gm}}_{I,J,L}}.}

So we have proved the identity in Equation~\ref{Eq:FourPointFive}. Hence Theorem~\ref{theorem:LambdaMuTwo} follows.
\end{proof}
Now we mention an example which supports Equation~\ref{Eq:FiveFive} in Theorem~\ref{theorem:LambdaNu}.
\begin{example}
	\label{example:RepeatedCase}
	Let $\ul{\gl}=(\gl_1=4^1>\gl^2_2=\gl^2_k=2^2), \ul{\gm}=(4^1>2^1),\ul{\gn}=((\gl_1-2)^1=2^1>(\gl_k-1)^2=1^2)$. We have $n_{\ul{\gl}}(q)=q^4 + 5q^3 + 13q^2 + 16q + 10$,
	$n_{\ul{\gm}}(q)=q^4 + 5q^3 + 12q^2 + 11q + 4$ and $n_{\ul{\gn}}(q)=q^2 + 5q + 6$.
	We have $\{I\in \mcl{J}(P_{\ul{\gl}})\mid \max(I)\cap P_{(\gl_k)}\neq \es\}=\{I^1=\langle\{(1,4),(0,2)\}\rangle >I^2=\langle\{(0,2)\}\rangle>I^3=\langle \{(2,4),(1,2)\}\rangle>I^4=\langle\{(1,2)\}\rangle\}$. We also have 
	$\mcl{J}(P_{\ul{\gn}})=\{I_1=\langle\{(0,2)\}\rangle>I_2=\langle\{(0,1)\}\rangle>I_3=\langle\{(1,2)\}\rangle>I_4=\es\}$.
	Then map $I^i\lra I_i$ from the set $\{I\in \mcl{J}(P_{\ul{\gl}})\mid \max(I)\cap P_{(\gl_k)}\neq \es\}$ to $\mcl{J}(P_{\ul{\gn}})$ is an isomorphism.
	\fo{9}{9}{
		\begin{center}
			\begin{tabular}{|c|c|c|c|c|c|c|}
				\hline\\
				Ideal & $J \in \mcl{J}(P_{\ul{\gl^{'}/I}})$ & $K\in \mcl{J}(P_{\ul{\gl}''})$ & $\max([J]_{(\gl_k)})$ & $\ga^{\ul{\gl}}_{I,J,K}$ & $\mid X^{\ul{\gl}}_{I,J,K}\mid$& Ratio\\ 
				&$\max(J)$ & $\max(K)\cap P_{(\gl_k)}\neq \es$ &$<\max([K]_{(\gl_k)})$ & & &$\frac{\mid X^{\ul{\gl}}_{I,J,K}\mid}{\ga^{\ul{\gl}}_{I,J,K}}$ \\
				\hline\\
				\multirow{5}{4.5em}{$I^1,\ul{\gl}'=(4>2)$, $\ul{\gl}''=(2)$, $\ul{\gl^{'}/I^1}=(3)$} & $\{(1,3)\}$ & $\{(0,2)\}$ &$(1,2)<(0,2)$&$q^5(q-1)$ &$q^5(q-1)^2$ &$(q-1)$\\ 
				& $\{(2,3)\}$ & $\{(0,2)\}$ & $\es<(0,2)$&$q^5(q-1)$ &$q^4(q-1)^2$ &$\frac{q-1}{q}$\\ 
				& $\es$ & $\{(0,2)\}$ & $\es<(0,2)$ & $q^5(q-1)$ & $q^4(q-1)$ & $\frac{1}{q}$\\ 
				&$\{(2,3)\}$& $\{(1,2)\}$ & $\es<(1,2)$ & $q^2(q-1)$ & $q^3(q-1)^2$ & $q(q-1)$\\ 
				& $\es$ & $\{(1,2)\}$ & $\es<(1,2)$ & $q^2(q-1)$ & $q^3(q-1)$ & $q$\\ 
				\hline
			\end{tabular}
		\end{center}
		\begin{center}
			\begin{tabular}{|c|c|c|c|c|c|}
				\hline\\
				Ideal & $J \in \mcl{J}(P_{\ul{\gn^{'}/I}})$ & $K\in \mcl{J}(P_{\ul{\gn}''})$ & $\ga^{\ul{\gn}}_{I,J,K}$ & $\mid X^{\ul{\gn}}_{I,J,K}\mid$& Ratio\\ 
				&$\max(J)$ & $\max(K)$ & & &$\frac{\mid X^{\ul{\gn}}_{I,J,K}\mid}{\ga^{\ul{\gn}}_{I,J,K}}$ \\
				\hline\\
				\multirow{2}{3cm}{$I_1,\ul{\gn}'=(2), \ul{\gn}''=(1,1),\ul{\gn^{'}/I_1}=\es$} & $\es$ & $\{(0,1)\}$ &$q(q^2-1)$ &$q^2(q^2-1)$ &$q$\\ 
				& $\es$ & $\es$ &$1$ &$q^2$ &$q^2$\\ 
				\hline
			\end{tabular}
		\end{center}
		\begin{center}
			\begin{tabular}{|c|c|c|c|c|c|c|}
				\hline\\
				Ideal & $J \in \mcl{J}(P_{\ul{\gl^{'}/I}})$ & $K\in \mcl{J}(P_{\ul{\gl}''})$ & $\max([J]_{(\gl_k)})$ & $\ga^{\ul{\gl}}_{I,J,K}$ & $\mid X^{\ul{\gl}}_{I,J,K}\mid$& Ratio\\ 
				&$\max(J)$ & $\max(K)\cap P_{(\gl_k)}\neq \es$ &$<\max([K]_{(\gl_k)})$ & & &$\frac{\mid X^{\ul{\gl}}_{I,J,K}\mid}{\ga^{\ul{\gl}}_{I,J,K}}$ \\
				\hline\\
				\multirow{4}{5em}{$I^2,\ul{\gl}'=(2)$, $\ul{\gl}''=(4>2)$, $\ul{\gl^{'}/I^2}=\es$} & $\es$ & $\{(1,4),(0,2)\}$ &$\es<(0,2)$&$q^5(q-1)^2$ &$q^5(q-1)^2$ &$1$\\ 
				& $\es$ & $\{(0,2)\}$ & $\es<(0,2)$&$q^5(q-1)$ &$q^5(q-1)$ &$1$\\ 
				& $\es$ & $\{(2,4),(1,2)\}$ & $\es<(1,2)$ & $q^2(q-1)^2$ & $q^3(q-1)^2$ & $q$\\ 
				&$\es$& $\{(1,2)\}$ & $\es<(1,2)$ & $q^2(q-1)$ & $q^3(q-1)$ & $q$\\ 
				\hline
			\end{tabular}
		\end{center}
		\begin{center}
			\begin{tabular}{|c|c|c|c|c|c|}
				\hline\\
				Ideal & $J \in \mcl{J}(P_{\ul{\gn^{'}/I}})$ & $K\in \mcl{J}(P_{\ul{\gn}''})$ & $\ga^{\ul{\gn}}_{I,J,K}$ & $\mid X^{\ul{\gn}}_{I,J,K}\mid$& Ratio\\ 
				&$\max(J)$ & $\max(K)$ & & &$\frac{\mid X^{\ul{\gn}}_{I,J,K}\mid}{\ga^{\ul{\gn}}_{I,J,K}}$ \\
				\hline\\
				\multirow{4}{3cm}{$I_2,\ul{\gn}'=(1), \ul{\gn}''=(2,1),\ul{\gn^{'}/I_2}=\es$} & $\es$ & $\{(0,2)\}$ &$q^3(q-1)$ &$q^3(q-1)$ &$1$\\ 
				& $\es$ & $\{(0,1)\}$ &$q^2(q-1)$ &$q^2(q-1)$ &$1$\\ 
				& $\es$ & $\{(1,2)\}$ &$(q-1)$ &$q(q-1)$ &$q$\\ 
				& $\es$ & $\es$ &$1$ &$q$ &$q$\\ 
				\hline
			\end{tabular}
		\end{center}
		\begin{center}
			\begin{tabular}{|c|c|c|c|c|c|c|}
				\hline\\
				Ideal & $J \in \mcl{J}(P_{\ul{\gl^{'}/I}})$ & $K\in \mcl{J}(P_{\ul{\gl}''})$ & $\max([J]_{(\gl_k)})$ & $\ga^{\ul{\gl}}_{I,J,K}$ & $\mid X^{\ul{\gl}}_{I,J,K}\mid$& Ratio\\ 
				&$\max(J)$ & $\max(K)\cap P_{(\gl_k)}\neq \es$ &$<\max([K]_{(\gl_k)})$ & & &$\frac{\mid X^{\ul{\gl}}_{I,J,K}\mid}{\ga^{\ul{\gl}}_{I,J,K}}$ \\
				\hline\\
				\multirow{7}{5em}{$I^3,\ul{\gl}'=(4>2)$, $\ul{\gl}''=(2)$, $\ul{\gl^{'}/I^3}=(3>1)$} & $\{(1,3),(0,1)\}$ & $\{(0,2)\}$ &$(1,2)<(0,2)$&$q^5(q-1)$ &$q^4(q-1)^3$ &$\frac{(q-1)^2}{q}$\\ 
				& $\{(0,1)\}$ & $\{(0,2)\}$ &$(1,2)<(0,2)$&$q^5(q-1)$ &$q^4(q-1)^2$ &$\frac{q-1}{q}$\\
				& $\{(1,3)\}$ & $\{(0,2)\}$ & $(1,2)<(0,2)$&$q^5(q-1)$ &$q^4(q-1)^2$ &$\frac{q-1}{q}$\\ 
				& $\{(2,3)\}$ & $\{(0,2)\}$ & $\es<(0,2)$ & $q^5(q-1)$ & $q^3(q-1)^2$ & $\frac{q-1}{q^2}$\\ 
				&$\es$& $\{(0,2)\}$ & $\es<(0,2)$ & $q^5(q-1)$ & $q^3(q-1)$ & $\frac 1{q^2}$\\ 
				&$(2,3)$& $\{(1,2)\}$ & $\es<(1,2)$ & $q^2(q-1)$ & $q^2(q-1)^2$ & $q-1$\\ 
				&$\es$& $\{(1,2)\}$ & $\es<(1,2)$ & $q^2(q-1)$ & $q^2(q-1)$ & $1$\\ 
				\hline
			\end{tabular}
		\end{center}
		\begin{center}
			\begin{tabular}{|c|c|c|c|c|c|}
				\hline\\
				Ideal & $J \in \mcl{J}(P_{\ul{\gn^{'}/I}})$ & $K\in \mcl{J}(P_{\ul{\gn}''})$ & $\ga^{\ul{\gn}}_{I,J,K}$ & $\mid X^{\ul{\gn}}_{I,J,K}\mid$& Ratio\\ 
				&$\max(J)$ & $\max(K)$ & & &$\frac{\mid X^{\ul{\gn}}_{I,J,K}\mid}{\ga^{\ul{\gn}}_{I,J,K}}$ \\
				\hline\\
				\multirow{4}{3cm}{$I_3,\ul{\gn}'=(2), \ul{\gn}''=(1,1),\ul{\gn^{'}/I_3}=(1)$} & $\{(0,1)\}$ & $\{(0,1)\}$ &$q^3$ &$q(q-1)^2(q+1)$ &$q-1-\frac 1q+\frac 1{q^2}$\\ 
				& $\es$ & $\{(0,1)\}$ &$q(q-1)(q+1)$ &$q(q-1)(q+1)$ &$1$\\ 
				& $\{(0,1)\}$ & $\es$ &$q^3$ &$q(q-1)$ &$\frac{(q-1)}{q^2}$\\ 
				& $\es$ & $\es$ &$1$ &$q$ &$q$\\ 
				\hline
			\end{tabular}
		\end{center}
		\begin{center}
			\begin{tabular}{|c|c|c|c|c|c|c|}
				\hline\\
				Ideal & $J \in \mcl{J}(P_{\ul{\gl^{'}/I}})$ & $K\in \mcl{J}(P_{\ul{\gl}''})$ & $\max([J]_{(\gl_k)})$ & $\ga^{\ul{\gl}}_{I,J,K}$ & $\mid X^{\ul{\gl}}_{I,J,K}\mid$& Ratio\\ 
				&$\max(J)$ & $\max(K)\cap P_{(\gl_k)}\neq \es$ &$<\max([K]_{(\gl_k)})$ & & &$\frac{\mid X^{\ul{\gl}}_{I,J,K}\mid}{\ga^{\ul{\gl}}_{I,J,K}}$ \\
				\hline\\
				\multirow{6}{5em}{$I^4,\ul{\gl}'=(2)$, $\ul{\gl}''=(4>2)$, $\ul{\gl^{'}/I^4}=(1)$} & $\{(0,1)\}$ & $\{(1,4),(0,2)\}$ &$(1,2)<(0,2)$&$q^5(q-1)^2$ &$q^4(q-1)^3$ &$\frac{q-1}{q}$\\ 
				& $\es$ & $\{(1,4),(0,2)\}$ &$\es<(0,2)$&$q^5(q-1)^2$ &$q^4(q-1)^2$ &$\frac{1}{q}$\\
				& $\{(0,1)\}$ & $\{(0,2)\}$ & $(1,2)<(0,2)$&$q^5(q-1)$ &$q^4(q-1)^2$ &$\frac{q-1}{q}$\\ 
				& $\es$ & $\{(0,2)\}$ & $\es<(0,2)$ & $q^5(q-1)$ & $q^4(q-1)$ & $\frac1{q}$\\ 
				&$\es$& $\{(2,4),(1,2)\}$ & $\es<(1,2)$ & $q^2(q-1)^2$ & $q^2(q-1)^2$ & $1$\\ 
				&$\es$& $\{(1,2)\}$ & $\es<(1,2)$ & $q^2(q-1)$ & $q^2(q-1)$ & $1$\\ 
				\hline
			\end{tabular}
		\end{center}
		\begin{center}
			\begin{tabular}{|c|c|c|c|c|c|}
				\hline\\
				Ideal & $J \in \mcl{J}(P_{\ul{\gn^{'}/I}})$ & $K\in \mcl{J}(P_{\ul{\gn}''})$ & $\ga^{\ul{\gn}}_{I,J,K}$ & $\mid X^{\ul{\gn}}_{I,J,K}\mid$& Ratio\\ 
				&$\max(J)$ & $\max(K)$ & & &$\frac{\mid X^{\ul{\gn}}_{I,J,K}\mid}{\ga^{\ul{\gn}}_{I,J,K}}$ \\
				\hline\\
				\multirow{4}{3cm}{$I_4,\ul{\gn}'=\es, \ul{\gn}''=(2>1^2),\ul{\gn^{'}/I_4}=\es$} & $\es$ & $\{(0,2)\}$ &$q^3(q-1)$ &$q^3(q-1)$ &$1$\\ 
				& $\es$ & $\{(0,1)\}$ &$q(q-1)(q+1)$ &$q(q-1)(q+1)$ &$1$\\ 
				& $\es$ & $\{(1,2)\}$ &$(q-1)$ &$(q-1)$ &$1$\\ 
				& $\es$ & $\es$ &$1$ &$1$ &$1$\\ 
				\hline
			\end{tabular}
		\end{center}
	}
	We have the sum of ratios $\frac{\mid X^{\ul{\gl}}_{I,J,K}\mid}{\ga^{\ul{\gl}}_{I,J,K}}$ equals
	\equa{&(q-1)+\frac{q-1}{q}+\frac1q+q(q-1)+q+1+1+q+q+\\
		&\frac{(q-1)^2}{q}+\frac{q-1}{q}+\frac{q-1}{q}+\frac{q-1}{q^2}+\frac 1{q^2}+(q-1)+1+\frac{q-1}{q}+\frac 1q+\frac{q-1}{q}+\frac 1q+1+1\\
		&=q^2+5q+6=n_{\ul{\gn}}(q).}
	Hence Example~\ref{example:RepeatedCase} agrees with Equation~\ref{Eq:FiveFive}.
\end{example}
Following Propositions~\ref{prop:LambdaandNU},~\ref{prop:LowerLatticeHomoLamdaPPPtoNuModI},~\ref{prop:UpperLatticeHomoLamdaPPPtoNuModI},~\ref{prop:UpperLower},~\ref{prop:ReverseLatticeMap},~\ref{prop:SandS1},~\ref{prop:ExistenceofJK},~\ref{prop:JKJ1K1Multiple} are useful in the proof of Theorem~\ref{theorem:LambdaNu} especially in establishing the identity in  Equation~\ref{Eq:FiveFive}. They may appear slightly unmotivated and technical because they appear before the proof of Theorem~\ref{theorem:LambdaNu}. The proofs of these propositions use lattice theory especially the theory of lattice of characteristic submodules $\mcl{J}(P_{\ul{\gl}})$ of any finite $\R$-module $\grpp$.
\begin{prop}
	\label{prop:LambdaandNU}
	Let $\ugl\in \Gl_0$ and $\gl_{i_0}$ be a part of $\ul{\gl}$ for some $1\leq i_0\leq k$. Let  $\ul{\gn}=\big((\gl_1-2)^{\gr_1}>(\gl_2-2)^{\gr_2}>\ldots>(\gl_{i_0-1}-2)^{\gr_{i_0-1}}\geq (\gl_{i_0}-1)^{\gr_{i_0}}\geq \gl_{i_0+1}^{\gr_{i_0+1}}>\ldots>\gl_{k-1}^{\gr_{k-1}}>\gl_k^{\gr_k}\big)$.
	Let  $\mcl{M}=\{I\in \mcl{J}(P_{\ul{\gl}})\mid \max(I)\cap P_{(\gl_{i_0})}\neq\es\}\subseteq \mcl{J}(P_{\ul{\gl}})$.
	For $I\in \mcl{M},I_1=\gc(I)$ where $\gc$ is as defined in Theorem~\ref{theorem:LatticeIso}(2), 
	\begin{itemize}
		\item if $\gm-2$ is a part of $\ul{\gn'/I_1}$ and $\gm>\gl_{i_0}+1$ then $\gm$ is a part of $\ul{\gl'/I}$. Also if $\gm< \gl_{i_0}-1$ and $\gm$ is a part of $\ul{\gn'/I_1}$ then $\gm$ is a part of $\ul{\gl'/I}$. 
		
		\item If $\gl$ is a part of $\ul{\gl'/I}$ and $\gl>\gl_{i_0}+1$ then $\gl-2$ is a part of $\ul{\gn'/I_1}$. If $\gl$ is a part of $\ul{\gl'/I}$ and $\gl<\gl_{i_0}-1$ then $\gl$ is a part of $\ul{\gn'/I_1}$. If $\gl$ is a part of $\ul{\gl'/I}$ and $\gl=\gl_{i_0}+1$ then $\gl-2$ need not be a part of $\ul{\gn'/I_1}$. If $\gl$ is a part of $\ul{\gl'/I}$ and $\gl=\gl_{i_0}-1$ then $\gl$ need not be a part of $\ul{\gn'/I_1}$.
		\item If $\max(I_1)\cap P_{(\gl_0-1)}\neq\es$ then $\gl_{i_0}+1$ is not a part of $\ul{\gl'/I}$. Also $\gl_{i_0}-1$ is a part of $\ul{\gl'/I}$ if and only if $\gl_{i_0}-1$ is a part of $\ul{\gn'/I}$.
		\item If $\max(I_1)\cap P_{(\gl_0-1)}=\es$ then either $\gl_{i_0}+1$ is a part of $\ul{\gl'/I}$ or $\gl_{i_0}-1$ is a part of $\ul{\gl'/I}$. Both $\gl_{i_0}\pm 1$ are parts of $\ul{\gl'/I}$ if and only if $\gl_{i_0}-1$ is a part of $\ul{\gn'/I}$.
		\item (More obvious conclusion:) $\gl_{i_0}$ is a never a part of $\ul{\gl'/I}$ for any $I\in \mcl{M}$.
	\end{itemize}
\end{prop}
\begin{proof}
	Let $\max(I)=\{(v_1,t_1),(v_2,t_2),\cdots,(v_{i-1},t_{i-1}),(v_i,t_i)=(v_i,\gl_{i_0}),(v_{i+1},t_{i+1})$ $\cdots,(v_s,t_s)\mid t_1>t_2>\cdots>t_{i-1}>t_i=\gl_{i_0}>t_{i+1}>\cdots>t_s\}$. Then $t_{i-1}>\gl_{i_0}+1\geq s-i+2,v_{i-1}\geq v_i+1\geq s-i+1,\gl_{i_0}-1>t_{i+1}\geq s-i,v_{i+1}\geq s-i-1$ and 
	\fo{10}{10}{
		\equa{v_1>v_2>\cdots>v_{i-1}>&v_i>v_{i+1}>\cdots>v_{s-1}>v_s\geq 0\\
			t_1-v_1>t_2-v_2>\cdots>t_{i-1}-v_{i-1}>&t_i-v_i>t_{i+1}-v_{i+1}>\cdots>t_{s-1}-v_{s-1}>t_s-v_s\geq 1.}
	}
	\equan{One}{\ul{\gl}'&=(t_1>t_2>\cdots>t_{i-1}>t_i=\gl_{i_0}>t_{i+1}>\cdots>t_{s-1}>t_s),\\ \ul{\gl^{'}/I}
		&=\big((v_1+t_2-v_2)>(v_2+t_3-v_3)>\cdots(v_{i-2}+t_{i-1}-v_{i-1})\\&>(v_{i-1}+t_i-v_i=v_{i-1}+\gl_{i_0}-v_i)>(v_i+t_{i+1}-v_{i+1})>\cdots\\&>(v_{s-2}+t_{s-1}-v_{s-1})>(v_{s-1}+t_s-v_s)>v_s\big).} 
	
	\begin{center}
		*****************************************************************
	\end{center}
	
	In the first case, if $\max(I_1)=\{(v_1-1,t_1-2),(v_2-1,t_2-2),\cdots,(v_{i-1}-1,t_{i-1}-2),(v_i,t_i-1)=(v_i,\gl_{i_0}-1),(v_{i+1},t_{i+1}),\cdots,(v_{s-1},t_{s-1}),(v_s,t_s)\}$ then $(v_i,t_i-1)=(v_i,\gl_{i_0}-1)\in \max(I_1)$,
	\fo{9}{9}{
		\equa{v_1-1>v_2-1>\cdots>v_{i-1}-1>&v_i>v_{i+1}>\cdots>v_{s-1}>v_s\geq 0\\
			t_1-v_1-1>t_2-v_2-1>\cdots>t_{i-1}-v_{i-1}-1>&t_i-v_i-1>t_{i+1}-v_{i+1}>\cdots>t_{s-1}-v_{s-1}>t_s-v_s\geq 1.}
	}
	and
	\equan{Two}{\ul{\gn}'&=(t_1-2>t_2-2>\cdots>t_{i-1}-2>t_i-1=\gl_{i_0}-1>t_{i+1}>\cdots>t_{s-1}>t_s),\\
		\ul{\gn^{'}/I_1}&=\big((v_1+t_2-v_2-2)>(v_2+t_3-v_3-2)>\cdots>(v_{i-2}+t_{i-1}-v_{i-1}-2)\\&>(v_{i-1}+t_i-v_i-2=v_{i-1}+\gl_{i_0}-v_i-2)>(v_i+t_{i+1}-v_{i+1})\\&>(v_{i+1}+t_{i+2}-v_{i+2})\cdots>(v_{s-2}+t_{s-1}-v_{s-1})>(v_{s-1}+t_s-v_s)>v_s\big).} 
	
	The number of parts of $\ul{\gl}'$ is $s$. The number of parts of $\ul{\gn}'$ is also $s$. The number of parts of $\ul{\gl^{'}/I}$ is $s$ if $v_s>0$ and it is $s-1$ if $v_s=0$. Note that the part $v_{i-1}+t_i-v_i-2=v_{i-1}+\gl_{i_0}-v_i-2$ of $\ul{\gn^{'}/I_1}$ is greater than zero, since $v_{i-1}-v_i\geq 2$ and $\gl_{i_0}\geq 1$. The number of parts of $\ul{\gn^{'}/I_1}$ is $s$ if $v_s>0$, it is $s-1$ if $v_s=0$. It cannot be $s-2$ in this case because the case $i=s,v_s=0,v_{s-1}+t_s=2$, that is, $i_0=k,i=s,v_s=0,v_{s-1}=1,t_s=\gl_k=1$ is not allowed here since $i_0=k,i=s,t_s=\gl_k,(v_{s-1}-1,t_{s-1}-2),(v_s,t_s-1)\in \max(I_1)\Ra v_{s-1}-1>v_s$.

	\begin{center}
		*****************************************************************
	\end{center}
	
	In the second case, if  $\max(I_1)=\{(v_1-1,t_1-2),(v_2-1,t_2-2),\cdots,(v_{i-1}-1,t_{i-1}-2),(v_{i+1},t_{i+1}),\cdots,(v_{s-1},t_{s-1}),(v_s,t_s)\}$ then $(v_i,t_i-1)=(v_i,\gl_{i_0}-1)\nin \max(I_1)$,
	\fo{10}{10}{
		\equa{v_1-1>v_2-1>\cdots>v_{i-1}-1>&v_{i+1}>\cdots>v_{s-1}>v_s\geq 0\\
			t_1-v_1-1>t_2-v_2-1>\cdots>t_{i-1}-v_{i-1}-1>&t_{i+1}-v_{i+1}>\cdots>t_{s-1}-v_{s-1}>t_s-v_s\geq 1.}
	}
	\equan{Three}{\ul{\gn}'&=(t_1-2>t_2-2>\cdots>t_{i-1}-2>t_{i+1}>\cdots>t_{s-1}>t_s),\\
		\ul{\gn^{'}/I_1}&=\big((v_1+t_2-v_2-2)>(v_2+t_3-v_3-2)>\cdots>(v_{i-2}+t_{i-1}-v_{i-1}-2)\\&>(v_{i-1}+t_{i+1}-v_{i+1}-1)>(v_{i+1}+t_{i+2}-v_{i+2})>\cdots\\&>(v_{s-2}+t_{s-1}-v_{s-1})>(v_{s-1}+t_s-v_s)>v_s\big).} 
	
	The number of parts of $\ul{\gl}'$ is $s$. The number of parts of $\ul{\gn}'$ is $s-1$. The number of parts of $\ul{\gl^{'}/I}$ is $s$ if $v_s>0$ and it is $s-1$ if $v_s=0$. If $i=s$ (for example if $i_0=k,t_s=\gl_{i_0}=\gl_k$) then the number of parts of $\ul{\gn^{'}/I_1}$ is $s-1$ if $v_{s-1}>1$ and it is $s-2$ if $v_{s-1}=1$. If $i<s$ then the number of parts of $\ul{\gn^{'}/I_1}$ is $s-1$ if $v_s>0$ and it is $s-2$ if $v_s=0$.
	
	This case occurs if   
	\begin{enumerate}[label=(\alph*)]
		\item either $t_i-v_i-1=t_{i+1}-v_{i+1}$,
		\item or $v_{i-1}-1=v_i$.
	\end{enumerate}
	In case (a) we have 
	\equan{Four}{\ul{\gl^{'}/I}
		&=\big((v_1+t_2-v_2)>(v_2+t_3-v_3)>\cdots(v_{i-2}+t_{i-1}-v_{i-1})\\&>(v_{i-1}+t_i-v_i=v_{i-1}+\gl_{i_0}-v_i)>(v_i+t_{i+1}-v_{i+1}=t_i-1=\gl_{i_0}-1)\\&>(v_{i+1}+t_{i+2}-v_{i+2})>\cdots>(v_{s-2}+t_{s-1}-v_{s-1})>(v_{s-1}+t_s-v_s)>v_s\big),\\
		\ul{\gn^{'}/I_1}&=\big((v_1+t_2-v_2-2)>(v_2+t_3-v_3-2)>\cdots>(v_{i-2}+t_{i-1}-v_{i-1}-2)\\&>(v_{i-1}+t_{i+1}-v_{i+1}-1=v_{i-1}+t_i-v_i-2=v_{i-1}+\gl_{i_0}-v_i-2)\\&>(v_{i+1}+t_{i+2}-v_{i+2})>\cdots>(v_{s-2}+t_{s-1}-v_{s-1})>(v_{s-1}+t_s-v_s)>v_s\big).} 
	In case (b) we have
	\equan{Five}{\ul{\gl^{'}/I}
		&=\big((v_1+t_2-v_2)>(v_2+t_3-v_3)>\cdots(v_{i-2}+t_{i-1}-v_{i-1})\\&>(v_{i-1}+t_i-v_i=t_i+1=\gl_{i_0}+1)>(v_i+t_{i+1}-v_{i+1})>(v_{i+1}+t_{i+2}-v_{i+2})\\&>\cdots>(v_{s-2}+t_{s-1}-v_{s-1})>(v_{s-1}+t_s-v_s)>v_s\big),\\
		\ul{\gn^{'}/I_1}&=\big((v_1+t_2-v_2-2)>(v_2+t_3-v_3-2)>\cdots>(v_{i-2}+t_{i-1}-v_{i-1}-2)\\&>(v_{i-1}+t_{i+1}-v_{i+1}-1=v_i+t_{i+1}-v_{i+1})>(v_{i+1}+t_{i+2}-v_{i+2})>\cdots\\&>(v_{s-2}+t_{s-1}-v_{s-1})>(v_{s-1}+t_s-v_s)>v_s\big).} 
	\begin{center}
		*****************************************************************
	\end{center}
	From Equations~[\ref{Eq:One}-\ref{Eq:Five}], we observe that the proposition follows. 	
\end{proof}

\begin{prop}
	\label{prop:LowerLatticeHomoLamdaPPPtoNuModI}
	Let $\ugl\in \Gl_0$ and $\gl_{i_0}$ be a part of $\ul{\gl}$ for some $1\leq i_0\leq k$. Let  $\ul{\gn}=\big((\gl_1-2)^{\gr_1}>(\gl_2-2)^{\gr_2}>\ldots>(\gl_{i_0-1}-2)^{\gr_{i_0-1}}\geq (\gl_{i_0}-1)^{\gr_{i_0}}\geq \gl_{i_0+1}^{\gr_{i_0+1}}>\ldots>\gl_{k-1}^{\gr_{k-1}}>\gl_k^{\gr_k}\big)$.
	Let  $\mcl{M}=\{I\in \mcl{J}(P_{\ul{\gl}})\mid \max(I)\cap P_{(\gl_{i_0})}\neq\es\}$.
	For $I\in \mcl{M}\subseteq \mcl{J}(P_{\ul{\gl}})$ let $I_1=\gc(I)$ where $\gc$ is as defined in Theorem~\ref{theorem:LatticeIso}(2). Let $\ul{\gl}'''$ be the partition obtained by inserting a new part $\gl_{i_0}$ into the partition $\ul{\gl'/I}$. Let  $\mcl{N}=\{\ti{J}\in \mcl{J}(P_{\ul{\gl}'''})\mid \max(\ti{J})\cap P_{(\gl_{i_0})}\neq\es\}$. Then there exists a lattice homomorphism $\gf: \mcl{N}=\{\ti{J}\in \mcl{J}(P_{\ul{\gl}'''})\mid \max(\ti{J}) \cap P_{(\gl_{i_0})}\neq \es\}\lra \mcl{J}(P_{\ul{\gn'/I_1}})$ such that $\gf(\ti{J})=J_1$ where 
	\equa{\max(J_1)&=\{(v-1,\gl-2)\mid (v,\gl)\in \max(\ti{J}), \gl>\gl_{i_0}\}\\ &\cup\{(v,\gl)\mid (v,\gl)\in \max(\ti{J}),\gl<\gl_{i_0}\}.}  
\end{prop}
\begin{proof}
	If $\max(\ti{J}) \cap P_{(\gl_{i_0})}\neq \es$ then $(v,\gl)\in \max(\ti{J}), \gl>\gl_{i_0}\Ra \gl>\gl_{i_0}+1\Ra \gl-2$ is a part of $\ul{\gn'/I_1}$	using Proposition~\ref{prop:LambdaandNU}. We observe that if $\max(\ti{J}) \cap P_{(\gl_{i_0})}=\{(u,\gl_{i_0})\}$ and $(v,\gl)\in \max(\ti{J}), \gl>\gl_{i_0}$ then $v>u\Ra v\geq 1,\gl-v>\gl_{i_0}-u\geq 1\Ra \gl-v\geq 2\Ra v-1<\gl-2$. So $(v-1,\gl-2)\in P_{\ul{\gn'/I_1}}$.
	
	Also if $\max(\ti{J}) \cap P_{(\gl_{i_0})}\neq \es$ then $(v,\gl)\in \max(\ti{J}), \gl<\gl_{i_0}\Ra \gl<\gl_{i_0}-1\Ra \gl$ is a part of $\ul{\gn'/I_1}$	using Proposition~\ref{prop:LambdaandNU}. 
	
	Hence $\{(v-1,\gl-2)\mid (v,\gl)\in \max(\ti{J}), \gl>\gl_{i_0}\}\cup\{(v,\gl)\mid (v,\gl)\in \max(\ti{J}),\gl<\gl_{i_0}\}$ is a subset of $P_{\ul{\gn'/I_1}}$. Now we prove that this set is an antichain.
	Clearly the set $\{(v-1,\gl-2)\mid (v,\gl)\in \max(\ti{J}), \gl>\gl_{i_0}\}$ is an antichain and the set $\{(v,\gl)\mid (v,\gl)\in \max(\ti{J}),\gl<\gl_{i_0}\}$ is also an antichain. Now if $(v,\gl)\in \max(\ti{J}),\gl>\gl_{i_0}$ and $(w,\gm)\in \max(\ti{J}),\gm<\gl_{i_0}$ then we have $v>u>w,\gl-v>\gl_{i_0}-u>\gm-w$. This implies that $v-1>w,\gl-v-1>\gm-w$ which further implies that $(v-1,\gl-2),(w,\gm)$ are not comparable. Hence $\max(J_1)$ is an antichain and defines an ideal $J_1\in \mcl{J}(P_{\ul{\gn'/I_1}})$.
	
	Now it is easy to see that the map $\gf:\mcl{N}\lra \mcl{J}(P_{\ul{\gn'/I_1}})$ is a lattice homomorphism. This proves the proposition.
\end{proof}
In the following proposition we define another lattice homomorphism $\gx:\mcl{N}\lra \mcl{J}(P_{\ul{\gn'/I_1}})$.
\begin{prop}
	\label{prop:UpperLatticeHomoLamdaPPPtoNuModI}
	Let $\ugl\in \Gl_0$ and $\gl_{i_0}$ be a part of $\ul{\gl}$ for some $1\leq i_0\leq k$. Let  $\ul{\gn}=\big((\gl_1-2)^{\gr_1}>(\gl_2-2)^{\gr_2}>\ldots>(\gl_{i_0-1}-2)^{\gr_{i_0-1}}\geq (\gl_{i_0}-1)^{\gr_{i_0}}\geq \gl_{i_0+1}^{\gr_{i_0+1}}>\ldots>\gl_{k-1}^{\gr_{k-1}}>\gl_k^{\gr_k}\big)$.
	Let  $\mcl{M}=\{I\in \mcl{J}(P_{\ul{\gl}})\mid \max(I)\cap P_{(\gl_{i_0})}\neq\es\}$.
	For $I\in \mcl{M}\subseteq \mcl{J}(P_{\ul{\gl}})$ let $I_1=\gc(I)$ where $\gc$ is as defined in Theorem~\ref{theorem:LatticeIso}(2). Let $\ul{\gl}'''$ be the partition obtained by inserting a new part $\gl_{i_0}$ into the partition $\ul{\gl'/I}$. Let  $\mcl{N}=\{\ti{J}\in \mcl{J}(P_{\ul{\gl}'''})\mid \max(\ti{J})\cap P_{(\gl_{i_0})}\neq\es\}$. Then there exists a lattice homomorphism $\gx: \mcl{N}=\{\ti{J}\in \mcl{J}(P_{\ul{\gl}'''})\mid \max(\ti{J}) \cap P_{(\gl_{i_0})}\neq \es\}\lra \mcl{J}(P_{\ul{\gn'/I_1}})$ such that $\gx(\ti{J})=J_3$ where $\ti{J}\in \mcl{J}(P_{\ul{\gl}'''}),J_3 \in \mcl{J}(P_{\ul{\gn'/I_1}})$ are given as follows. Let $\max([\ti{J}]_{(\gm_i)})=\{(a_i,\gm_i)\}$ if it is nonempty, for $\gm_i$ a part of $\ul{\gl}'''$.
	Let $\max([\ti{J}]_{(\gm_i=\gl_{i_0})})=\{(a_i=u,\gm_i=\gl_{i_0})\}\in P_{(\gl_{i_0})}$. Then 
	\begin{enumerate}
		\item $\max([\ti{J}]_{(\gl_{i_0}+1)})=\{(u+1,\gl_{i_0}+1)\},\max([\ti{J}]_{(\gl_{i_0}-1)})=\{(u,\gl_{i_0}-1)\}$ if $u<\gl_{i_0}-1$ and $\max([\ti{J}]_{(\gl_{i_0}-1)})$ is empty if $u=\gl_{i_0}-1$.
		\item $\max([J_3])_{(\gm_i-2)}=\{(a_i-1,\gm_i-2)\}$ for $\gm_i> \gl_{i_0}+1$ provided $a_i-1<\gm_i-2$. Otherwise if $a_i=\gm_i-1$ then $\max([J_3])_{(\gm_i-2)}$ is empty. 
		\item $\max([J_3])_{(\gm_i)}=\max([\ti{J}]_{(\gm_i)})$ for $\gm_i<\gl_{i_0}-1$, that is, $\max([J_3])_{(\gm_i)}=\{(a_i,\gm_i)\}$ for $\gm_i< \gl_{i_0}-1$ and if $\max([\ti{J}]_{(\gm_i)})$ is nonempty, $\max([J_3])_{(\gm_i)}$ is empty for $\gm_i< \gl_{i_0}-1$ and if $\max([\ti{J}]_{(\gm_i)})$ is empty. 
		\item If $\gl_{i_0}-1$ is a part of $\ul{\gn'/I_1}$ then $\max([J_3])_{(\gl_{i_0}-1)}=\max([\ti{J}]_{(\gl_{i_0}-1)})$.	
	\end{enumerate}
	Let $\ul{\gl}''''$ be a partition obtained by inserting both $\gl_{i_0}\pm 1$ into $\ul{\gl}'''$ (if either of them is not a part of $\ul{\gl}'''$). Let $\ul{\gn}'''$ be a partition obtained by inserting $\gl_{i_0}-1$ into $\ul{\gn'/I_1}$ (if it is not a part of $\ul{\gn'/I_1}$). Let $\mcl{T}=\{J'\in \mcl{J}(P_{\ul{\gl}''''})\mid \max(J')\cap P_{(\gl_{i_0})}\neq \es\}$. Then we have that $\gx$ is the composition of the following maps.
	\equan{Comp1}{&\mcl{N}\us{\cong}{\lra} \mcl{T} \os{\gc}{\us{\cong}{\lra}} \mcl{J}(P_{\ul{\gn}'''})\lra \mcl{J}(P_{\ul{\gn'/I_1}})\\
		&\ti{J}\lra [\ti{J}]_{\ul{\gl}''''}\lra \gc([\ti{J}]_{\ul{\gl}''''})\lra [\gc([\ti{J}]_{\ul{\gl}''''})]_{\ul{\gn'/I_1}}}  
	where $\gc$ is the map defined in Theorem~\ref{theorem:LatticeIso}(2) corresponding to $\ul{\gl}''''$ and $\ul{\gn}'''$.
\end{prop}
\begin{proof}
	We prove (1). Since $(u,\gl_{i_0})\in \max(\ti{J})$, (1) follows.
	
	We prove (2),(3),(4). If $\gm_i$ is a part of $\ul{\gl'/I}$ then using Proposition~\ref{prop:LambdaandNU}, we observe that $\gm_i-2$ is a part of $\ul{\gn'/I_1}$ if $\gm_i>\gl_{i_0}+1$ and $\gm_i$ is a part of $\ul{\gn'/I_1}$ if $\gm_i<\gl_{i_0}-1$. If $\gm_i>\gl_{i_0}+1$ then $a_i>u\Ra a_i\geq 1, \gm_i-a_i>\gl_{i_0}-u\Ra \gm_i-a_i\geq 2\Ra 0\leq a_i-1\leq \gm_i-2$. So $(a_i-1,\gm_i-2)\in P_{(\gm_i-2)}$ if $a_i<\gm_i-1$. Now define $J_3$ such that $\max([J_3]_{(\gm)})$ is as given in (2),(3),(4) for $\gm$ a part of $\ul{\gn'/I_1}$. It is easy to check that this actually defines a well defined ideal in $\mcl{J}(P_{\ul{\gn'/I_1}})$. This proves (2),(3),(4).
	
	Now it is easy to see that the map $\gx:\mcl{N}\lra \mcl{J}(P_{\ul{\gn'/I_1}})$ is a lattice homomorphism and is the composition of the maps given in~\ref{Eq:Comp1}. Hence the proposition follows.
\end{proof}
\begin{prop}
	\label{prop:UpperLower}
	Let $\ugl\in \Gl_0$ and $\gl_{i_0}$ be a part of $\ul{\gl}$ for some $1\leq i_0\leq k$. Let  $\ul{\gn}=\big((\gl_1-2)^{\gr_1}>(\gl_2-2)^{\gr_2}>\ldots>(\gl_{i_0-1}-2)^{\gr_{i_0-1}}\geq (\gl_{i_0}-1)^{\gr_{i_0}}\geq \gl_{i_0+1}^{\gr_{i_0+1}}>\ldots>\gl_{k-1}^{\gr_{k-1}}>\gl_k^{\gr_k}\big)$.
	Let  $\mcl{M}=\{I\in \mcl{J}(P_{\ul{\gl}})\mid \max(I)\cap P_{(\gl_{i_0})}\neq\es\}$.
	For $I\in \mcl{M}\subseteq \mcl{J}(P_{\ul{\gl}})$ let $I_1=\gc(I)$ where $\gc$ is as defined in Theorem~\ref{theorem:LatticeIso}(2). Let $\ul{\gl}'''$ be the partition obtained by inserting a new part $\gl_{i_0}$ into the partition $\ul{\gl'/I}$. Let  $\mcl{N}=\{\ti{J}\in \mcl{J}(P_{\ul{\gl}'''})\mid \max(\ti{J})\cap P_{(\gl_{i_0})}\neq\es\}$. Let $\gf:\mcl{N}\lra \mcl{J}(P_{\ul{\gn'/I_1}}),\gx:\mcl{N}\lra \mcl{J}(P_{\ul{\gn'/I_1}})$ be the lattice homomorphisms as defined in Proposition~\ref{prop:LowerLatticeHomoLamdaPPPtoNuModI} and Proposition~\ref{prop:UpperLatticeHomoLamdaPPPtoNuModI} respectively. Then we have for any $\ti{J}\in \mcl{N}$ with $J_1=\gf(\ti{J}),J_3=\gx(\ti{J})$, \equ{\max(J_1)=\max(\gf(\ti{J}))\subseteq \max(\gx(\ti{J}))=\max(J_3).}
\end{prop}
\begin{proof}
	Let $\max([\ti{J}]_{(\gm_i)})=\{(a_i,\gm_i)\}$ if it is nonempty, for $\gm_i$ a part of $\ul{\gl}'''$.
	Let $\max([\ti{J}]_{(\gm_i=\gl_{i_0})})=\{(a_i=u,\gm_i=\gl_{i_0})\}\in P_{(\gl_{i_0})}$. Then both the sets $\max(J_1)$ and $\max(J_3)$ are contained in the set 
	\equa{U_1&=\{(a_i-1,\gm_i-2)\mid \gm_i>\gl_{i_0}+1, a_i<\gm_i-1\}\\&\cup\{(a_i,\gm_i)
		\mid \gm_i<\gl_{i_0}-1\}\cup \{(u,\gl_{i_0}-1)\}}
	if $u<\gl_{i_0}-1$ and both the sets $\max(J_1)$ and $\max(J_3)$ are contained in the set 
	\equ{U_2=\{(a_i-1,\gm_i-2)\mid \gm_i>\gl_{i_0}+1,a_i<\gm_i-1\}}
	if $u_0=\gl_{i_0}-1$. Moreover any element in the set $U_1$ or $U_2$ is less than or equal to some element in the set $\max(J_3)$. Now let $J\in \mcl{J}(P_{\ul{\gl'/I}})$ be such that $\max(J)=\max(\ti{J})\bs\{(u,\gl_{i_0})\}$. Now $(a,\gm)\in \max(J),\gm>\gl_{i_0}+1$ if and only if $(a-1,\gm-2)\in \max(J_1)$.  Also $(a,\gm)\in \max(J),\gm<\gl_{i_0}-1$ if and only if $(a,\gm)\in \max(J_1)$.	So in the initial case if $(a-1,\gm-2)\in \max(J_1)$ implies $(a-1,\gm-2)\in U_1$ or $U_2$ and $(a,\gm)\in \max(\ti{J})$ and in the latter case if $(a,\gm)\in \max(J_1)$ implies $(a,\gm)\in U_1$ and $(a,\gm)\in \max(\ti{J})$. 
	
	In the initial case for $(a-1,\gm-2)\in \max(J_1),\gm>\gl_{i_0}+1$ if $(a-1,\gm-2)\leq (b-1,\gl-2)\in U_2$ then $(a,\gm)\leq (b,\gl)\in \ti{J}\Ra (a,\gm)=(b,\gl)$ since $(a,\gm)\in \max(\ti{J})\Ra (a-1,\gm-2)=(b-1,\gl-2)$. If $\gm>\gl_{i_0}+1, \gl\leq \gl_{i_0}-1, (a-1,\gm-2)\in \max(J_1), (b,\gl)\in U_1$, then the case $(a-1,\gm-2)\leq (b,\gl)$ does not arise.
	This is because:
	\begin{enumerate}[label=Case (\alph*)]
		\item $(a,\gm)$ and $(u,\gl_{i_0})$ are not comparable as they are both in $\max(\ti{J})$, $(u,\gl_{i_0})$ and $(b,\gl)$ are not comparable in which 
		case we have $a>u>b,\gm-a>\gl_{i_0}-u>\gl-b\Ra $ the elements $(a-1,\gm-2)$ and $(b,\gl)$ are not comparable.
		\item  $(a,\gm)$ and $(u,\gl_{i_0})$ are not comparable as they are both in $\max(\ti{J})$, $(u,\gl_{i_0})\geq (b,\gl)$ and $u=b$ since $(b,\gl)\in U_1,\gl\leq \gl_{i_0}-1\Ra \max([\ti{J}]_{(\gl)})=\{(b,\gl)\}$. In this case we have $\gm-a>\gl_{i_0}-u=\gl_{i_0}-b>\gl-b,a>u=b$. So we have $a-1\geq b$ and $\gm-a-1> \gl-b \Ra$ either $(a-1,\gm-2),(b,\gl)$ are not comparable or $(a-1,\gm-2)$ and $(b,\gl)$ are comparable and $(b,\gl)< (a-1,\gm-2)$. 
	\end{enumerate}
	So we have $(a-1,\gm-2)\in \max(J_1),\gm>\gl_{i_0}+1 \Ra (a-1,\gm-2)$ is a maximal element in $U_1$ which implies that $(a-1,\gm-2)\in \max(J_3)$.
	
	In the latter case for $(a,\gm)\in \max(J_1),\gm< \gl_{i_0}-1$ we have $(a,\gm)\in \max(\ti{J})$. So $(a,\gm)\leq (b,\gl)\in U_1$ for some $\gl\leq \gl_{i_0}-1\Ra (b,\gl)\in [\ti{J}]_{(\gl)}$ and  $(a,\gm)=(b,\gl)$. If $\gl>\gl_{i_0}+1$ and $(b-1,\gl-2)\in U_1$ then the case $(a,\gm)\leq (b-1,\gl-2)$ does not arise. This is because:
	\begin{enumerate}[label=Case (\Alph*)]
		\item $(u,\gl_{i_0}),(a,\gm)$ are not comparable as they are both in $\max(\ti{J})$ and $(u,\gl_{i_0})$, $(b,\gl)$ are not comparable in which case we have $b>u>a, \gl-b>\gl_{i_0}-u>\gm-a\Ra $ the elements $(b-1,\gl-2),(a,\gm)$ are not comparable.  
		\item $(u,\gl_{i_0}),(a,\gm)$ are not comparable as they are both in $\max(\ti{J}), (u,\gl_{i_0})\geq (b,\gl)$ and $\gl-b=\gl_{i_0}-u$ since $(b-1,\gl-2)\in U_1,\gl>\gl_{i_0}+1\Ra\max([\ti{J}]_{(\gl)})=\{(b,\gl)\}$. In this case we have $\gl-b=\gl_{i_0}-u>\gm-a\Ra \gl-b-1\geq \gm-a$. Now $\gl-b=\gl_{i_0}-u$ and $\gl>\gl_{i_0}\Ra b>u$. So we have $b>u>a\Ra b-1>a$. So $\gl-b-1\geq \gm-a,b-1>a$ implies either $(b-1,\gl-2),(a,\gm)$ are not comparable or $(b-1,\gl-2),(a,\gm)$ are  comparable and $(b-1,\gl-2)<(a,\gm)$. 
	\end{enumerate} 
	So we have $(a,\gm)\in \max(J_1),\gm< \gl_{i_0}-1 \Ra (a,\gm)$ is a maximal element in $U_1$ which implies that $(a,\gm)\in \max(J_3)$.
	
	This proves the proposition that $\max(J_1)=\max(\gf(\ti{J}))\subseteq \max(\gx(\ti{J}))=\max(J_3)$.
\end{proof}
\begin{remark}
	The map $\gx:\mcl{N}\lra \mcl{J}(P_{\ul{\gn'/I_1}})$ is called the upper lattice homomorphism and $\gf:\mcl{N}\lra \mcl{J}(P_{\ul{\gn'/I_1}})$ is called the lower lattice homomorphism since for every $\ti{J}\in \mcl{N}$ we have $\max(\gf(\ti{J}))\subseteq \max(\gx(\ti{J}))$ implying $\gf(\ti{J})\subseteq \gx(\ti{J})$.
\end{remark}
\begin{remark}
	\label{remark:J1J3}
	For $(u,\gl_{i_0})\in \max(\ti{J}),\gf(\ti{J})=J_1,\gx(\ti{J})=J_3=[\gc([\ti{J}]_{\ul{\gl}''''})]_{\ul{\gn'/I_1}}$, 
	we have $[J_1\cup \langle \{(u,\gl_{i_0}-1)\}\rangle]_{\ul{\gn}'''}=\gc([\ti{J}]_{\ul{\gl}''''})$ if $u<\gl_{i_0}-1$ and if $u=\gl_{i_0}-1$ then $J_1=J_3$.
	This is because of the following. 
	
	In case $u<\gl_{i_0}-1$, if $(a,\gm)\in \max(J_3)\bs \max(J_1)$ and $\gm<\gl_{i_0}-1$ then $(a,\gm)\in \ti{J}$ and $(a,\gm)\leq (u,\gl_{i_0}) \Ra u\leq a, \gm-a<\gl_{i_0}-a\leq \gl_{i_0}-u\Ra \gm-a\leq \gl_{i_0}-u-1$. So $(a,\gm)\leq (u,\gl_{i_0}-1)$. 
	Also if $(a-1,\gm-2)\in \max(J_3)\bs \max(J_1)$ and $\gm>\gl_{i_0}+1$ then $(a,\gm)\in \ti{J}, (a,\gm)\leq (u,\gl_{i_0})$. So $u\leq a,\gm-a\leq \gl_{i_0}-u$. Now we can get better inequalities using the fact that $\{(u+1,\gl_{i_0}+1)\}= \max([\ti{J}]_{(\gl_{i_0}+1)}),\{(a,\gm)\}=\max([\ti{J}]_{(\gm)})$. Here we have $a\geq u+1,\gm-a\geq \gl_{i_0}-u$. So we get $\gm-a=\gl_{i_0}-u, u\leq a-1\Ra (a-1,\gm-2)\leq (u,\gl_{i_0}-1)$. This proves that $\max(J_3)\subs[J_1\cup \langle \{(u,\gl_{i_0}-1)\}\rangle]_{\ul{\gn}'''}\Ra [J_1\cup \langle \{(u,\gl_{i_0}-1)\}\rangle]_{\ul{\gn}'''}=\gc([\ti{J}]_{\ul{\gl}''''})$.
	
	In case $u=\gl_{i_0}-1$ then if $(a,\gm)\leq (u=\gl_{i_0}-1,\gl_{i_0})$ and $\gm\geq \gl_{i_0}+1$ then $a=\gm-1$. Also if $\gm<\gl_{i_0}$ then $\max([\ti{J}]_{(\gm)})=\es$. So $\max(J_1)=\max(J_3)\Ra J_1=J_3$. 
\end{remark}
Now we prove that there exists a lattice homomorphism $\gz:\mcl{J}(P_{\ul{\gn'/I_1}})\lra \mcl{N}$ in the reverse direction satisfying certain nice properties.
\begin{prop}
	\label{prop:ReverseLatticeMap}
	Let $\ugl\in \Gl_0$ and $\gl_{i_0}$ be a part of $\ul{\gl}$ for some $1\leq i_0\leq k$. Let  $\ul{\gn}=\big((\gl_1-2)^{\gr_1}>(\gl_2-2)^{\gr_2}>\ldots>(\gl_{i_0-1}-2)^{\gr_{i_0-1}}\geq (\gl_{i_0}-1)^{\gr_{i_0}}\geq \gl_{i_0+1}^{\gr_{i_0+1}}>\ldots>\gl_{k-1}^{\gr_{k-1}}>\gl_k^{\gr_k}\big)$.
	Let  $\mcl{M}=\{I\in \mcl{J}(P_{\ul{\gl}})\mid \max(I)\cap P_{(\gl_{i_0})}\neq\es\}$.
	For $I\in \mcl{M}\subseteq \mcl{J}(P_{\ul{\gl}})$ let $I_1=\gc(I)$ where $\gc$ is as defined in Theorem~\ref{theorem:LatticeIso}(2). Let $\ul{\gl}'''$ be the partition obtained by inserting a new part $\gl_{i_0}$ into the partition $\ul{\gl'/I}$. Let  $\mcl{N}=\{\ti{J}\in \mcl{J}(P_{\ul{\gl}'''})\mid \max(\ti{J})\cap P_{(\gl_{i_0})}\neq\es\}$. Then there exists a lattice homomorphism $\gz:\mcl{J}(P_{\ul{\gn'/I_1}})\lra \mcl{N}$ such that $\gz(J_3)=J_5$ where 
	$J_3\in \mcl{J}(P_{\ul{\gn'/I_1}})$ with (say) $\max([J_3]_{(\gd_i)})=\{(b_i,\gd_i)\}$ if nonempty, for $\gd_i$ a part of $\ul{\gn'/I_1}$ and (say) $\max([J_3]_{(\gl_{i_0}-1)})=\{(b,\gl_{i_0}-1)\}$ if nonempty, irrespective of whether $\gl_{i_0}-1$ is a part of $\ul{\gn'/I_1}$. Then $J_5$ is given as follows.
	\begin{enumerate}
		\item If $\gd_i+2>\gl_{i_0}+1$ and $\max([J_3]_{(\gd_i)})=\{(b_i,\gd_i)\}$ then $\max([J_5]_{(\gd_i+2)})=\{(b_i+1,\gd_i+2)\}$.
		\item If $\gd_i+2>\gl_{i_0}+1$ and $\max([J_3]_{(\gd_i)})$ is empty then $\max([J_5]_{(\gd_i+2)})=\{(\gd_i+1,\gd_i+2)\}$.
		\item If $\gd_i<\gl_{i_0}-1$ and $\max([J_5]_{(\gd_i)})=\max([J_3]_{(\gd_i)})$, that is, $\max([J_3]_{(\gd_i)})=\{(b_i,\gd_i)\}$ then $\max([J_5]_{(\gd_i)})=\{(b_i,\gd_i)\}$ and
		if $\max([J_3]_{(\gd_i)})$ is empty then $\max([J_5]_{(\gd_i)})$ is empty.
		\item If $\gl_{i_0}+1$ is a part of $\ul{\gl'/I}$ then 
		$\max([J_5]_{(\gl_{i_0}+1)})=\{(b+1,\gl_{i_0}+1)\}$ if $\max([J_3]_{(\gl_{i_0}-1)})=\{(b,\gl_{i_0}-1)\}$.
		If $\max([J_3]_{(\gl_{i_0}-1)})$ is empty  
		then $\max([J_5]_{(\gl_{i_0}+1)})$ $=\{(\gl_{i_0},\gl_{i_0}+1)\}$
		\item If $\gl_{i_0}-1$ is a part of $\ul{\gl'/I}$ then $\max([J_5]_{(\gl_{i_0}-1)})=\max([J_3]_{(\gl_{i_0}-1)})$.
		\item $\max([J_5]_{(\gl_{i_0})})=\{(b,\gl_{i_0})\}$ if $\max([J_3]_{(\gl_{i_0}-1)})=\{(b,\gl_{i_0}-1)\}$. 
		If $\max([J_3]_{(\gl_{i_0}-1)})$ is empty  
		then $\max([J_5]_{(\gl_{i_0})})=\{(\gl_{i_0}-1,\gl_{i_0})\}$.
	\end{enumerate} 
	Let $\ul{\gl}''''$ be a partition obtained by inserting both $\gl_{i_0}\pm 1$ into $\ul{\gl}'''$ (if either of them is not a part of $\ul{\gl}'''$). Let $\ul{\gn}'''$ be a partition obtained by inserting $\gl_{i_0}-1$ into $\ul{\gn'/I_1}$ (if it is not a part of $\ul{\gn'/I_1}$). Let $\mcl{T}=\{J'\in \mcl{J}(P_{\ul{\gl}''''})\mid \max(J')\cap P_{(\gl_{i_0})}\neq \es\}$. Then we have that $\gz$ is a composition of the following maps.
	\equan{Comp2}{&\mcl{J}(P_{\ul{\gn'/I_1}})\lra  \mcl{J}(P_{\ul{\gn}'''}) \os{\gc^{-1}}{\us{\cong}{\lra}}\mcl{T}\us{\cong}{\lra} \mcl{N}  \\
		&J_3\lra [J_3]_{\ul{\gn}'''}\lra \gc^{-1}([J_3]_{\ul{\gn}'''})\lra [\gc^{-1}([J_3]_{\ul{\gn}'''})]_{\ul{\gl}'''}} 
	where $\gc$ is the map defined in Theorem~\ref{theorem:LatticeIso}(2) corresponding to $\ul{\gl}''''$ and $\ul{\gn}'''$.
\end{prop}
\begin{proof}
	We observe that $J_5$ is an ideal in $\mcl{J}(P_{\ul{\gl}'''})$.	
	If $\max([J_3]_{(\gl_{i_0}-1)})=\{(b,\gl_{i_0}-1)\}$ then $(b,\gl_{i_0})\in \max(J_5)$. If $\max([J_3]_{(\gl_{i_0}-1)})$ is empty then $(\gl_{i_0}-1,\gl_{i_0})\in \max(J_5)$. So we have $\max(J_5)\cap P_{(\gl_{i_0})}\neq \es\Ra J_5\in \mcl{N}$.
	
	It is clear that the map $\gz:\mcl{J}(P_{\ul{\gn'/I_1}})\lra \mcl{N}$ is a lattice homomorphism and is a composition of maps given in~\ref{Eq:Comp2}. This proves the proposition. 
\end{proof}
\begin{remark}
	If $\gl_{i_0}-1$ is a part of $\ul{\gn'/I_1}$, that is, $\ul{\gn}'''=\ul{\gn'/I_1}$ then the maps $\gx:\mcl{N}\lra \mcl{J}(P_{\ul{\gn'/I_1}})$ and $\gz:\mcl{J}(P_{\ul{\gn'/I_1}})\lra \mcl{N}$ are inverses of each other.
\end{remark}
\begin{prop}
	\label{prop:SandS1}
	Let $\ugl\in \Gl_0$ and $\gl_{i_0}$ be a part of $\ul{\gl}$ for some $1\leq i_0\leq k$ with $\gr_{i_0}=2$. Let  $\ul{\gn}=\big((\gl_1-2)^{\gr_1}>(\gl_2-2)^{\gr_2}>\ldots>(\gl_{i_0-1}-2)^{\gr_{i_0-1}}\geq (\gl_{i_0}-1)^{\gr_{i_0}}\geq \gl_{i_0+1}^{\gr_{i_0+1}}>\ldots>\gl_{k-1}^{\gr_{k-1}}>\gl_k^{\gr_k}\big)$.
	Let  $\mcl{M}=\{I\in \mcl{J}(P_{\ul{\gl}})\mid \max(I)\cap P_{(\gl_{i_0})}\neq\es\}$.
	For $I\in \mcl{M}\subseteq \mcl{J}(P_{\ul{\gl}})$ let $I_1=\gc(I)$ where $\gc$ is as defined in Theorem~\ref{theorem:LatticeIso}(2). Correspondingly we have partitions $\ul{\gl}',\ul{\gl}'',\ul{\gl'/I},\ul{\gn}',\ul{\gn}'',\ul{\gn'/I_1}$. Let $\ul{\gl}'''$ be the partition obtained by inserting a new part $\gl_{i_0}$ into the partition $\ul{\gl'/I}$. Let  $\mcl{N}=\{\ti{J}\in \mcl{J}(P_{\ul{\gl}'''})\mid \max(\ti{J})\cap P_{(\gl_{i_0})}\neq\es\}$.
	Let $L\in \mcl{M},K\in \mcl{T}=\{T\in \mcl{J}(P_{\ul{\gl}''})\mid \max(T)\cap P_{(\gl_{i_0})}\neq\es\}$ and $\ti{J}\in \mcl{N}$ be such that 
	$J\in \mcl{J}(P_{\ul{\gl^{'}/I}})$ with $\max(J)=\max(\ti{J})\bs P_{(\gl_{i_0})}$,
	$[J\cup K]_{\ul{\gl}}=L$ and $[\ti{J}]_{(\gl_{i_0})}=[K]_{(\gl_{i_0})}$. Let $L_1=\gc(L)\in \mcl{J}(P_{\ul{\gn}})$ 
	and let $K_1=\gc(K)\in \mcl{J}(P_{\ul{\gn}''})$. Let $J_1=\gf(\ti{J})\in \mcl{J}(P_{\ul{\gn^{'}/I_1}})$ where $\gf:\mcl{N}\lra \mcl{J}(P_{\ul{\gn^{'}/I_1}})$ is the lower lattice homomorphism as defined in  Proposition~\ref{prop:LowerLatticeHomoLamdaPPPtoNuModI}.
	Let $S=\max(K)\cap [J]_{\ul{\gl}''}$ and let $S_1=\{(a-1,\gl-2)\mid (a,\gl)\in S,\gl>\gl_{i_0}+1\}\cup\{(a,\gl)\mid (a,\gl)\in S,\gl<\gl_{i_0}+1\}$. Let $J_3=\gx(\ti{J})\in  \mcl{J}(P_{\ul{\gn^{'}/I_1}})$ where $\gx:\mcl{N}\lra \mcl{J}(P_{\ul{\gn^{'}/I_1}})$ is the upper lattice homomorphism as defined in Proposition~\ref{prop:UpperLatticeHomoLamdaPPPtoNuModI}.
	Then we have 
	\begin{enumerate}
		\item $[K_1]_{\ul{\gn}} \cup [J_1]_{\ul{\gn}}=[K_1]_{\ul{\gn}} \cup [J_3]_{\ul{\gn}}=L_1$.
		\item $\max(K_1)\cap [J_3]_{\ul{\gn}''}\supseteq \max(K_1)\cap [J_1]_{\ul{\gn}''} \supseteq S_1$.
		\item $S_1\cup\{(u,\gl_{i_0}-1)\}\supseteq \max(K_1)\cap [J_1]_{\ul{\gn}''} \supseteq S_1$ when $u<\gl_{i_0}-1$ and $\{(u,\gl_{i_0})\}=[\ti{J}]_{(\gl_{i_0})}=[K]_{(\gl_{i_0})}$.
		\item When $u=\gl_{i_0}-1$ and  $\{(u,\gl_{i_0})\}=[\ti{J}]_{(\gl_{i_0})}=[K]_{(\gl_{i_0})}$ then 
		$\max(K_1)\cap [J_1]_{\ul{\gn}''}=S_1$
	\end{enumerate}
\end{prop}
\begin{proof}
	We prove (1). We observe that if $\{(u,\gl_{i_0})\}=[\ti{J}]_{(\gl_{i_0})}=[K]_{(\gl_{i_0})}$ and $u<\gl_{i_0}-1$ then $(u,\gl_{i_0}-1)\in J_1\cup K_1$. Hence using Remark~\ref{remark:J1J3} we get that $J_1\cup K_1=J_3\cup K_1$, that is, they generate the same ideal in the big fundamental poset P defined in~\ref{Eq:FP}. If $u=\gl_{i_0}-1$ then we have $J_1=J_3\Ra J_1\cup K_1=J_3\cup K_1$.  Now we have $[\ti{J}]_{\ul{\gl}}\cup [K]_{\ul{\gl}}=L$. We observe that $\gc([\ti{J}]_{\ul{\gl}})=[J_3\cup \langle \{(u,\gl_{i_0}-1)\}\rangle]_{\ul{\gn}}, \gc([K]_{\ul{\gl}})=[K_1]_{\ul{\gn}}$. So we obtain $L_1=\gc(L)=\gc([\ti{J}]_{\ul{\gl}}\cup [K]_{\ul{\gl}})=\gc([\ti{J}]_{\ul{\gl}})\cup \gc([K]_{\ul{\gl}})=[J_3\cup K_1]_{\ul{\gn}}=[J_1\cup K_1]_{\ul{\gn}}$. Hence (1) follows.
	
	We prove (2). Let $(a,\gl)\in S,\gl>\gl_{i_0}+1$. Then $(a,\gl)\in \max(K)$. If $(a,\gl)\leq (v,\gd)\in \max(J)$ then either $\gd>\gl_{i_0}+1$ or $\gd<\gl_{i_0}-1$. We show that $\gd$ cannot be smaller than $\gl_{i_0}$. Since $(a,\gl),(u,\gl_{i_0})\in \max(K)$ we have $a>u,\gl-a>\gl_{i_0}-u$. Also $(v,\gd),(u,\gl_{i_0})\in \max(\ti{J})\Ra u>v,\gl_{i_0}-u>\gd-v$. So we have $a>u>v,\gl-a>\gl_{i_0}-u>\gd-v$ which implies that $(a,\gl)$ and $(v,\gd)$ are not comparable with each other. So $\gd>\gl_{i_0}$. If $\gd>\gl_{i_0}$ then $(a-1,\gl-2)\leq (v-1,\gd-2)$. Here $(a-1,\gl-2)\in \max(K_1)$ and $(v-1,\gd-2)\in \max(J_1)$. So $(a-1,\gl-2)\in \max(K_1)\cap [J_1]_{\ul{\gn}''}$.
	Now consider the case $(a,\gl)\in S,\gl<\gl_{i_0}-1$. Here also $(a,\gl)\in \max(K)$. If $(a,\gl)\leq (v,\gd)\in \max(J)$ then either $\gd>\gl_{i_0}+1$ or $\gd<\gl_{i_0}-1$. We show that $\gd$ cannot be bigger than $\gl_{i_0}$. Since $(a,\gl),(u,\gl_{i_0})\in \max(K)$ they are not comparable with each other. Since $(v,\gd),(u,\gl_{i_0})\in \max(\ti{J})$ they are not also comparable with each other. So if $\gd>\gl_{i_0}$ then $(a,\gl)$ and $(v,\gd)$ are not comparable with each other.
	So $\gd<\gl_{i_0}$. If $\gd<\gl_{i_0}$ and $(a,\gl)\leq (v,\gd)$, then we have $(a,\gl)\in \max(K_1),(v,\gd)\in \max(J_1)$ and therefore $(a,\gl)\in \max(K_1)\cap [J_1]_{\ul{\gn}''}\subseteq \max(K_1)\cap [J_3]_{\ul{\gn}''}$. Hence $S_1\subseteq \max(K_1)\cap [J_1]_{\ul{\gn}''}$. This proves (2).
	
	We prove (3),(4). Now suppose $(a,\gl)\in \max(K_1)\cap [J_1]_{\ul{\gn}''},\gl>\gl_{i_0}-1$ and $(a,\gl)\leq (b,\gd)\in \max(J_1)$.
	Now either $\gd<\gl_{i_0}-1$ or $\gd>\gl_{i_0}-1$. If $\gd>\gl_{i_0}-1$ then $(a+1,\gl+2)\leq (b+1,\gd+2)\in \max(J)$ and $(a+1,\gl+2)\in \max(K)$. So $(a,\gl)\in S_1$. If $\gd<\gl_{i_0}-1$ then $(b,\gd)\in \max(J)\subseteq \max(\ti{J})$.
	So $u>b,\gl_{i_0}-u>\gd-b$. Also $(a+1,\gl+2)\in \max(K)\Ra a+1>u,\gl-a+1>\gl_{i_0}-u$. So $a+1>u>b,\gl-a+1>\gl_{i_0}-u>\gd-b\Ra a>b,\gl-a>\gd-b$ which implies that $(a,\gl)$ and $(b,\gd)$ are not comparable. So the case $\gd<\gl_{i_0}-1$ does not occur. Now assume that $\gl<\gl_{i_0}-1$ and $(a,\gl)\leq (b,\gd)\in \max(J_1)$. Here again either $\gd<\gl_{i_0}-1$ or $\gd>\gl_{i_0}-1$. If $\gd<\gl_{i_0}-1$ then $(a,\gl)\in \max(K),(b,\gd)\in \max(J)$ and $(a,\gl)\leq (b,\gd)\Ra (a,\gl)\in S\Ra (a,\gl)\in S_1$. If $\gd>\gl_{i_0}-1$ then $(b+1,\gd+2)\in \max(J)\subseteq \max(\ti{J})$. So $b+1>u,\gd-b+1>\gl_{i_0}-u$. Also $(a,\gl)\in \max(K)\Ra u>a,\gl_{i_0}-u>\gl-a$. Therefore we have 
	$b+1>u>a,\gd-b+1>\gl_{i_0}-u>\gl-a\Ra b>a,\gd-b>\gl-a$ which implies that $(a,\gl)$ and $(b,\gd)$ are not comparable. So the case $\gd>\gl_{i_0}-1$ does not occur. Now we assume that $\gl=\gl_{i_0}-1$. 
	Then $(a,\gl)=(u,\gl_{i_0}-1)\in \max(K_1)\Ra u<\gl_{i_0}-1$. So we get that
	\equ{S_1\cup\{(u,\gl_{i_0}-1)\}\supseteq \max(K_1)\cap [J_1]_{\ul{\gn}''} \supseteq S_1\text{ if }u<\gl_{i_0}-1.}
	Otherwise if $u=\gl_{i_0}-1$  then $S_1=\max(K_1)\cap [J_1]_{\ul{\gn}''}$. This proves (3),(4).
	
	Hence the proposition follows.	
\end{proof}
\begin{prop}
	\label{prop:ExistenceofJK}
	Let $\ugl\in \Gl_0$ and $\gl_{i_0}$ be a part of $\ul{\gl}$ for some $1\leq i_0\leq k$ with $\gr_{i_0}=2$. Let  $\ul{\gn}=\big((\gl_1-2)^{\gr_1}>(\gl_2-2)^{\gr_2}>\ldots>(\gl_{i_0-1}-2)^{\gr_{i_0-1}}\geq (\gl_{i_0}-1)^{\gr_{i_0}}\geq \gl_{i_0+1}^{\gr_{i_0+1}}>\ldots>\gl_{k-1}^{\gr_{k-1}}>\gl_k^{\gr_k}\big)$.
	Let  $\mcl{M}=\{I\in \mcl{J}(P_{\ul{\gl}})\mid \max(I)\cap P_{(\gl_{i_0})}\neq\es\}$.
	For $I\in \mcl{M}\subseteq \mcl{J}(P_{\ul{\gl}})$ let $I_1=\gc(I)$ where $\gc$ is as defined in Theorem~\ref{theorem:LatticeIso}(2). Correspondingly we have partitions $\ul{\gl}',\ul{\gl}'',\ul{\gl'/I},\ul{\gn}',\ul{\gn}'',\ul{\gn'/I_1}$. Let $\ul{\gl}'''$ be the partition obtained by inserting a new part $\gl_{i_0}$ into the partition $\ul{\gl'/I}$. Let  $\mcl{N}=\{\ti{J}\in \mcl{J}(P_{\ul{\gl}'''})\mid \max(\ti{J})\cap P_{(\gl_{i_0})}\neq\es\}$. Let $\mcl{T}=\{T\in \mcl{J}(P_{\ul{\gl}''})\mid \max(T)\cap P_{(\gl_{i_0})}\neq\es\}$.
	Let $J_4\in \mcl{J}(P_{\ul{\gn^{'}/I_1}})$ and $K_2\in \mcl{J}(P_{\ul{\gn}''})$. Let $L_1=[J_4]_{\ul{\gn}}\cup [K_2]_{\ul{\gn}}\in \mcl{J}(P_{\ul{\gn}})$. Let $L\in \mcl{M}$ be such that $\gc(L)=L_1$ where the lattice homomorphism $\gc:\mcl{M}\lra \mcl{J}(P_{\ul{\gn}})$ is as defined in Theorem~\ref{theorem:LatticeIso}(2). Let $[L]_{(\gl_{i_0})}=\{(u,\gl_{i_0})\}\in \max(L)$. Let $K\in \mcl{T}\subseteq \mcl{J}(P_{\ul{\gl}''})$ be such that $K=\gc^{-1}(K_2)\cup [\langle \{(u,\gl_{i_0})\}\rangle]_{\ul{\gl}''}$. Let $K_1=\gc(K)\in \mcl{J}(P_{\ul{\gn}''})$ where the lattice homomorphism $\gc:\mcl{T}\lra \mcl{J}(P_{\ul{\gn}''})$ is as defined in Theorem~\ref{theorem:LatticeIso}(2) for the partitions $\ul{\gl}'',\ul{\gn}''$. Let $J_5=\gz(J_4)\in \mcl{N}$ where the lattice homomorphism $\gz:\mcl{J}(P_{\ul{\gn^{'}/I_1}})\lra \mcl{N}$ is as defined in Proposition~\ref{prop:ReverseLatticeMap}. Let $\ti{J}=J_5\cup [\langle \{(u,\gl_{i_0})\}\rangle]_{\ul{\gl}'''}$. Then $\ti{J}\in \mcl{N}$. Let $J\in \mcl{J}(P_{\ul{\gl^{'}/I}})$ be such that $\max(J)=\max(\ti{J})\bs \{(u,\gl_{i_0})\}$. Let $J_1\in \mcl{J}(P_{\ul{\gn^{'}/I_1}})$ be such that $\gf(\ti{J})=J_1$. Let $J_3=\gx(\ti{J})\in \mcl{J}(P_{\ul{\gn^{'}/I_1}})$ where the upper lattice homomorphism $\gx:\mcl{N}\lra \mcl{J}(P_{\ul{\gn^{'}/I_1}})$ is as defined in Proposition~\ref{prop:UpperLatticeHomoLamdaPPPtoNuModI}. Let $S=\max(K)\cap [J]_{\ul{\gl}''}$ and let $S_1=\{(a-1,\gl-2)\mid (a,\gl)\in S,\gl>\gl_{i_0}\}\cup \{(a,\gl)\mid (a,\gl)\in S,\gl<\gl_{i_0}\}$.
	Then we have 
	\begin{enumerate}
		\item $(u,\gl_{i_0})\in \max(K)$.
		\item $(u,\gl_{i_0})\in \max(\ti{J})$.
		\item $J_1\subseteq J_4\subseteq J_3$.
		\item $K_2\subseteq K_1$.
		\item $\max(K_2)\cap[J_4]_{\ul{\gn}''}\supseteq S_1$.
		\item $[K_1]_{\ul{\gn}}\cup [J_1]_{\ul{\gn}}=L_1$.
		\item Additional consequences are 
		\begin{itemize}
			\item $[K_1]_{\ul{\gn}}\cup [J_4]_{\ul{\gn}}=
			[K_1]_{\ul{\gn}}\cup [J_3]_{\ul{\gn}}=L_1$.
			\item $\max(K_1)\cap[J_3]_{\ul{\gn}''}\supseteq \max(K_1)\cap[J_4]_{\ul{\gn}''}\supseteq \max(K_1)\cap[J_1]_{\ul{\gn}''}\supseteq S_1$.
			\item $\max(K_2)\cap[J_3]_{\ul{\gn}''}\supseteq S_1$.
		\end{itemize}
	\end{enumerate}
\end{prop}
\begin{proof}
	We prove (1).
	$(u,\gl_{i_0})\in \max(L), L=\gc^{-1}([J_4]_{\ul{\gn}})\cup \gc^{-1}([K_2]_{\ul{\gn}})$. Hence $(u,\gl_{i_0})\in \max(\gc^{-1}([K_2]_{\ul{\gn}}) \cup [\langle \{(u,\gl_{i_0})\}\rangle]_{\ul{\gl}})=\max([K]_{\ul{\gl}})$. Now $\gl_{i_0}$ is a part of $\ul{\gl}''$ as well. Hence $(u,\gl_{i_0})\in \max(K)$. This proves (1).
	
	We prove (2). First observe that $\gl_{i_0}-1$ is a part of $\ul{\gn}$. Hence $[[J_4]_{\ul{\gn}}]_{(\gl_{i_0}-1)}=[J_4]_{(\gl_{i_0}-1)}$. Suppose $u<\gl_{i_0}-1$ and $\max([L_1]_{(\gl_{i_0}-1)})=\{(u,\gl_{i_0}-1)\}$. Then either $\max([J_4]_{(\gl_{i_0}-1)})$ is empty  or $\{(b,\gl_{i_0}-1)\}=\max([J_4]_{(\gl_{i_0}-1)})\leq \{(u,\gl_{i_0}-1)\}$, that is, $b\geq u$.  
	If $u=\gl_{i_0}-1$ then $\max([L_1]_{(\gl_{i_0}-1)})$ is empty and also $\max([J_4]_{(\gl_{i_0}-1)})$ is empty. 
	So if $\max([J_4]_{(\gl_{i_0}-1)})$ is empty then $\max([J_5])_{(\gl_{i_0})}=\{(\gl_{i_0}-1,\gl_{i_0})\}$ and $u\leq \gl_{i_0}-1$. If $\max([J_4]_{(\gl_{i_0}-1)})=\{(b,\gl_{i_0}-1)\}$ then $\max([J_5])_{(\gl_{i_0})}$ $=\{(b,\gl_{i_0})\}$ and $u\leq b$. Now $J_5\in \mcl{N}$. Hence $\ti{J}\in \mcl{N}$ and $(u,\gl_{i_0})\in \max(\ti{J})$. This proves (2).

	We prove (3). Since $J_5\subs \ti{J}$, we have from Equations~\ref{Eq:Comp1},~\ref{Eq:Comp2}
	\equa{\gz(J_4)\subseteq \ti{J}&\Llra [\gc^{-1}([J_4]_{\ul{\gn}'''})]_{\ul{\gl}'''} \subseteq \ti{J}\text{ by definition}\\
		&\Llra \gc^{-1}([J_4]_{\ul{\gn}'''}) \subseteq [\ti{J}]_{\ul{\gl}''''} \text{ under isomorphism } \\
		&\Llra [J_4]_{\ul{\gn}'''} \subseteq \gc([\ti{J}]_{\ul{\gl}''''}) \text{ under isomorphism}\\
		&\Ra [[J_4]_{\ul{\gn}'''}]_{\ul{\gn^{'}/I_1}} \subseteq [\gc([\ti{J}]_{\ul{\gl}''''})]_{\ul{\gn'/I_1}}\\
		&\Ra J_4\subseteq J_3. }
	The last implication follows because $J_4\in \mcl{J}(P_{\ul{\gn^{'}/I_1}}),\ul{\gn^{'}/I_1}\subseteq \ul{\gn}'''$ and hence $J_4=[[J_4]_{\ul{\gn}'''}]_{\ul{\gn^{'}/I_1}}$ and we have $J_3=[\gc([\ti{J}]_{\ul{\gl}''''})]_{\ul{\gn'/I_1}}$ by definition.   
	
	Now for any part $\gm_i$ of $\ul{\gn'/I_1}$, let $\max([J_4]_{(\gm_i)})=\{(b_i,\gm_i)\}$ if nonempty and chooose $b_i=\gm_i$ if
	$\max([J_4]_{(\gm_i)})$ is empty. Then we have for the ideal $J_5\in \mcl{J}(P_{\ul{\gl}'''})$, if $\gm_i\geq \gl_{i_0}-1$
	then $\max([J_5]_{(\gm_i+2)})=\{(b_i+1,\gm_i+2)\}$ and if $\gm_i \leq \gl_{i_0}-1$ then $\max([J_5]_{(\gm_i)})=\max([J_4]_{(\gm_i)})$. We also have $\max([J_5]_{(\gl_{i_0})})=\{(b,\gl_{i_0})\}$ where $b=b_i$ such that $\gm_i=\gl_{i_0}-1$. Therefore we have 
	\fo{10}{10}{
		\equ{\max(J_5)\subseteq \{(b_i+1,\gm_i+2)\mid \gm_i\geq \gl_{i_0}-1\}\cup \{(b,\gl_{i_0})\}\cup  \{(b_i,\gm_i)\mid \gm_i\leq \gl_{i_0}-1,b_i<\gm_i\} \subseteq J_5.}
	}
	So we have 
	\fo{9}{9}{
		\equ{(u,\gl_{i_0})\in \max(\ti{J})\subseteq \{(b_i+1,\gm_i+2)\mid \gm_i> \gl_{i_0}-1\}\cup \{(u,\gl_{i_0})\} \cup \{(b_i,\gm_i)\mid \gm_i< \gl_{i_0}-1,b_i<\gm_i\}\subseteq \ti{J}.}
	}
	So 
	\fo{10}{10}{
		\equ{\max(J)=\max(\ti{J})\bs\{(u,\gl_{i_0})\}\subseteq \{(b_i+1,\gm_i+2)\mid \gm_i> \gl_{i_0}-1\}\cup \{(b_i,\gm_i)\mid \gm_i< \gl_{i_0}-1,b_i<\gm_i\}}
	}
	which generates an ideal contained in $J_5$. Therefore \equ{\max(J_1)\subseteq \{(b_i,\gm_i)\mid \gm_i>\gl_{i_0}-1\} \cup \{(b_i,\gm_i)\mid \gm_i<\gl_{i_0}-1,b_i<\gm_i\}} which generates an ideal contained in $J_4$. So $J_1\subseteq J_4$. This proves (3).
	
	We prove (4). Since $K=\gc^{-1}(K_2)\cup [\langle \{(u,\gl_{i_0})\}\rangle]_{\ul{\gl}''}$ and $K_1=\gc(K)$ we have $K_2\subseteq K_1$. This proves (4).
	
	We prove (5). Let $(a,\gl)\in S,\gl>\gl_{i_0}+1$. Then $(a,\gl)\in [J]_{\ul{\gl}''}\Ra \gl$ is a part of $\ul{\gl}''$ and $(a,\gl)\leq (b,\gm)$ for some $(b,\gm)\in \max(J)$. Now $(a,\gl),(u,\gl_{i_0})\in \max(K)$ and hence are not comparable. Also $(u,\gl_{i_0}),(b,\gm)\in \max(\ti{J})$ are not comparable. Hence if $\gm<\gl_{i_0}-1$ then we have $a>u>b,\gl-a>\gl_{i_0}-u>\gm-b$. So $(a,\gl)$ and $(b,\gm)$ are not  comparable. Therefore $\gm>\gl_{i_0}+1$ and $(b-1,\gm-2)\in \max(J_1)$. Then $(a-1,\gl-2)\leq (b-1,\gm-2)$. So $(a-1,\gl-2)\in [J_1]_{\ul{\gn}''}$. If $(a,\gl)\in \max(K)$ and $\gl>\gl_{i_0}+1$ then $(a,\gl)$ is not comparable with $(u,\gl_{i_0})$. We have $K=\gc^{-1}(K_2)\cup [\langle \{(u,\gl_{i_0})\}\rangle]_{\ul{\gl}''}$. So $(a,\gl)\in \gc^{-1}(K_2)\Ra  (a-1,\gl-2)\in \max(K_1)\cap K_2\subseteq \max(K_2)$. Hence $(a-1,\gl-2)\in \max(K_2)\cap [J_4]_{\ul{\gn}''}$.
	Now let $(a,\gl)\in S,\gl<\gl_{i_0}-1$. Then $(a,\gl)\in [J]_{\ul{\gl}''}\Ra \gl$ is a part of $\ul{\gl}''$ and $(a,\gl)\leq (b,\gm)$ for some $(b,\gm)\in \max(J)$.  Now $(a,\gl),(u,\gl_{i_0})\in \max(K)$ and hence are not comparable. Also $(u,\gl_{i_0}),(b,\gm)\in \max(\ti{J})$ are not comparable. Hence if $\gm>\gl_{i_0}-1$ then we have 
	$b>u>a,\gm-b>\gl_{i_0}-u>\gl-a$.  So $(a,\gl)$ and $(b,\gm)$ are not  comparable. Therefore $\gm<\gl_{i_0}-1$ and $(b,\gm)\in \max(J_1)$. So $(a,\gl)\leq (b,\gm)$ implies $(a,\gl)\in [J_1]_{\ul{\gn}''}$. If $(a,\gl)\in \max(K)$ and $\gl<\gl_{i_0}+1$ then $(a,\gl)$ is not comparable with $(u,\gl_{i_0})$. We have $K=\gc^{-1}(K_2)\cup [\langle \{(u,\gl_{i_0})\}\rangle]_{\ul{\gl}''}$. So $(a,\gl)\in \gc^{-1}(K_2)\Ra  (a,\gl)\in \max(K_1)\cap K_2\subseteq \max(K_2)$. Hence $(a,\gl)\in \max(K_2)\cap [J_4]_{\ul{\gn}''}$.
	This implies that $S_1\subseteq \max(K_2)\cap [J_4]_{\ul{\gn}''}$. This proves (5).	
	
	We prove (6). We observe that $[K]_{\ul{\gl}}\cup [J]_{\ul{\gl}}=[K]_{\ul{\gl}}\cup [\ti{J}]_{\ul{\gl}}=L$. So using Proposition~\ref{prop:SandS1}(1), we obtain $[K_1]_{\ul{\gn}}\cup [J_1]_{\ul{\gn}}=L_1$. 
	
	Now the additional consequences in (7) also follow. This completes the proof of the proposition.
\end{proof}
\begin{prop}
	\label{prop:JKJ1K1Multiple}
	Let $\ugl\in \Gl_0$ and $\gl_{i_0}$ be a part of $\ul{\gl}$ for some $1\leq i_0\leq k$ with $\gr_{i_0}=2$. Let  $\ul{\gn}=\big((\gl_1-2)^{\gr_1}>(\gl_2-2)^{\gr_2}>\ldots>(\gl_{i_0-1}-2)^{\gr_{i_0-1}}\geq (\gl_{i_0}-1)^{\gr_{i_0}}\geq \gl_{i_0+1}^{\gr_{i_0+1}}>\ldots>\gl_{k-1}^{\gr_{k-1}}>\gl_k^{\gr_k}\big)$.
	Let  $\mcl{M}=\{I\in \mcl{J}(P_{\ul{\gl}})\mid \max(I)\cap P_{(\gl_{i_0})}\neq\es\}$.
	For $I\in \mcl{M}\subseteq \mcl{J}(P_{\ul{\gl}})$ let $I_1=\gc(I)$ where $\gc$ is the lattice homomorphism defined in Theorem~\ref{theorem:LatticeIso}(2). Correspondingly we have partitions $\ul{\gl}',\ul{\gl}'',\ul{\gl'/I},\ul{\gn}',\ul{\gn}'',\ul{\gn'/I_1}$.
	\begin{enumerate}
		\item \equ{\frac{\mid \grppp\mid }{\mid \grpppI\mid}=q\frac{\mid \grnppp\mid }{\mid \grnpppI\mid}.}
		\item For $L\in \mcl{M}$ and $L_1=\gc(L)\in \mcl{J}(P_{\ul{\gn}})$, we have \equ{\mid(\grpp)_L\mid=q^{\#(\text{parts of }\ul{\gl} \text{ greater than equal to }\gl_{i_0})}\mid (\grnpp)_{L_1}\mid.} 
	\end{enumerate}	
\end{prop}
\begin{proof}
	We prove (1).
	In Equations~[\ref{Eq:One}-\ref{Eq:Five}], we observe that 
	\equ{\mid \ul{\gl}' \mid-\mid \ul{\gl^{'}/I}\mid=t_1-v_1=1+\mid \ul{\gn}'\mid -\mid \ul{\gn^{'}/I_1} \mid.}	
	
	We prove (2). We observe that if $L$ corresponds to \equ{(\grpp)_L=\us{i=1}{\os{i_0-1}{\oplus}}\gp^{r_i}(\RR {\gl_i})^{\gr_i}\oplus \gp^{r_{i_0}}(\RR {\gl_{i_0}})^{\gr_{i_0}}\oplus \us{i=i_0+1}{\os{k}{\oplus}}\gp^{r_i}(\RR {\gl_i})^{\gr_i}} then 
	$L_1$ corresponds to \equ{(\grnpp)_{L_1}=\us{i=1}{\os{i_0-1}{\oplus}}\gp^{r_i-1}(\RR {\gl_i-2})^{\gr_i}\oplus \gp^{r_{i_0}}(\RR {\gl_{i_0}-1})^{\gr_{i_0}}\oplus \us{i=i_0+1}{\os{k}{\oplus}}\gp^{r_i}(\RR {\gl_i})^{\gr_i}.}
	So \equ{\mid(\grpp)_L\mid=q^{\us{i=1}{\os{i_0}{\sum}}\gr_i}\mid (\grnpp)_{L_1}\mid=q^{\#(\text{parts of }\ul{\gl} \text{ greater than equal to }\gl_{i_0})}\mid (\grnpp)_{L_1}\mid.}
	
	This proves the proposition. 
\end{proof}
\begin{proof}[Proof of Theorem~\ref{theorem:LambdaNu}]
	We consider the case when $\max(I)\cap P_{(\gl_{i_0})}\neq \es, I\in \mcl{J}(P_{\ul{\gl}})$. These ideals are in bijection with the ideals in $\mcl{J}(P_{\ul{\gn}})$ using Theorem~\ref{theorem:LatticeIso}(2). We have to prove the identity in Equation~\ref{Eq:FiveFive}.

\begin{center}
	*****************************************************************
\end{center}
We consider the scenario \equ{\{(v,\gl_{i_0})\}=\max([J]_{(\gl_{i_0})})< \max([K]_{(\gl_{i_0})})=\{(u,\gl_{i_0})\},\text{ that is, } v>u}
or  \equ{\es=\max([J]_{(\gl_{i_0})})< \max([K]_{(\gl_{i_0})})=\{(u,\gl_{i_0})\}}
where $J\in \mcl{J}(P_{\ul{\gl^{'}/I}}), K\in \mcl{J}(P_{\ul{\gl}''})$ is such that $\max(K)\cap P_{(\gl_{i_0})}\neq \es$. Let $\ul{\gl}'''$ be the partition obtained by inserting the part $\gl_{i_0}$ into $\ul{\gm^{'}/I}$. Let $\ti{J}\in \mcl{J}(P_{\ul{\gl}'''})$ be an ideal such that $(u,\gl_{i_0})\in \max(\ti{J})$. Let $J\in \mcl{J}(P_{\ul{\gl^{'}/I}})$ be such that $\max(J)=\max(\ti{J})\bs \{(u,\gl_{i_0})\}$. Then clearly $\max([J]_{(\gl_{i_0})})< \max([K]_{(\gl_{i_0})})=\{(u,\gl_{i_0})\}$. Now consider ideals $J_2\in \mcl{J}(P_{\ul{\gl^{'}/I}})$ such that $J\subseteq J_2\subseteq \ti{J}\bs P_{(\gl_{i_0})}\sbnq \ti{J}$. Note that for every ideal $J_2\in \mcl{J}(P_{\ul{\gl^{'}/I}})$ such that $\max([J_2]_{(\gl_{i_0})})< \max([K]_{(\gl_{i_0})})=\{(u,\gl_{i_0})\}$ there is a unique ideal $\ti{J}$ with $(u,\gl_{i_0})\in \max(\ti{J})$ such that  $J\subseteq J_2\subseteq \ti{J}\bs P_{(\gl_{i_0})}\sbnq \ti{J}$ where $\max(J)=\max(\ti{J})\bs \{(u,\gl_{i_0})\}$.
The unique ideal $\ti{J}\in \mcl{J}(P_{\ul{\gl}'''})$ is given by $\ti{J}=[J_2]_{\ul{\gl}'''}\cup[\langle \{(u,\gl_{i_0})\}\rangle]_{\ul{\gl}'''}$.

For such an ideal $J_2$ we have $[J_2\cup K]_{\ul{\gl}'}=[J\cup K]_{\ul{\gl}'}$ since $J\cup K,J_2\cup K$ define the same ideal in the big fundamental poset $P$ defined in Equation~\ref{Eq:FP}. Also we have $\max(K)\bs [J_2]_{\ul{\gl}''}=\max(K)\bs[J]_{\ul{\gl}''}$. So for all $J\subseteq J_2\subseteq \ti{J}\bs P_{(\gl_{i_0})}$ \equ{\ga^{\ul{\gl}}_{I,J_2,K}(q)=\ga^{\ul{\gl}}_{I,J,K}(q)=\mid(\grpp)_{J\cup K}\mid\us{\gl\in \ul{\gl}'' \text{ such that there exists }(u,\gl)\in \max(K)\bs[J]_{\ul{\gl}''}}{\prod}(1-\frac1{q^{m(\gl)}}).}
where $m(\gl)$ is the multiplicity of the part $\gl$ in $\ul{\gl}''$.

Now we observe that 
\equa{\us{J\subseteq J_2 \subseteq \ti{J}\bs P_{(\gl_{i_0})}}{\bigsqcup}\mid X^{\ul{\gl}}_{I,J_2,K}\mid&=\frac{\mid \grppp\mid }{\mid \grpppI\mid}\mid(\grpppp)^*_K\mid \bigg(\us{J\subseteq J_2 \subseteq \ti{J}\bs P_{(\gl_{i_0})}}{\sum}\mid (\grpppI)^*_{J_2}\mid\bigg)\\
	&=\frac{\mid \grppp\mid }{\mid \grpppI\mid}\mid(\grpppp)^*_K\mid \frac{\mid (\mcl{A}_{\ul{\gl}'''})^*_{\ti{J}}\mid}{(q^{\gl_{i_0}-u}-q^{\gl_{i_0}-u-1})}.} 

Therefore we get that 
\equa{&\us{J\subseteq J_2 \subseteq \ti{J}\bs P_{(\gl_{i_0})}}{\sum}\frac{\mid X^{\ul{\gl}}_{I,J_2,K}\mid}{\ga^{\ul{\gl}}_{I,J_2,K}(q)}=\frac{\frac{\mid \grppp\mid }{\mid \grpppI\mid}\frac1{\mid(\grpp)_{J\cup K}\mid}\frac{\mid (\mcl{A}_{\ul{\gl}'''})^*_{\ti{J}}\mid}{(q^{\gl_{i_0}-u}-q^{\gl_{i_0}-u-1})}\mid(\grpppp)^*_K\mid}{\us{\gl\in \ul{\gl}'' \text{ such that there exists }(u,\gl)\in \max(K)\bs[J]_{\ul{\gl}''}}{\prod}(1-\frac1{q^{m(\gl)}})}\\&=
	\frac{\mid \grppp\mid }{\mid \grpppI\mid}\frac{1}{\mid(\grpp)_{J\cup K}\mid}\frac{\mid (\mcl{A}_{\ul{\gl}'''})^*_{\ti{J}}\mid}{(q^{\gl_{i_0}-u}-q^{\gl_{i_0}-u-1})}\mid(\grpppp)_K\mid\us{\os{\gl\in \ul{\gl}''}{\os{ \text{ such that there exists }}{(u,\gl)\in \max(K)\cap[J]_{\ul{\gl}''}}}}{\prod}(1-\frac1{q^{m(\gl)}}).}
\begin{center}
	*****************************************************************
\end{center}
Let $I_1=\gc(I)\in \mcl{J}(P_{\ul{\gn}})$ where the lattice homomorphism $\gc:\mcl{M}=\{T\in \mcl{J}(P_{\ul{\gl}})\mid \max(T)\cap P_{(\gl_{i_0})}\neq\es\} \lra \mcl{J}(P_{\ul{\gn}})$ is as defined in Theorem~\ref{theorem:LatticeIso}(2). Let $J_4\in  \mcl{J}(P_{\ul{\gn^{'}/I_1}}), K_2\in  \mcl{J}(P_{\ul{\gn}''})$ be any two ideals. Then using Proposition~\ref{prop:ExistenceofJK}, there are unique ideals $\ti{J}\in \mcl{N}=\{T\in \mcl{J}(P_{\ul{\gl}'''})\mid \max(T)\cap P_{(\gl_{i_0})}\neq\es\}\sbnq \mcl{J}(P_{\ul{\gl}'''})$ and $K\in \mcl{T}=\{T\in \mcl{J}(P_{\ul{\gl}''})\mid \max(T)\cap P_{(\gl_{i_0})}\neq\es\} \sbnq \mcl{J}(P_{\ul{\gl}''})$ such that $\max([\ti{J}]_{\gl_{i_0}})=\max([K]_{(\gl_{i_0})})$.
\begin{center}
	*****************************************************************
\end{center}
Now we begin with an ideal $\ti{J}\in \mcl{N}$ and an ideal  $K\in \mcl{T}$ such that $\max([\ti{J}]_{\gl_{i_0}})=\max([K]_{(\gl_{i_0})})=\{(u,\gl_{i_0})\}$.
Let $J\in \mcl{J}(P_{\ul{\gl^{'}/I}})$ be such that $\max(J)=\max(\ti{J})\bs P_{(\gl_{i_0})}$,
$[J\cup K]_{\ul{\gl}}=L$. Let $L_1=\gc(L)\in \mcl{J}(P_{\ul{\gn}})$ and $K_1=\gc(K)\in \mcl{J}(P_{\ul{\gn}''})$. Let $J_1=\gf(\ti{J})\in \mcl{J}(P_{\ul{\gn^{'}/I_1}})$. Let $S=\max(K)\cap [J]_{\ul{\gl}''}$ and let $S_1=\{(a-1,\gl-2)\mid (a,\gl)\in S,\gl>\gl_{i_0}\}\cup \{(a,\gl)\mid (a,\gl)\in S,\gl<\gl_{i_0}\}$. Let $J_3=\gx(\ti{J})$ where the upper lattice homomophism $\gx:\mcl{N}\lra \mcl{J}(P_{\ul{\gn^{'}/I_1}})$ is as defined in Proposition~\ref{prop:UpperLatticeHomoLamdaPPPtoNuModI}. Let $J_4\in  \mcl{J}(P_{\ul{\gn^{'}/I_1}})$ denote  any ideal such that $J_1\subseteq J_4\subseteq J_3$. Let $K_2\in  \mcl{J}(P_{\ul{\gn}''})$ denote any ideal such that $K_2\subseteq K_1$. Then we prove the following.

\fo{10}{10}{ 
	\equan{SixSix}{&
		\frac{\mid \grppp\mid }{\mid \grpppI\mid}\frac{1}{\mid(\grpp)_{L}\mid}\frac{\mid (\mcl{A}_{\ul{\gl}'''})^*_{\ti{J}}\mid}{(q^{\gl_{i_0}-u}-q^{\gl_{i_0}-u-1})}\mid(\grpppp)_K\mid\us{\os{\gl\in \ul{\gl}''}{\os{ \text{ such that there exists }}{(u,\gl)\in \max(K)\cap[J]_{\ul{\gl}''}}}}{\prod}(1-\frac1{q^{m(\gl)}})=\\
		&\frac{\mid \grnppp\mid }{\mid \grnpppI\mid}\frac{1}{\mid(\grnpp)_{L_1}\mid}\us{\os{J_1\subseteq J_4\subseteq J_3,K_2\subseteq K_1}{\us{\max(K_2)\cap [J_4]_{\ul{\gn}''}\supseteq S_1}{[J_4\cup K_2]_{\ul{\gn}}=L_1}}}{\sum}\bigg(\mid (\grnpppI)^*_{J_4}\mid \mid (\grnpppp)_{K_2}\mid\us{\os{\gl\in \ul{\gn}''}{\os{ \text{ such that there exists }}{(u,\gl)\in \max(K_2)\cap[J_4]_{\ul{\gn}''}}}}{\prod}(1-\frac1{q^{m(\gl)}})\bigg).}
}

$[J_4\cup K_2]_{\ul{\gn}}=L_1\Ra [J_4]_{\ul{\gn}''}\cup K_2=[L_1]_{\ul{\gn}''}$. Also we have $[J_4]_{\ul{\gn}''}\cup K_1=[L_1]_{\ul{\gn}''}$. Hence $[L_1]_{\ul{\gn}''}\bs [J_4]_{\ul{\gn}''}=K_1\bs [J_4]_{\ul{\gn}''}=K_2\bs [J_4]_{\ul{\gn}''}$. Now we show that $\max(K_1)\bs [J_4]_{\ul{\gn}''}=\max(K_2)\bs [J_4]_{\ul{\gn}''}$. Since $K_2\subseteq K_1$ we have $\max(K_1)\cap K_2\subseteq \max(K_2)$. Now $(c,\gm)\in \max(K_1)\bs [J_4]_{\ul{\gn}''}\Ra (c,\gm)\in K_1\bs [J_4]_{\ul{\gn}''}=K_2\bs [J_4]_{\ul{\gn}''}$. So $(c,\gm)\in \max(K_1)\cap K_2\subseteq \max(K_2)\Ra (c,\gm)\in \max(K_2)\bs [J_4]_{\ul{\gn}''}$. Hence $\max(K_1)\bs [J_4]_{\ul{\gn}''}$ $\subseteq \max(K_2)\bs [J_4]_{\ul{\gn}''}$. Now $(c,\gm)\in \max(K_2)\bs [J_4]_{\ul{\gn}''}$ and $[J_4]_{\ul{\gn}''}\cup K_2=[L_1]_{\ul{\gn}''} \Ra (c,\gm)\in \max([L_1]_{\ul{\gn}''})\bs [J_4]_{\ul{\gn}''}$. Now $[J_4]_{\ul{\gn}''}\cup K_1=[L_1]_{\ul{\gn}''} \Ra (c,\gm)\in \max(K_1)\bs [J_4]_{\ul{\gn}''}$. So $\max(K_1)\bs [J_4]_{\ul{\gn}''}= \max(K_2)\bs [J_4]_{\ul{\gn}''}$.	

Now let $K_3\in \mcl{J}(P_{\ul{\gn}''})$ be such that $\max(K_3)=\max(K_1)\bs [J_4]_{\ul{\gn}''} \cup S_1$.
Note that $S_1\subseteq \max(K_1)\cap [J_4]_{\ul{\gn}''}$ using Proposition~\ref{prop:ExistenceofJK}(7). Hence $\max(K_3)\subseteq \max(K_1)$ is indeed an antichain.
Since $L_1$ is generated by $[J_4]_{\ul{\gn}}$ and $\max(K_1)\bs [J_4]_{\ul{\gn}}= \max(K_1)\bs [J_4]_{\ul{\gn}''} =\max(K_3)\bs [J_4]_{\ul{\gn}''}=\max(K_3)\bs [J_4]_{\ul{\gn}}$ we have $L_1=[K_3]_{\ul{\gn}}\cup [J_4]_{\ul{\gn}}$.
Also we observe that $\max(K_3)\cap [J_4]_{\ul{\gn}''}\supseteq S_1$, in fact equal to $S_1$.
Now let $K_2\in \mcl{J}(P_{\ul{\gn}''})$ be any ideal such that $K_3\subseteq K_2\subseteq K_1$. Then we have $L_1=[K_2]_{\ul{\gn}}\cup [J_4]_{\ul{\gn}}$ and $(\grnpppp)^*_{K_3}\subseteq (\grnpppp)^*_{K_2}\subseteq (\grnpppp)^*_{K_1}$. So $\max(K_3)\subseteq \max(K_1) \Ra \max(K_3)\subseteq \max(K_2)$. So $\max(K_2)\cap [J_4]_{\ul{\gn}''}\supseteq \max(K_3)\cap [J_4]_{\ul{\gn}''}\supseteq S_1$. 

Conversely if $K_2\subseteq K_1$ such that $L_1=[K_2]_{\ul{\gn}}\cup [J_4]_{\ul{\gn}}$ and $\max(K_2)\cap [J_4]_{\ul{\gn}''}\supseteq S_1$ then $\max(K_2)\supseteq S_1\cup \max(K_1)\bs [J_4]_{\ul{\gn}''}=\max(K_3)$. So we have 
for a fixed $J_4$ such that $J_1\subseteq J_4\subseteq J_3$,

\equa{&\us{\os{K_2\subseteq K_1}{\us{\max(K_2)\cap [J_4]_{\ul{\gn}''}\supseteq S_1}{[J_4\cup K_2]_{\ul{\gn}}=L_1}}}{\sum}\bigg(\mid (\grnpppI)^*_{J_4}\mid \mid (\grnpppp)_{K_2}\mid\us{\os{\gl\in \ul{\gn}''}{\os{ \text{ such that there exists }}{(u,\gl)\in \max(K_2)\cap[J_4]_{\ul{\gn}''}}}}{\prod}(1-\frac1{q^{m(\gl)}})\bigg)}
\equa{&=\mid (\grnpppI)^*_{J_4}\mid\us{K_3\subseteq K_2\subseteq K_1}{\sum}\bigg(\mid (\grnpppp)_{K_2}\mid\us{\os{\gl\in \ul{\gn}''}{\os{ \text{ such that there exists }}{(u,\gl)\in \max(K_2)\cap[J_4]_{\ul{\gn}''}}}}{\prod}(1-\frac1{q^{m(\gl)}})\bigg)\\
&= \mid (\grnpppI)^*_{J_4}\mid\us{K_3\subseteq K_2\subseteq K_1}{\sum}\bigg(\frac{\mid (\grnpppp)^*_{K_2}\mid}{\us{\os{\gl\in \ul{\gn}''}{\os{ \text{ such that there exists }}{(u,\gl)\in \max(K_2)\bs[J_4]_{\ul{\gn}''}}}}{\prod}(1-\frac1{q^{m(\gl)}})}\bigg)\\&= \mid (\grnpppI)^*_{J_4}\mid\us{K_3\subseteq K_2\subseteq K_1}{\sum}\bigg(\frac{\mid (\grnpppp)^*_{K_2}\mid}{\us{\os{\gl\in \ul{\gn}''}{\os{ \text{ such that there exists }}{(u,\gl)\in \max(K_1)\bs[J_4]_{\ul{\gn}''}}}}{\prod}(1-\frac1{q^{m(\gl)}})}\bigg)\\&=
	\frac{\mid (\grnpppI)^*_{J_4}\mid}{\bigg(\us{\os{\gl\in \ul{\gn}''}{\os{ \text{ such that there exists }}{(u,\gl)\in \max(K_1)\bs[J_4]_{\ul{\gn}''}}}}{\prod}(1-\frac1{q^{m(\gl)}})\bigg)}\us{K_3\subseteq K_2\subseteq K_1}{\sum}\bigg(\mid (\grnpppp)^*_{K_2}\mid\bigg).}
Now $\max(K_3)\subseteq \max(K_2)\cap \max(K_1)$ implies that 
\equa{\us{K_3\subseteq K_2\subseteq K_1}{\sum}\mid (\grnpppp)^*_{K_2}\mid=\mid (\grnpppp)_{K_1}\mid\bigg(\us{\os{\gl\in \ul{\gn}''}{\os{ \text{ such that there exists }}{(u,\gl)\in \max(K_3)}}}{\prod}(1-\frac1{q^{m(\gl)}})\bigg).}
Hence 
\equa{&\frac{\mid (\grnpppI)^*_{J_4}\mid}{\bigg(\us{\os{\gl\in \ul{\gn}''}{\os{ \text{ such that there exists }}{(u,\gl)\in \max(K_1)\bs[J_4]_{\ul{\gn}''}}}}{\prod}(1-\frac1{q^{m(\gl)}})\bigg)}\us{K_3\subseteq K_2\subseteq K_1}{\sum}\bigg(\mid (\grnpppp)^*_{K_2}\mid\bigg)\\&=\frac{\mid (\grnpppI)^*_{J_4}\mid\mid (\grnpppp)_{K_1}\mid}{\bigg(\us{\os{\gl\in \ul{\gn}''}{\os{ \text{ such that there exists }}{(u,\gl)\in \max(K_1)\bs[J_4]_{\ul{\gn}''}}}}{\prod}(1-\frac1{q^{m(\gl)}})\bigg)}\bigg(\us{\os{\gl\in \ul{\gn}''}{\os{ \text{ such that there exists }}{(u,\gl)\in \max(K_3)}}}{\prod}(1-\frac1{q^{m(\gl)}})\bigg)\\&=\mid (\grnpppI)^*_{J_4}\mid\mid (\grnpppp)_{K_1}\mid\bigg(\us{\os{\gl\in \ul{\gn}''}{\os{ \text{ such that there exists }}{(u,\gl)\in S_1}}}{\prod}(1-\frac1{q^{m(\gl)}})\bigg).}

Now $J_1\subseteq J_4\subseteq J_3$ and $\max(J_1)\subseteq \max(J_3)\Ra \max(J_1)\subseteq \max(J_4)$. Hence 
\equa{\us{J_1\subseteq J_4\subseteq J_3}{\sum}\mid (\grnpppI)^*_{J_4}\mid=\mid (\grnpppI)_{J_3}\mid\bigg(\us{\os{\gl\in \ul{\gn'/I_1}}{\os{ \text{ such that there exists }}{(u,\gl)\in \max(J_1)}}}{\prod}(1-\frac1{q})\bigg).}
So
\fo{9}{9}{\equa{&\frac{\mid \grnppp\mid }{\mid \grnpppI\mid}\frac{1}{\mid(\grnpp)_{L_1}\mid}\us{\os{J_1\subseteq J_4\subseteq J_3,K_2\subseteq K_1}{\us{\max(K_2)\cap [J_4]_{\ul{\gn}''}\supseteq S_1}{[J_4\cup K_2]_{\ul{\gn}}=L_1}}}{\sum}\bigg(\mid (\grnpppI)^*_{J_4}\mid \mid (\grnpppp)_{K_2}\mid\us{\os{\gl\in \ul{\gn}''}{\os{ \text{ such that there exists }}{(u,\gl)\in \max(K_2)\cap[J_4]_{\ul{\gn}''}}}}{\prod}(1-\frac1{q^{m(\gl)}})\bigg)\\&=\frac{\mid \grnppp\mid }{\mid \grnpppI\mid}\frac{1}{\mid(\grnpp)_{L_1}\mid}\mid (\grnpppI)_{J_3}\mid\bigg(\us{\os{\gl\in \ul{\gn'/I_1}}{\os{ \text{ such that there exists }}{(u,\gl)\in \max(J_1)}}}{\prod}(1-\frac1{q})\bigg)\mid (\grnpppp)_{K_1}\mid\bigg(\us{\os{\gl\in \ul{\gn}''}{\os{ \text{ such that there exists }}{(u,\gl)\in S_1}}}{\prod}(1-\frac1{q^{m(\gl)}})\bigg)\\&=\frac{\mid \grppp\mid }{\mid \grpppI\mid}\frac{q^{\os{-1+\#(\text{parts of }\ul{\gl}}{ \text{ greater than equal to }\gl_{i_0})}}}{\mid(\grpp)_{L}\mid}\mid (\grnpppI)_{J_3}\mid\bigg(\us{\os{\gl\in \ul{\gn'/I_1}}{\os{ \text{ such that there exists }}{(u,\gl)\in \max(J_1)}}}{\prod}(1-\frac1{q})\bigg)\mid (\grnpppp)_{K_1}\mid\bigg(\us{\os{\gl\in \ul{\gn}''}{\os{ \text{ such that there exists }}{(u,\gl)\in S_1}}}{\prod}(1-\frac1{q^{m(\gl)}})\bigg).}}
The last equality follows using Proposition~\ref{prop:JKJ1K1Multiple}.

Now $(u,\gl)\in S \Llra (u-1,\gl-2)\in S_1$ if $\gl>\gl_{i_0}$ and $(u,\gl)\in S_1$ if $\gl<\gl_{i_0}$. Here we observe that if $\gl>\gl_{i_0}$ then $m_{\ul{\gl}''}(\gl)$, the multiplicity of $\gl\in \ul{\gl}''$ equals $m_{\ul{\gn}''}(\gl-2)$, the mulitplicity of $(\gl-2)\in \ul{\gn}''$ and if $\gl<\gl_{i_0}$ then $m_{\ul{\gl}''}(\gl)=m_{\ul{\gn}''}(\gl)$. 

Hence we have 
\equ{\bigg(\us{\os{\gl\in \ul{\gn}''}{\os{ \text{ such that there exists }}{(u,\gl)\in S_1}}}{\prod}(1-\frac1{q^{m(\gl)}})\bigg)=\mid\bigg(\us{\os{\gl\in \ul{\gl}''}{\os{ \text{ such that there exists }}{(u,\gl)\in S}}}{\prod}(1-\frac1{q^{m(\gl)}})\bigg).}

We also see that $(u,\gl)\in \max(\ti{J}),\gl\neq \gl_{i_0},\gl>\gl_{i_0}+1$ if and only if $(u-1,\gl-2)\in \max(J_1)$ and $(u,\gl)\in \max(\ti{J}),\gl\neq \gl_{i_0},\gl<\gl_{i_0}-1$ if and only if $(u,\gl)\in \max(J_1)$.
Hence we have 
\equ{\bigg(\us{\os{\gl\in \ul{\gn'/I_1}}{\os{ \text{ such that there exists }}{(u,\gl)\in \max(J_1)}}}{\prod}(1-\frac1{q})\bigg)=\bigg(\us{\os{\gl\in \ul{\gl'/I}}{\os{ \text{ such that there exists }}{(u,\gl)\in \max(J)}}}{\prod}(1-\frac1{q})\bigg)=\frac{\mid (\mcl{A}_{\ul{\gl}'''})^*_{\ti{J}}\mid}{\mid (\mcl{A}_{\ul{\gl}'''})_{\ti{J}}\mid (1-\frac 1q)}.}

So to prove Equation~\ref{Eq:SixSix}, it is enough to prove the following.
\equan{SevenSeven}{q^{\os{-1+\#(\text{parts of }\ul{\gl}}{ \text{ greater than equal to }\gl_{i_0})}}\mid (\grnpppI)_{J_3}\mid\mid (\grnpppp)_{K_1}\mid=\frac{\mid (\mcl{A}_{\ul{\gl}'''})_{\ti{J}}\mid}{(q^{\gl_{i_0}-u})}\mid(\grpppp)_K\mid.}

From Equations~\ref{Eq:One},~\ref{Eq:Two} we get that the first $i-1$ parts of $\ul{\gl'/I}$ are greater than $t_i=\gl_{i_0}$
and the remaining parts are less than $\gl_{i_0}$. Also in $\ul{\gn'/I}$ the first $i-1$ parts are greater than $t_i-1=\gl_{i_0}-1$ and the remaining parts are less than $\gl_{i_0}-1$. If $\max([J_3]_{(\gm)})=\{b,\gm\}$ with $\gm>\gl_{i_0}-1$ then $\max([\ti{J}]_{(\gm+2)})=\{b+1,\gm+2\}$. If $\max([J_3]_{(\gm)})$ is empty for $\gm>\gl_{i_0}-1$ then $\max([\ti{J}]_{(\gm+2)})=\{(\gm+1,\gm+2)\}$. If $\gm\leq \gl_{i_0}-1$ then  $\max([J_3]_{(\gm)})=\max([\ti{J}]_{(\gm)})$.
So there is a decrement by one that occurs in exactly $i-1$ parts of $J_3$ with respect to $\ti{J}$ because $\gm-b+1=(\gm+2)-(b+1)$. Now we observe that the first $i$ parts of $\ul{\gl}'$ are greater than equal to $\gl_{i_0}$ and the remaining parts of $\ul{\gl}'$ are less than $\gl_{i_0}$. For the ideals $K_1$ and $K$, there is decrement by one that occurs in those parts of $\ul{\gl}''$ which are greater than equal to $\gl_{i_0}$. So the total decrement is \equa{&(i-1)+\#(\text{parts of }\ul{\gl}'' \text{ greater than equal to }\gl_{i_0})=\\ &\#(\text{parts of }\ul{\gl}' \text{ greater than equal to }\gl_{i_0})-1+\#(\text{parts of }\ul{\gl}'' \text{ greater than equal to }\gl_{i_0})\\&=-1+\#(\text{parts of }\ul{\gl} \text{ greater than equal to }\gl_{i_0}).}
In $(\mcl{A}_{\ul{\gl}'''})_{\ti{J}}$ there is an extra component $\gp^u(\R/\gp^{\gl_{i_0}}\R)$ corresponding to $(u,\gl_{i_0})$.
Hence in this case from Equations~\ref{Eq:One},~\ref{Eq:Two}, Equation~\ref{Eq:SevenSeven} follows.

In Equation~\ref{Eq:Three}, a similar computation yields that 
\equan{EightEight}{\frac{\mid (\grnpppp)_{K_1}\mid}{q^{\gl_{i_0}-u-1}}=\frac{\mid (\grpppp)_{K}\mid}{q^{\#(\text{parts of }\ul{\gl}''\text{ greater than equal to }\gl_{i_0})}}}
as there is an extra component $\gp^u(\R/\gp^{\gl_{i_0}-u-1}\R)$ in $K_1$ since $\gl_{i_0}-1=t_i-1$ is not a part of $\ul{\gn}'$.

In Equation~\ref{Eq:Four} $\gl_{i_0}-1$ is a part of $\ul{\gl'/I}$. So a similar computation yields that 
\equan{NineNine}{\mid (\grnpppI)_{J_3}\mid=\frac{\mid(\mcl{A}_{\ul{\gl}'''})_{\ti{J}} \mid }{q^{i-1}q^{\gl_{i_0}-u}q^{\gl_{i_0}-u-1}}}
where $i$ is the number of parts of $\ul{\gl}'$ that are greater than equal to $\gl_{i_0}=t_i$. Here in this case there are extra components $\gp^u(\R/\gp^{\gl_{i_0}}\R)\oplus \gp^u(\R/\gp^{\gl_{i_0}-1}\R)$ in $(\mcl{A}_{\ul{\gl}'''})_{\ti{J}}$.

In Equation~\ref{Eq:Five}, $\gl_{i_0}+1$ is a part of $\ul{\gl'/I}$. So a similar computation yields that 
\equan{TenTen}{\mid (\grnpppI)_{J_3}\mid=\frac{\mid(\mcl{A}_{\ul{\gl}'''})_{\ti{J}} \mid }{q^{i-2}q^{\gl_{i_0}-u}q^{\gl_{i_0}-u}}}
where $i$ is the number of parts of $\ul{\gl}'$ that are greater than equal to $\gl_{i_0}=t_i$. Here in this case there are extra components $\gp^{u+1}(\R/\gp^{\gl_{i_0}+1}\R)\oplus \gp^u(\R/\gp^{\gl_{i_0}}\R)$ in $(\mcl{A}_{\ul{\gl}'''})_{\ti{J}}$.

Hence by combining Equations~\ref{Eq:NineNine},~\ref{Eq:TenTen} individually with Equation~\ref{Eq:EightEight} we get that Equation~\ref{Eq:SevenSeven} follows in the case where $t_i-1=\gl_{i_0}-1$ is not a part of $\ul{\gn}'$. This completes the proof of Equation~\ref{Eq:SixSix}.
\begin{center}
	*****************************************************************
\end{center} 
First we observe that for $J_4\in \mcl{J}(P_{\ul{\gn'/I_1}}),K_2\in \mcl{J}(P_{\ul{\gn}''})$, we have
\equa{&\frac{\mid X^{\ul{\gn}}_{I_1,J_4,K_2}\mid}{\ga^{\ul{\gn}}_{I_1,J_4,K_2}}=\\&\frac{\mid \grnppp\mid }{\mid \grnpppI\mid}\frac{1}{\mid(\grnpp)_{L_1}\mid}\bigg(\mid (\grnpppI)^*_{J_4}\mid \mid (\grnpppp)_{K_2}\mid\us{\os{\gl\in \ul{\gn}''}{\os{ \text{ such that there exists }}{(u,\gl)\in \max(K_2)\cap[J_4]_{\ul{\gn}''}}}}{\prod}(1-\frac1{q^{m(\gl)}})\bigg).}
Therefore we have 
\fo{10}{10}{
	\equa{&\us{\os{J_1\subseteq J_4\subseteq J_3,K_2\subseteq K_1}{\us{\max(K_2)\cap [J_4]_{\ul{\gn}''}\supseteq S_1}{[J_4\cup K_2]_{\ul{\gn}}=L_1}}}{\sum}\frac{\mid X^{\ul{\gn}}_{I_1,J_4,K_2}\mid}{\ga^{\ul{\gn}}_{I_1,J_4,K_2}}=\\&\frac{\mid \grnppp\mid }{\mid \grnpppI\mid}\frac{1}{\mid(\grnpp)_{L_1}\mid}\us{\os{J_1\subseteq J_4\subseteq J_3,K_2\subseteq K_1}{\us{\max(K_2)\cap [J_4]_{\ul{\gn}''}\supseteq S_1}{[J_4\cup K_2]_{\ul{\gn}}=L_1}}}{\sum}\bigg(\mid (\grnpppI)^*_{J_4}\mid \mid (\grnpppp)_{K_2}\mid\us{\os{\gl\in \ul{\gn}''}{\os{ \text{ such that there exists }}{(u,\gl)\in \max(K_2)\cap[J_4]_{\ul{\gn}''}}}}{\prod}(1-\frac1{q^{m(\gl)}})\bigg).}}
Using Proposition~\ref{prop:ExistenceofJK}, we have that, every $J_4\in \mcl{J}(P_{\ul{\gn'/I_1}}),K_2\in \mcl{J}(P_{\ul{\gn}''})$ appears in a unique such sum exactly once. Now we observe that summing all such equations over all possibilities for a fixed $I_1$ we get  
\equ{\bigg(\us{J\in \mcl{J}(P_{\ul{\gn^{'}/I_1}}),K\in \mcl{J}(P_{\ul{\gn}''})}{\sum}\frac{\mid X^{\ul{\gn}}_{I_1,J,K}\mid}{\ga^{\ul{\gn}}_{I_1,J,K}}\bigg).}

Also correspondingly summing the below expressions 
\fo{10}{10}{
	\equ{\frac{\mid \grppp\mid }{\mid \grpppI\mid}\frac{1}{\mid(\grpp)_{L}\mid}\frac{\mid (\mcl{A}_{\ul{\gl}'''})^*_{\ti{J}}\mid}{(q^{\gl_{i_0}-u}-q^{\gl_{i_0}-u-1})}\mid(\grpppp)_K\mid\us{\os{\gl\in \ul{\gl}''}{\os{ \text{ such that there exists }}{(u,\gl)\in \max(K)\cap[J]_{\ul{\gl}''}}}}{\prod}(1-\frac1{q^{m(\gl)}})}}
over all possibilities for a fixed $I\in \mcl{J}(P_{\ul{\gl}})$ such that $\max(I)\cap P_{(\gl_{i_0})}\neq \es$, we get 
\equ{\bigg(\us{\us{\max([J]_{(\gl_{i_0})})<\max([K]_{(\gl_{i_0})}),\max(K)\cap P_{(\gl_{i_0})}\neq\es}{J\in \mcl{J}(P_{\ul{\gl^{'}/I}}),K\in \mcl{J}(P_{\ul{\gl}''})}}{\sum}\frac{\mid X^{\ul{\gl}}_{I,J,K}\mid}{\ga^{\ul{\gl}}_{I,J,K}}\bigg).}
\begin{center}
	*****************************************************************
\end{center} 
This gives us, under the association $I\in \mcl{J}(P_{\ul{\gl}})$ such that $\max(I)\cap P_{(\gl_{i_0})}\neq \es$ with $\gc(I)=I_1\in \mcl{J}(P_{\ul{\gn}})$,
\equan{ElevenEleven}{\bigg(\us{\us{\max([J]_{(\gl_{i_0})})<\max([K]_{(\gl_{i_0})}),\max(K)\cap P_{(\gl_{i_0})}\neq\es}{J\in \mcl{J}(P_{\ul{\gl^{'}/I}}),K\in \mcl{J}(P_{\ul{\gl}''})}}{\sum}\frac{\mid X^{\ul{\gl}}_{I,J,K}\mid}{\ga^{\ul{\gl}}_{I,J,K}}\bigg)=\bigg(\us{J\in \mcl{J}(P_{\ul{\gn^{'}/I_1}}),K\in \mcl{J}(P_{\ul{\gn}''})}{\sum}\frac{\mid X^{\ul{\gn}}_{I_1,J,K}\mid}{\ga^{\ul{\gn}}_{I_1,J,K}}\bigg).}

So upon summation over $I\in \mcl{J}(P_{\ul{\gl}})$ such that $\max(I)\cap P_{(\gl_{i_0})}\neq \es$, we get the identity in Equation~\ref{Eq:FiveFive}. This proves Theorem~\ref{theorem:LambdaNu}.
\end{proof}

\section{\bf{The polynomial $n^0_{\ul{\gl}}(q)=\mid \autgp\bs\big( \{Height\ Zero\ Elements\}\times \grpp\big)\mid$ for a partition $\ul{\gl}$}}
We say an element $a\in \grpp$ is of height zero if the equation $a=px$ has no solution for $x\in \grpp$. We say an ideal $I\in \mcl{J}(P_{\ul{\gl}})$ is a height zero ideal if one of the co-ordinates of $e_I$ is $1$, that is, $v_i=0$ for some $1\leq i\leq s$ in~\ref{Eq:CanonicalForm}. Note in this case all elements of $(\grpp)^*_{I}$ are of height zero.
Let $\mcl{J}^0(P_{\ul{\gl}})\subseteq \mcl{J}(P_{\ul{\gl}})$ be the set of height zero ideals in $\mcl{J}(P_{\ul{\gl}})$. Let \equ{\hgrpp=\us{I\in \mcl{J}^0(P_{\ul{\gl}})}{\bigsqcup} (\grpp)^*_{I}.}
If $\ul{\gl}=\es$ then $\hgrpp=\es$. Define for $\ul{\gl}\neq \es$, \equ{n^0_{\ul{\gl}}(q)=\mid \autgp\bs\big( \hgrpp\times \grpp\big)\mid.}
If $\ul{\gl}=\es$ then define $n^0_{\ul{\gl}}(q)=1$ and not zero for the sake of some consistency which we will see in this section.
The following facts about $n^0_{\ul{\gl}}(q)$ can be deduced immediately.
\begin{enumerate}
	\item The number $n^0_{\ul{\gl}}(q)$ is a polynomial in $q$ with integer coefficients (C.~P.~Anil Kumar, Amritanshu Prasad~\cite{MR3261812}, Theorem 6.4, Page 21).
	\item The degree of $n^0_{\ul{\gl}}(q)$ is $\gl_1$ the highest part of $\ul{\gl}$. This can be deduced by adapting the proof of 
	C.~P.~Anil Kumar, Amritanshu Prasad~\cite{MR3261812}, Theorem 5.11, Page 17, using the maximal ideal generated by $P_{\ul{\gl}}$ which is a height zero principal ideal. The polynomial $n^0_{\ul{\gl}}(q)$ is monic as well.
	\item The polynomial $n^0_{\ul{\gl}}(q)=n^0_{\ul{\gl}^2}(q)$ where \equ{\ul{\gl}^2=(\gl_1^{\min(\gr_1,2)}>\gl_2^{\min(\gr_2,2)}>\ldots>\gl_k^{\min(\gr_k,2)})}
	(C.~P.~Anil Kumar, Amritanshu Prasad~\cite{MR3261812}, Corollary 4.5  and Remark 4.7 on Page 12). So, for the computation of the polynomial $n^0_{\ul{\gl}}(q)$, any part of the partition which repeats more than twice can be reduced to two. 
	\item An alternative expression for $n^0_{\ul{\gl}}(q)$ for $\ul{\gl}\neq \es$ is given as 
	\equan{Alternative}{n^0_{\ul{\gl}}(q)=\us{I\in \mcl{J}^0(P_{\ul{\gl}}) }{\sum}\mid (\autgp)_I\bs\grpp\mid= \us{I\in \mcl{J}^0(P_{\ul{\gl}}) }{\sum}\bigg( \us{J\in \mcl{J}(P_{\ul{\gl^{'}/I}}),K\in \mcl{J}(P_{\ul{\gl}''})}{\sum}\frac{\mid X^{\ul{\gl}}_{I,J,K}\mid}{\ga^{\ul{\gl}}_{I,J,K}}\bigg)}
	where $\ga^{\ul{\gl}}_{I,J,K}$ is cardinality of the valuative set $(\grppp)_{J\cup K}\oplus \big((\grpppp)^*_K+(\grpppp)_J\big)$ $\subseteq \grppp\oplus \grpppp=\grpp$ and $X^{\ul{\gl}}_{I,J,K}=\{(x',x'')\in \grpp\mid I(\ol{x}')=J,I(x'')=K\}$.
	\item We have $\mid\autgp\bs(\gp\grpp\times \gp\grpp)\mid=\mid\autmgp\bs(\grmpp\times \grmpp)\mid$ where $\ul{\gm}$ is the partition corresponding to the finite $\R$-module $\gp\grpp$.
	Therefore \equ{\mid\autgp\bs(\grpp\times \grpp)\mid=\mid\autmgp\bs(\grmpp\times \grmpp)\mid+\mid\autgp\bs(\hgrpp\times \grpp)\mid+\mid\autgp\bs(\gp\grpp\times \hgrpp)\mid.}
	We also have  \equa{&\mid\autgp\bs(\gp\grpp\times \hgrpp)\mid=\us{I\in \big(\mcl{J}(P_{\ul{\gl}})\bs \mcl{J}^0(P_{\ul{\gl}})\big) }{\sum}\mid (\autgp)_I\bs\hgrpp\mid\\
	&=\us{I\in \big(\mcl{J}(P_{\ul{\gl}})\bs \mcl{J}^0(P_{\ul{\gl}})\big)}{\sum}
	\bigg( \us{\os{\text{either }J\in \mcl{J}(P_{\ul{\gl^{'}/I}}) \text{ is of height zero or }}{K\in \mcl{J}(P_{\ul{\gl}''})\text{ is of height zero}}}{\sum}\frac{\mid X^{\ul{\gl}}_{I,J,K}\mid}{\ga^{\ul{\gl}}_{I,J,K}}\bigg).}
So we have 
\equan{ThirteenThirteen}{n_{\ul{\gl}}(q)&=n_{\ul{\gm}}(q)+n^0_{\ul{\gl}}(q)+\mid\autgp\bs(\gp\grpp\times \hgrpp)\mid\\&=n_{\ul{\gm}}(q)+n^0_{\ul{\gl}}(q)\\&+\us{I\in \big(\mcl{J}(P_{\ul{\gl}})\bs \mcl{J}^0(P_{\ul{\gl}})\big)}{\sum}
\bigg( \us{(J,K)\in \mcl{J}^0(P_{\ul{\gl^{'}/I}})\times \mcl{J}(P_{\ul{\gl}''}) \bigcup \mcl{J}(P_{\ul{\gl^{'}/I}})\times \mcl{J}^0(P_{\ul{\gl}''})}{\sum}\frac{\mid X^{\ul{\gl}}_{I,J,K}\mid}{\ga^{\ul{\gl}}_{I,J,K}}\bigg).}
\end{enumerate} 
\begin{remark}
Here below is the table of values of $n^0_{\ul{\gl}}(q)$ for partitions whose parts sum upto $6$. From the table, we observe that $n^0_{(\gl)}=q^{\gl}$. So it is always true that the polynomial $n_{\ul{\gl}}(q)$ can be expressed as a sum, in terms of $n^0_{\ul{\gm}}(q)$ for some partitions $\ul{\gm}=(\gm), \gm\geq 0$. We are not interested in the trivial manner of expression where we do not know apriori, the coefficients of the summands whether they are positive or negative or zero. 
\end{remark}
\begin{center}  
\begin{tabular}{|c|c|c|}
\hline\\
The Partition & The Height Zero Polynomial & The Orbit Polynomial\\	
$\ul{\gl}$ & $n^0_{\ul{\gl}}(q)$ & $n_{\ul{\gl}}(q)$\\
\hline\\
$\es$ & $1(\text{using definition})$ & $1$\\
$(1)$ & $q$ & $q+2$\\
$(1^2)$ & $q+1$ & $q+3$\\
$(2)$ & $q^2$ & $q^2+2q+2$\\
$(2>1)$ & $q^2+3q+1$ & $q^2+5q+5$\\
$(2>1^2)$ & $q^2+3q+2$ & $q^2+5q+6$\\
$(2^2)$ & $q^2+q+1$ & $q^2+3q+5$\\
$(2^2>1)$ & $q^2+4q+2$ & $q^2+6q+8$\\
$(2^2>1^2)$ & $q^2+4q+3$ & $q^2+6q+9$\\
$(3)$ & $q^3$ & $q^3+2q^2+2q+2$\\
$(3>1)$ & $q^3+3q^2+3q$ & $q^3+5q^2+7q+4$\\
$(3>1^2)$ & $q^3+3q^2+4q+2$ & $q^3+5q^2+8q+6$\\
$(3>2)$ & $q^3+3q^2+2q+1$ & $q^3+5q^2+10q+7$\\
$(3>2>1)$ & $q^3+6q^2+9q+3$ & $q^3+8q^2+19q+13$\\
$(3^2)$ & $q^3+q^2+q+1$ & $q^3+3q^2+5q+7$\\
$(4)$ & $q^4$ & $q^4+2q^3+2q^2+2q+2$\\
$(4>1)$ & $q^4+3q^3+3q^2+2q$ & $q^4+5q^3+7q^2+6q+4$\\
$(4>1^2)$ & $q^4+3q^3+4q^2+4q+2$ & $q^4+5q^3+8q^2+8q+6$\\
$(4>2)$ & $q^4+3q^3+4q^2+q$ & $q^4+5q^3+12q^2+11q+4$\\
$(5)$ & $q^5$ & $q^5+2q^4+2q^3+2q^2+2q+2$\\
$(5>1)$ & $q^5+3q^4+3q^3+2q^2+2q$ & $q^5+5q^4+7q^3+6q^2+6q+4$\\
$(6)$ & $q^6$ & $q^6+2q^5+2q^4+2q^3+2q^2+2q+2$\\
\hline	
\end{tabular}
\end{center}

Now we state a theorem, using which we can express the polynomial $n_{\ul{\gl}}(q)$ as a sum, in terms of $n^0_{\ul{\gm}}(q)$ for some suitable partitions $\ul{\gm}$ (not the trivial manner) where we know that the coefficients of the summands are positive. If $\ul{\gl}$ has distinct parts then parts of the suitable partitions $\ul{\gm}$ are also distinct. 
\begin{theorem}
\label{theorem:OrbitPolytoHeightZeroPoly}
Let $\ugl\in \Gl_0$ be a nonempty partition. 
\begin{enumerate}
	\item Suppose $\gl_k>1$. Let $\ul{\gm}=\big((\gl_1-1)^{\gr_1}>\cdots>(\gl_k-1)^{\gr_k}\big)$. Then we have 
	\equ{n_{\ul{\gl}}(q)=n^0_{\ul{\gl}}(q)+n_{\ul{\gm}}(q)+n^0_{\ul{\gm}}(q).}
	\item Suppose $\gl_k=1$. Let $\ul{\gm}=\big((\gl_1-1)^{\gr_1}>\cdots>(\gl_{k-1}-1)^{\gr_{k-1}}\big)$. Then we have 
	\equ{n_{\ul{\gl}}(q)=n^0_{\ul{\gl}}(q)+2n_{\ul{\gm}}(q).}
\end{enumerate}
In particular if $\ul{\gl}$ has distinct parts then $\ul{\gm}$ also has distinct parts in both cases.
\end{theorem}
We prove Theorem~\ref{theorem:OrbitPolytoHeightZeroPoly} after illustrating it in the following example.
\begin{example}
\equa{n_{(3^2>2^2)}(q)&=n^0_{(3^2>2^2)}(q)+n_{(2^2>1^2)}(q)+n^0_{(2^2>1^2)}(q)\\
&=n^0_{(3^2>2^2)}(q)+n^0_{(2^2>1^2)}(q)+2n_{(1^2)}(q)+n^0_{(2^2>1^2)}(q)\\
&=n^0_{(3^2>2^2)}(q)+2n^0_{(2^2>1^2)}(q)+2n_{(1^2)}(q)\\
&=n^0_{(3^2>2^2)}(q)+2n^0_{(2^2>1^2)}(q)+2(n^0_{(1^2)}(q)+2n_{\es}(q))\\
&=n^0_{(3^2>2^2)}(q)+2n^0_{(2^2>1^2)}(q)+2n^0_{(1^2)}(q)+4n^0_{\es}(q).}
This identity can be verified by the actual polynomials $n_{(3^2>2^2)}(q)=q^3+6q^2+14q+15,n^0_{(3^2>2^2)}(q)=q^3+4q^2+4q+3,n^0_{(2^2>1^2)}(q)=q^2+4q+3,n^0_{(1^2)}(q)=q+1,n_{\es}(q)=1=n^0_{\es}(q)$.
\end{example}
\begin{proof}[Proof of Theorem~\ref{theorem:OrbitPolytoHeightZeroPoly}]
To prove this theorem we use Equation~\ref{Eq:ThirteenThirteen}. To prove (1), it is enough to show that
\equ{n^0_{\ul{\gm}}(q)=\mid \autmgp\bs\big( \grmpp\times \hgrmpp\big)\mid=\mid\autgp\bs(\gp\grpp\times \hgrpp)\mid.}
We first show that 	$\mid\autgp\bs(\gp\grpp\times \hgrpp)\mid=\mid\autgp\bs(\gp\grpp\times \gp\hgrpp)\mid$. 
Let $x_1,y_1\in \grpp,x_2,y_2\in \hgrpp$. Suppose $g\in \autgp$ is such that $g\gp x_1=\gp y_1,gx_2=y_2$. Then we have 
$g \gp x_1=\gp y_1,g \gp x_2=\gp y_2$. Conversely suppose $g \gp x_1=\gp y_1, g \gp x_2=\gp y_2$. Let $z=gx_2-y_2$. Then we have $\gp z=0$.
So $z\in (\grpp)_J$ where $J\in \mcl{J}(P_{\ul{\gl}})$ is an ideal such that $(\grpp)_J=\us{i=1}{\os{k}{\oplus}}\gp^{\gl_i-1}(\R/\gp^{\gl_i}\R)^{\gr_i}\subs \grpp$. Note that $x_2$ is an element of  height zero in $\grpp$. Any height zero ideal, say $K$, (here $K$ is the ideal of $x_2$) in the big fundamental poset $P$ defined in Equation~\ref{Eq:FP} contains the ideal generated by $J$ in the big fundamental poset $P$. So there is an endomorphism $h\in \Hom(\grpp,\us{i=1}{\os{k}{\oplus}}\gp^{\gl_i-1}(\R/\gp^{\gl_i}\R)^{\gr_i})$ such that $hx_2=z$. Here we have $\gp h=0, (g-h)(x_2)=gx_2-hx_2=gx_2-z=y_2$. Now note that, since $\gl_i\geq \gl_k>1$ for $1\leq i\leq k$ we also have $h\mod \gp=0$. Hence we have $g\mod \gp=(g-h)\mod \gp$. Now any endomorphism $k\in \End(\grpp)$ is invertible if and only $k\mod \gp$ is invertible. So the endomorphism $g-h\in \End(\grpp)$ is invertible. Moreover $(g-h)\gp x_1=gx_1-\gp hx_1=gx_1=y_1$.  This proves that $(\gp x_1,x_2)$ and $(\gp y_1,y_2)$ are in the same $\autgp$ orbit if and only if $(\gp x_1,\gp x_2),(\gp y_1,\gp y_2)$ are in the same $\autgp$ orbit. Hence we have $\mid\autgp\bs(\gp\grpp\times \hgrpp)\mid=\mid\autgp\bs(\gp\grpp\times \gp\hgrpp)\mid$. Since $\gl_k>1$, the set $\gp\hgrpp$ is precisely the set of height zero elements in $\gp\grpp\cong \grmpp$, that is, $\gp\hgrpp\cong \hgrmpp$. Now we have $\mid\autgp\bs(\gp\grpp\times \gp\hgrpp)\mid=\mid \autmgp\bs\big( \grmpp\times \hgrmpp\big)\mid$ because the map $\autgp\lra \autmgp$ which takes $g$ to $g_{\mid_{\gp\grpp}}$ is surjective. This proves (1).

Now we prove (2). Here $\gl_k=1$. In this case we have $\gp\hgrpp=\gp\grpp\cong \grmpp$. Using Equation~\ref{Eq:ThirteenThirteen}, it is enough to show that
\equ{n_{\ul{\gm}}(q)=\mid \autmgp\bs\big( \grmpp\times \grmpp\big)\mid=\mid\autgp\bs(\gp\grpp\times \hgrpp)\mid.}

The proof is similar to that of the previous case except a bit longer. Let $x_1,y_1\in \grpp,x_2,y_2\in \hgrpp$. Suppose $g\in \autgp$ is such that $g\gp x_1=\gp y_1,gx_2=y_2$. Then we have $g \gp x_1=\gp y_1,g \gp x_2=\gp y_2$. Conversely suppose $g \gp x_1=\gp y_1, g \gp x_2=\gp y_2$. So $\gp x_2,\gp y_2$ are of the same height in $\gp\grpp\cong \grmpp$. Let $t$ stand for transpose. Let the column vectors $x_2=(x^t_{21},x^t_{22},\cdots,x^t_{2k})^t,y_2=(y^t_{21},y^t_{22},\cdots,y^t_{2k})^t\in \us{i=1}{\os{k}{\oplus}}(\R/\gp^{\gl_i}\R)^{\gr_i}=\grpp$
with the column vectors $x_{2i}=(x^1_{2i},x^2_{2i}\cdots,x^{\gr_i}_{2i})^t,y_{2i}=(y^1_{2i},y^2_{2i}\cdots,y^{\gr_i}_{2i})^t\in (\R/\gp^{\gl_i}\R)^{\gr_i},1\leq i\leq k$.  Let the column vector $gx_2=((gx)^t_{21},(gx)^t_{22},\cdots,(gx)^t_{2k})^t$ $\in \us{i=1}{\os{k}{\oplus}}(\R/\gp^{\gl_i}\R)^{\gr_i}=\grpp$
with the column vectors $(gx)_{2i}=((gx)^1_{2i},(gx)^2_{2i}\cdots,(gx)^{\gr_i}_{2i})^t\in (\R/\gp^{\gl_i}\R)^{\gr_i}$, $1\leq i\leq k$.
Let the column vector $z=gx_2-y_2=(z^t_1,z^t_2,\cdots,z^t_k)^t\in \us{i=1}{\os{k}{\oplus}}\gp^{\gl_i-1}(\R/\gp^{\gl_i}\R)^{\gr_i}$.

\begin{claim}
There exists $g_1\in \autgp$ such that $g_1\gp x_1=\gp y_1,g_1\gp x_2=\gp y_2$ and if the column vector $w=g_1x_2-y_2=(w^t_1,w^t_2,\cdots,w^t_k)^t\in \us{i=1}{\os{k}{\oplus}}\gp^{\gl_i-1}(\R/\gp^{\gl_i}\R)^{\gr_i}$ then $w_k=0$.
\end{claim}
\begin{proof}[Proof of Claim]
If $z_k$ is already zero then we take $g=g_1$. If $z_k\neq 0$ then there are two cases to consider.

We consider the first case.
Let $x_{2i_0}$ contain a unit in one of its coordinates for some $1\leq i_0\leq k-1$. We choose a homomorphism $e=[e^{ij}]_{k\times k}:\grpp\lra \grpp$ where $e^{ij}=[e_{mn}^{ij}]_{\gr_i\times \gr_j}:(\R/\gp^{\gl_j}\R)^{\gr_j}\lra (\R/\gp^{\gl_i}\R)^{\gr_i}$ such that $e^{ij}$ is nonzero only when $i=k,j=i_0$ and $e^{ki_0}(x_{2i_0})=z_k$. We observe that $\gp e=0$ because $\gp e^{ki_0}=0$. We have $i_0\neq k$ and the diagonal blocks of $g-e$ are all invertible modulo $\gp$. Hence $g-e$ is invertible. Moreover we have 
$(g-e)\gp x_1=g\gp x_1=\gp y_1,(g-e)\gp x_2=g\gp x_2=\gp y_2, (g-e)x_2-y_2=(z^t_1,z^t_2,\cdots,z^t_{k-1},\ul{0}^t)^t$. Hence choose $g_1=g-e$.

We consider the second case.
Let $x_{2i}$ does not contain a unit in any of its coordinates for all $1\leq i_0\leq k-1$. So $x_{2k}$ must contain a unit in one if its coordinates. Hence $(gx)_{2k}$ must also contain a unit in one of its coordinates. As $g\gp x_2=\gp y_2, \gp x_2,\gp y_2$ are of the same height in $\gp\grpp\cong \grmpp$ and therefore $y_{2i}$ also does not contain a unit in any of its coordinates for all $1\leq i_0\leq k-1$. So $y_{2k}$ must contain a unit in one if its coordinates. Hence we choose an automorphism $e=[e^{ij}]_{k\times k}:\grpp\lra \grpp$ where $e^{ij}=[e_{mn}^{ij}]_{\gr_i\times \gr_j}:(\R/\gp^{\gl_j}\R)^{\gr_j}\lra (\R/\gp^{\gl_i}\R)^{\gr_i}$ such that $e^{ij}$ is the identity block for $1\leq i=j\leq k-1$, $e^{ij}$ is the zero block for $1\leq i\neq j\leq k$ and $e^{kk}: (\R/\gp^{\gl_k}\R)^{\gr_k} \lra  (\R/\gp^{\gl_k}\R)^{\gr_k}$ is an automorphism such that $e^{kk}((gx)_{2k})=y_{2k}$. Now we have $eg(x_2)-y_2=(z^t_1,z^t_2,\cdots,z^t_{k-1},\ul{0}^t)^t, eg(\gp x_1)=\gp y_1,eg(\gp x_2)=\gp y_2$. Hence choose $g_1=eg$ which is an automophism.

This proves the claim.
\end{proof}
Continuing with the proof of (2), we assume without loss of generality that $g\gp x_1=\gp y_1,g \gp x_2=\gp y_2$ and  $gx_2-y_2=z=(z^t_1,z^t_2,\cdots,z^t_{k-1},z_k^t)^t$ with $z_k=0$ and $z\in (\grpp)_J=\us{i=1}{\os{k}{\oplus}}\gp^{\gl_i-1}(\R/\gp^{\gl_i}\R)^{\gr_i}\subs \grpp$ where $J$ is defined by its association to $(\grpp)_J$. Note that $x_2$ is an element of  height zero in $\grpp$. Any height zero ideal, say $K$, (here $K$ is the ideal of $x_2$) in the big fundamental poset $P$ defined in Equation~\ref{Eq:FP} contains the ideal generated by $J$ in the big fundamental poset $P$. So there is an endomorphism $h\in \Hom(\grpp,\us{i=1}{\os{k}{\oplus}}\gp^{\gl_i-1}(\R/\gp^{\gl_i}\R)^{\gr_i})$ such that $hx_2=z$ and for this $h$ we have $\gp h=0$. Now we left multiply $h$ by a homomorphism $f=[f^{ij}]_{k\times k}:\grpp\lra \grpp$ where $f^{ij}=[f_{mn}^{ij}]_{\gr_i\times \gr_j}:(\R/\gp^{\gl_j}\R)^{\gr_j}\lra (\R/\gp^{\gl_i}\R)^{\gr_i}$ such that $f^{ij}$ is the identity block for $1\leq i=j\leq k-1$, $f^{ij}$ is the zero block for $1\leq i\neq j\leq k$ and $f^{kk}: (\R/\gp^{\gl_k}\R)^{\gr_k} \lra  (\R/\gp^{\gl_k}\R)^{\gr_k}$ is also zero. Then we have $fh(x_2)=z$, the last row of blocks in $fh$ are all zero blocks, the remaining blocks of $fh$ are the same as that of $h$. So $\gp fh=0$ and also $fh \mod \gp = 0$. Hence we have $g\mod \gp=(g-fh)\mod \gp$. So $g-fh$ is invertible. Also $(g-fh)(\gp x_1)=g \gp x_1=\gp y_1,(g-fh)(x_2)-y_2=z-fh(x_2)=z-z=0$. 

This proves that $(\gp x_1,x_2)$ and $(\gp y_1,y_2)$ are in the same $\autgp$ orbit if and only if $(\gp x_1,\gp x_2),(\gp y_1,\gp y_2)$ are in the same $\autgp$ orbit. Hence we have $\mid\autgp\bs(\gp\grpp\times \hgrpp)\mid=\mid\autgp\bs(\gp\grpp\times \gp\hgrpp)\mid=\mid\autgp\bs(\gp\grpp\times \gp\grpp)\mid=\mid\autmgp\bs(\grmpp\times \grmpp)\mid=n_{\ul{\gm}}(q)$. This proves (2).

We have completed the proof of Theorem~\ref{theorem:OrbitPolytoHeightZeroPoly}.
\end{proof}
Now we state a theorem, using which we can express the polynomial $n^0_{\ul{\gl}}(q)$ for a partition $\ul{\gl}$ with a repeated part as a sum, in terms of the polynomials $n^0_{\ul{\gm}}(q)$ for some suitable partitions $\ul{\gm}$ which has all their parts distinct (not the trivial manner) where we know that the coefficients of the summands are positive.
\begin{theorem}
\label{theorem:HeightZeroPolyDistinctParts}
Let $\ugl\in \Gl_0$ be a partition such that $\gr_{i_0}=2$ for some $1\leq i_0\leq k$. Let $\ul{\gm} =\big(\gl_1^{\gr_1}>\gl_2^{\gr_2}>\ldots>\gl_{i_0-1}^{\gr_{i_0-1}}>\gl_{i_0}^{\gr_{i_0}-1}=\gl_{i_0}^1>\gl_{i_0+1}^{\gr_{i_0+1}}>\ldots>\gl_{k-1}^{\gr_{k-1}}>\gl_k^{\gr_k}\big)$.
Let $\ul{\gn}=\big((\gl_1-2)^{\gr_1}>(\gl_2-2)^{\gr_2}>\ldots>(\gl_{i_0-1}-2)^{\gr_{i_0-1}}\geq (\gl_{i_0}-1)^{\gr_{i_0}}\geq \gl_{i_0+1}^{\gr_{i_0+1}}>\ldots>\gl_{k-1}^{\gr_{k-1}}>\gl_k^{\gr_k}\big)$.
\begin{enumerate}
	\item Suppose $\gl_{i_0}>1$. Then we have 
	\equ{n^0_{\ul{\gl}}(q)=n^0_{\ul{\gm}}(q)+n^0_{\ul{\gn}}(q).}
	\item Suppose $\gl_{i_0}=1$ which implies in particular $i_0=k$. Then we have
	\equ{n^0_{\ul{\gl}}(q)=n^0_{\ul{\gm}}(q)+n_{\ul{\gn}}(q).}
\end{enumerate}
\end{theorem}
We prove Theorem~\ref{theorem:HeightZeroPolyDistinctParts} after illustrating it in the following example.
\begin{example}
\equa{n^0_{(3^2>2^2)}(q)&=n^0_{(3>2^2)}(q)+n^0_{(2^2)}(q)\\
&=n^0_{(3>2)}(q)+n^0_{(1^2)(q)}+n^0_{(2)}(q)+n^0_{(1^2)(q)}\\
&=n^0_{(3>2)}(q)+n^0_{(2)}(q)+2n^0_{(1^2)(q)}\\
&=n^0_{(3>2)}(q)+n^0_{(2)}(q)+2(n^0_{(1)(q)}+n_{\es}(q))\\
&=n^0_{(3>2)}(q)+n^0_{(2)}(q)+2n^0_{(1)}(q)+2n^0_{\es}(q).}
This identity can be verified by the actual polynomials $n^0_{(3^2>2^2)}(q)=q^3+4q^2+4q+3,n^0_{(3>2)}(q)=q^3+3q^2+2q+1 ,n^0_{(2)}(q)=q^2,n^0_{(1)}(q)=q,n_{\es}(q)=1=n^0_{\es}(q)$.
\end{example}
\begin{proof}[Proof of Theorem~\ref{theorem:HeightZeroPolyDistinctParts}]
To prove this theorem, we use the proof of Theorem~\ref{theorem:RepeatedPartCase}. 	We prove (1) and (2) simultaneously. We use the alternative expressions for $n^0_{\ul{\gl}}(q),n^0_{\ul{\gm}}(q),n^0_{\ul{\gn}}(q)$ as defined in Equation~\ref{Eq:Alternative}. 

In Equation~\ref{Eq:TwoPointFive}, if we sum over only height zero ideals $I\in \mcl{J}^0(P_{\ul{\gl}})=\mcl{J}^0(P_{\ul{\gm}})$ such that $\max(I)\cap P_{(\gl_{i_0})}=\es$, then we get 
\equan{FourteenFourteen}{\us{I\in \mcl{J}^0(P_{\ul{\gl}}),\max(I)\cap P_{(\gl_{i_0})}=\es }{\sum}&\bigg(\us{J\in \mcl{J}(P_{\ul{\gl^{'}/I}}),K\in \mcl{J}(P_{\ul{\gl}''})}{\sum}\frac{\mid X^{\ul{\gl}}_{I,J,K}\mid}{\ga^{\ul{\gl}}_{I,J,K}}\bigg)=\\ &\us{I\in \mcl{J}^0(P_{\ul{\gm}}), \max(I)\cap P_{(\gl_{i_0})}=\es }{\sum}\bigg(\us{J\in \mcl{J}(P_{\ul{\gm^{'}/I}}),K\in \mcl{J}(P_{\ul{\gm}''})}{\sum}\frac{\mid X^{\ul{\gm}}_{I,J,K}\mid}{\ga^{\ul{\gm}}_{I,J,K}}\bigg).} 

In Equation~\ref{Eq:FourPointFive}, if we sum over only height zero ideals $I\in \mcl{J}^0(P_{\ul{\gl}})=\mcl{J}^0(P_{\ul{\gm}})$ such that $\max(I)\cap P_{(\gl_{i_0})}\neq \es$, then we get

\equan{FifteenFifteen}{\us{I\in \mcl{J}^0(P_{\ul{\gl}}), \max(I)\cap P_{(\gl_{i_0})}\neq \es }{\sum}&\bigg(\us{\us{\us{\max([J]_{(\gl_{i_0})})<\max([K]_{(\gl_{i_0})}),\max(K)\cap P_{(\gl_{i_0})}=\es}{\max([J]_{(\gl_{i_0})})\geq \max([K]_{(\gl_{i_0})}) \text{ or }}}{J\in \mcl{J}(P_{\ul{\gl^{'}/I}}),K\in \mcl{J}(P_{\ul{\gl}''})}}{\sum}\frac{\mid X^{\ul{\gl}}_{I,J,K}\mid}{\ga^{\ul{\gl}}_{I,J,K}}\bigg)\\ &=\us{\os{I\in \mcl{J}^0(P_{\ul{\gm}})}{\max(I)\cap P_{(\gl_{i_0})}\neq \es} }{\sum}\bigg(\us{J\in \mcl{J}(P_{\ul{\gm^{'}/I}}),K\in \mcl{J}(P_{\ul{\gm}''})}{\sum}\frac{\mid X^{\ul{\gm}}_{I,J,K}\mid}{\ga^{\ul{\gm}}_{I,J,K}}\bigg).}

In Theorem~\ref{theorem:LatticeIso}(2), we have defined a lattice isomorphism $\gc:\mcl{M}=\{I\in \mcl{J}(P_{\ul{\gl}})\mid \max(I)\cap P_{(\gl_{i_0})}\neq \es\} \lra \mcl{J}(P_{\ul{\gn}})$. If $\gl_{i_0}=1$, then $\mcl{M}\subseteq \mcl{J}^0(P_{\ul{\gl}})$, since in this case, any ideal in $\mcl{M}$ is of height zero. If $\gl_{i_0}>1$ then the we have a restricted lattice isomorphism $\gc$ given as  
\equ{\gc_{\mid_{\mcl{M}\cap \mcl{J}^0(P_{\ul{\gl}})}}:\mcl{M}\cap \mcl{J}^0(P_{\ul{\gl}}) \lra \mcl{J}^0(P_{\ul{\gn}})}
onto the set $\mcl{J}^0(P_{\ul{\gn}})$ of height zero ideals. 

Hence in Equations~[\ref{Eq:ElevenEleven},\ref{Eq:FiveFive}], if we sum over only height zero ideals $I\in \mcl{J}^0(P_{\ul{\gl}})$ such that $\max(I)\cap P_{(\gl_{i_0})}\neq \es$, then we get for 
$\gl_{i_0}>1$,
\equan{SixteenSixteen}{\us{I\in \mcl{J}^0(P_{\ul{\gl}}), \max(I)\cap P_{(\gl_{i_0})}\neq \es }{\sum}&\bigg(\us{\us{\max([J]_{(\gl_{i_0})})<\max([K]_{(\gl_{i_0})}),\max(K)\cap P_{(\gl_{i_0})}\neq\es}{J\in \mcl{J}(P_{\ul{\gl^{'}/I}}),K\in \mcl{J}(P_{\ul{\gl}''})}}{\sum}\frac{\mid X^{\ul{\gl}}_{I,J,K}\mid}{\ga^{\ul{\gl}}_{I,J,K}}\bigg)\\&=\us{I_1\in \mcl{J}^0(P_{\ul{\gn}})}{\sum}\bigg(\us{J\in \mcl{J}(P_{\ul{\gn^{'}/I_1}}),K\in \mcl{J}(P_{\ul{\gn}''})}{\sum}\frac{\mid X^{\ul{\gn}}_{I_1,J,K}\mid}{\ga^{\ul{\gn}}_{I_1,J,K}}\bigg)}
and for $\gl_{i_0}=1$,
\equan{SeventeenSeventeen}{\us{I\in \mcl{J}(P_{\ul{\gl}}), \max(I)\cap P_{(1)}\neq \es }{\sum}&\bigg(\us{\us{\max([J]_{(1)})<\max([K]_{(1)}),\max(K)\cap P_{(1)}\neq\es}{J\in \mcl{J}(P_{\ul{\gl^{'}/I}}),K\in \mcl{J}(P_{\ul{\gl}''})}}{\sum}\frac{\mid X^{\ul{\gl}}_{I,J,K}\mid}{\ga^{\ul{\gl}}_{I,J,K}}\bigg)\\&=\us{I_1\in \mcl{J}(P_{\ul{\gn}})}{\sum}\bigg(\us{J\in \mcl{J}(P_{\ul{\gn^{'}/I_1}}),K\in \mcl{J}(P_{\ul{\gn}''})}{\sum}\frac{\mid X^{\ul{\gn}}_{I_1,J,K}\mid}{\ga^{\ul{\gn}}_{I_1,J,K}}\bigg).}

So upon summing Equations~\ref{Eq:FourteenFourteen},~\ref{Eq:FifteenFifteen},~\ref{Eq:SixteenSixteen}, we obtain (1) and upon summing Equations~\ref{Eq:FourteenFourteen},~\ref{Eq:FifteenFifteen},~\ref{Eq:SeventeenSeventeen}, we obtain (2). This proves Theorem~\ref{theorem:HeightZeroPolyDistinctParts}. 
\end{proof}
\begin{cor}
Let $\ugl\in \Gl_0$ be a partition. Then the polynomial $n_{\ul{\gl}}(q)$ is expressible as the sum of $n^0_{\ul{\gm}}(q)$ for partitions $\ul{\gm}$ which has all their parts distinct and $ \ul{\gm} \subseteq \ul{\gl}$ (not the trivial manner) where we know that the coefficients of the summands are positive.	
\end{cor}
\begin{proof}
This corollary follows from Theorem~\ref{theorem:OrbitPolytoHeightZeroPoly} and Theorem~\ref{theorem:HeightZeroPolyDistinctParts}.
\end{proof}
We illustrate the above corollary with an example.
\begin{example}
	\equa{n_{(2^2>1^2)}(q)&=n^0_{(2^2>1^2)}(q)+2n_{(1^2)}(q)\\
		&=n^0_{(2^2>1^2)}(q)+2(n^0_{(1^2)}(q)+2n_{\es}(q))\\
		&=(n^0_{(2>1^2)}(q)+n^0_{(1^2)}(q))+2n^0_{(1^2)}(q)+4n^0_{\es}(q)\\
		&=n^0_{(2>1^2)}(q)+3n^0_{(1^2)}(q)+4n^0_{\es}(q)\\
		&=(n^0_{(2>1)}(q)+n_{\es}(q))+3(n^0_{(1)}(q)+n_{\es}(q))+4n^0_{\es}(q)\\
		&=n^0_{(2>1)}(q)+3n^0_{(1)}(q)+8n^0_{\es}(q).}
	This identity can be verified by the actual polynomials $n_{(2^2>1^2)}(q)=q^2+6q+9,n^0_{(2>1)}(q)=q^2+3q+1,n^0_{(1)}(q)=q,n_{\es}(q)=1=n^0_{\es}(q)$.
\end{example}
Now we state a theorem using which we can express $n^0_{\ul{\gl}}(q)$ for a partition $\ul{\gl}$ which has all its parts distinct as a sum, in terms of $n^0_{\ul{\gm}}(q)$ for partitions $\ul{\gm}$ which has all its parts distinct  such that 
$\ul{\gm} \subsetneq \ul{\gl}$ with coefficients of the summands being polynomials in $q$ with nonnegative integer coefficients.
\begin{theorem}
\label{theorem:HeightZeroPartitiontoHeightZeroSmallerPartitions}
Let $\ul{\gl}=(\gl_1>\gl_2>\cdots>\gl_k)\in \Gl_0$ be a nonempty partition with distinct parts. Let $\ul{\gm}=(\gl_1-1>\gl_2-1>\cdots>\gl_k-1), \ul{\gn}=(\gl_1-\gl_k>\gl_2-\gl_k>\cdots>\gl_{k-1}-\gl_k),\ul{\gd}=(\gl_1-\gl_k-1>\gl_2-\gl_k-1>\cdots>\gl_{k-1}-\gl_k-1)$. 
\begin{enumerate}
	\item Suppose $k\geq 2$ and $\gl_{k-1}-\gl_k>1$. 
	Then we have \equa{n^0_{\ul{\gl}}(q)&=qn^0_{\ul{\gm}}(q)+n^0_{\ul{\gn}}(q)+n^0_{\ul{\gd}}(q) \text{ if }\gl_k>1\\
					   n^0_{\ul{\gl}}(q)&=qn_{\ul{\gm}}(q)+n^0_{\ul{\gn}}(q)+n^0_{\ul{\gd}}(q) \text{ if }\gl_k=1 (\text{ here }\ul{\gm}=\ul{\gn}).}
	\item Suppose $k\geq 2$ and $\gl_{k-1}-\gl_k=1$.  
	Then we have \equa{n^0_{\ul{\gl}}(q)&=qn^0_{\ul{\gm}}(q)+n^0_{\ul{\gn}}(q)+n_{\ul{\gd}}(q) \text{ if }\gl_k>1\\
					   n^0_{\ul{\gl}}(q)&=qn_{\ul{\gm}}(q)+n^0_{\ul{\gn}}(q)+n_{\ul{\gd}}(q) \text{ if }\gl_k=1(\text{ here }\ul{\gm}=\ul{\gn}).} 
	\item Suppose $k=1,\gl_1=\gl_k>1$. Let $\ul{\gm}=(\gl_1-1)$.
	Then we have \equ{n^0_{\ul{\gl}}(q)=qn^0_{\ul{\gm}}(q).} 
\end{enumerate}
\end{theorem}
We prove Theorem~\ref{theorem:HeightZeroPartitiontoHeightZeroSmallerPartitions} after illustrating it in the following example and proving Theorems~\ref{theorem:FirstSummand},~\ref{theorem:SecondSummandFirstPart}.
\begin{example}
\begin{enumerate}
	\item Examples for Theorem~\ref{theorem:HeightZeroPartitiontoHeightZeroSmallerPartitions}(1), $\gl_k>1,\gl_{k-1}-\gl_k>1$: 
	\begin{itemize}
	\item Let $\ul{\gl}=(4>2)$.
	We have $n^0_{(4>2)}(q)=qn^0_{(3>1)}(q)+n^0_{(2)}(q)+n^0_{(1)}(q)$. This identity holds as 
	$n^0_{(4>2)}(q)=q^4+3q^3+4q^2+q,n^0_{(3>1)}(q)=q^3+3q^2+3q,n^0_{(2)}(q)=q^2,n^0_{(1)}(q)=q$.
	\item Let $\ul{\gl}=(5>4>2)$. We $n^0_{(5>4>2)}(q)=qn^0_{(4>3>1)}(q)+n^0_{(3>2)}(q)+n^0_{(2>1)}(q)$.  This identity holds as 
	$n^0_{(5>4>2)}(q)=q^5+6q^4+15q^3+16q^2+7q+2,n^0_{(4>3>1)}(q)=q^4+6q^3+14q^2+12q+2,n^0_{(3>2)}(q)=q^3+3q^2+2q+1,n^0_{(2>1)}(q)=q^2+3q+1$.
	\item Let $\ul{\gl}=(6>4>2)$. We $n^0_{(6>4>2)}(q)=qn^0_{(5>3>1)}(q)+n^0_{(4>2)}(q)+n^0_{(3>1)}(q)$.  This identity holds as 
	$n^0_{(6>4>2)}(q)=q^6+6q^5+17q^4+22q^3+15q^2+4q,n^0_{(5>3>1)}(q)=q^5+6q^4+16q^3+18q^2+8q,n^0_{(4>2)}(q)=q^4+3q^3+4q^2+q,n^0_{(3>1)}(q)=q^3+3q^2+3q$.
	\end{itemize}
	\item Examples for Theorem~\ref{theorem:HeightZeroPartitiontoHeightZeroSmallerPartitions}(1), $\gl_k=1,\gl_{k-1}-\gl_k>1$: 
	\begin{itemize}
	\item Let $\ul{\gl}=(3>1)$.
	We have $n^0_{(3>1)}(q)=qn_{(2)}(q)+n^0_{(2)}(q)+n^0_{(1)}(q)$. This identity holds as 
	$n^0_{(3>1)}(q)=q^3+3q^2+3q, n_{(2)}(q)=q^2+2q+2,n^0_{(2)}(q)=q^2,n^0_{(1)}(q)=q$.
	\item Let $\ul{\gl}=(4>3>1)$.
	We have $n^0_{(4>3>1)}(q)=qn_{(3>2)}(q)+n^0_{(3>2)}(q)+n^0_{(2>1)}(q)$. This identity holds as 
	$n^0_{(4>3>1)}(q)=q^4+6q^3+14q^2+12q+2$, $n_{(3>2)}(q)=,q^3+5q^2+10q+7,n^0_{(3>2)}(q)=q^3+3q^2+2q+1,n^0_{(2>1)}(q)=q^2+3q+1$.
	\item Let $\ul{\gl}=(5>3>1)$.
	We have $n^0_{(5>3>1)}(q)=qn_{(4>2)}(q)+n^0_{(4>2)}(q)+n^0_{(3>1)}(q)$. This identity holds as 
	$n^0_{(5>3>1)}(q)=q^5+6q^4+16q^3+18q^2+8q, n_{(4>2)}(q)=q^4+5q^3+12q^2+11q+4, n^0_{(4>2)}(q)=q^4+3q^3+4q^2+q, n^0_{(3>1)}(q)=q^3+3q^2+3q$.
	\end{itemize}
	\item Examples for Theorem~\ref{theorem:HeightZeroPartitiontoHeightZeroSmallerPartitions}(2), $\gl_k>1,\gl_{k-1}-\gl_k=1$:  
	\begin{itemize}
	\item	Let $\ul{\gl}=(3>2)$.
	We have $n^0_{(3>2)}(q)=qn^0_{(2>1)}(q)+n^0_{(1)}(q)+n_{\es}(q)$. This identity holds as 
	$n^0_{(3>2)}(q)=q^3+3q^2+2q+1$, $n^0_{(2>1)}(q)=q^2+3q+1,n^0_{(1)}(q)=q,n_{\es}(q)=1$.
	\item Let $\ul{\gl}=(4>3>2)$.
	We have $n^0_{(4>3>2)}(q)=qn^0_{(3>2>1)}(q)+n^0_{(2>1)}(q)+n_{(1)}(q)$. This identity holds as 
	$n^0_{(4>3>2)}(q)=q^4+6q^3+10q^2+7q+3,n^0_{(3>2>1)}(q)=q^3+6q^2+9q+3,n^0_{(2>1)}(q)=q^2+3q+1,n_{(1)}(q)=q+2$.
	\item Let $\ul{\gl}=(5>3>2)$.
	We have $n^0_{(5>3>2)}(q)=qn^0_{(4>2>1)}(q)+n^0_{(3>1)}(q)+n_{(2)}(q)$. This identity holds as 
	$n^0_{(5>3>2)}(q)=q^5+6q^4+12q^3+13q^2+7q+2,n^0_{(4>2>1)}(q)=q^4+6q^3+11q^2+9q+2,n^0_{(3>1)}(q)$ $=q^3+3q^2+3q,n_{(2)}(q)=q^2+2q+2$.
	\end{itemize}
	\item Examples for Theorem~\ref{theorem:HeightZeroPartitiontoHeightZeroSmallerPartitions}(2), $\gl_k=1,\gl_{k-1}-\gl_k=1$:  
\begin{itemize}
	\item Let $\ul{\gl}=(2>1)$.
	We have $n^0_{(2>1)}(q)=qn_{(1)}(q)+n^0_{(1)}(q)+n_{\es}(q)$. This identity holds as 
	$n^0_{(2>1)}(q)=q^2+3q+1,n_{(1)}(q)=q+2,n^0_{(1)}(q)=q,n_{\es}(q)=1$.
	\item Let $\ul{\gl}=(3>2>1)$.
	We have $n^0_{(3>2>1)}(q)=qn_{(2>1)}(q)+n^0_{(2>1)}(q)+n_{(1)}(q)$. This identity holds as 
	$n^0_{(3>2>1)}(q)=q^3+6q^2+9q+3$, $n_{(2>1)}(q)=q^2+5q+5,n^0_{(2>1)}(q)=q^2+3q+1,n_{(1)}(q)=q+2$.
	\item Let $\ul{\gl}=(4>2>1)$.
	We have $n^0_{(4>2>1)}(q)=qn_{(3>1)}(q)+n^0_{(3>1)}(q)+n_{(2)}(q)$. This identity holds as 
	$n^0_{(4>2>1)}(q)=q^4+6q^3+11q^2+9q+2,n_{(3>1)}(q)=q^3+5q^2+7q+4, n^0_{(3>1)}(q)=q^3+3q^2+3q,n_{(2)}(q)=q^2+2q+2$.
\end{itemize}
	\item Example for Theorem~\ref{theorem:HeightZeroPartitiontoHeightZeroSmallerPartitions}(3): Let $\ul{\gl}=(2)$.
	We have $n^0_{(2)}(q)=qn^0_{(1)}(q)$. This identity holds as 
	$n^0_{(2)}(q)=q^2,n^0_{(1)}(q)=q$.
\end{enumerate}
\end{example}
We prove Theorem~\ref{theorem:HeightZeroPartitiontoHeightZeroSmallerPartitions} after proving the following theorem.
\begin{theorem}
\label{theorem:FirstSummand}
Let $\ul{\gl}=(\gl_1>\gl_2>\cdots>\gl_k)\in \Gl_0$ be a nonempty partition with distinct parts. Let $\ul{\gm}=(\gl_1-1>\gl_2-1>\cdots>\gl_k-1)$. Let $I\in \mcl{J}(P_{\ul{\gl}}),J\in \mcl{J}(P_{\ul{\gl'/I}}),K\in \mcl{J}(P_{\ul{\gl}''})$.
Define $X^{\ul{\gl}}_{I,J,K}=\{(x',x'')\in \grpp\mid I(\ol{x}')=J,I(x'')=K\}$  and
$\ga^{\ul{\gl}}_{I,J,K}=\mid(\grppp)_{J\cup K}\oplus\big((\grpppp)^*_K+(\grpppp)_J\big)\mid$.
Let 
\equ{X^{\ul{\gl}}_I=\us{J\in \big(\mcl{J}(P_{\ul{\gl^{'}/I}})\bs \mcl{J}^0(P_{\ul{\gl^{'}/I}})\big),K\in \big(\mcl{J}(P_{\ul{\gl}''})\bs \mcl{J}^0(P_{\ul{\gl}''})\big)}{\bigsqcup}{X^{\ul{\gl}}_{I,J,K}}}.
Then we have
\begin{enumerate}[label=(\alph*)]
	\item $X^{\ul{\gl}}_I=\R e_I+\gp\grpp$.
	\item \equan{EighteenEighteen}{\us{I\in \mcl{J}^0(P_{\ul{\gl}}) }{\sum}\mid (\autgp)_I\bs X^{\ul{\gl}}_I\mid =&\\ \us{I\in \mcl{J}^0(P_{\ul{\gl}}) }{\sum}&\bigg( \us{J\in \big(\mcl{J}(P_{\ul{\gl^{'}/I}})\bs \mcl{J}^0(P_{\ul{\gl^{'}/I}})\big),K\in \big(\mcl{J}(P_{\ul{\gl}''})\bs \mcl{J}^0(P_{\ul{\gl}''})\big)}{\sum}\frac{\mid X^{\ul{\gl}}_{I,J,K}\mid}{\ga^{\ul{\gl}}_{I,J,K}}\bigg).}
	\item $\mid (\autgp)_I\bs X^{\ul{\gl}}_I\mid = q\mid (\autgp)_I\bs \gp\grpp\mid$.
	\item If $\gl_k>1$ then \equa{\us{I\in \mcl{J}^0(P_{\ul{\gl}}) }{\sum}\mid (\autgp)_I\bs X^{\ul{\gl}}_I\mid & =q\us{I\in \mcl{J}^0(P_{\ul{\gl}}) }{\sum}\mid (\autgp)_I\bs \gp\grpp\mid= q \mid \autgp\bs (\hgrpp\times \gp \grpp)\mid\\&=q\mid \autgp\bs (\gp\hgrpp\times \gp \grpp)\mid=q\mid \autmgp\bs (\hgrmpp\times \grmpp)\mid=qn^0_{\ul{\gm}}(q).}
	\item	If $\gl_k=1$ then \equa{\us{I\in \mcl{J}^0(P_{\ul{\gl}}) }{\sum}\mid (\autgp)_I\bs X^{\ul{\gl}}_I\mid & =q\us{I\in \mcl{J}^0(P_{\ul{\gl}}) }{\sum}\mid (\autgp)_I\bs \gp\grpp\mid= q \mid \autgp\bs (\hgrpp\times \gp \grpp)\mid\\&=q\mid \autgp\bs (\gp\hgrpp\times \gp \grpp)\mid=q\mid \autmgp\bs (\grmpp\times \grmpp)\mid=qn_{\ul{\gm}}(q).}
\end{enumerate}
\end{theorem}
\begin{proof}
We prove (a). The set
\equa{X^{\ul{\gl}}_I=&\{(x',x'')\in \grppp\oplus \grpppp\mid& I(\ol{x}')=J\in \big(\mcl{J}(P_{\ul{\gl^{'}/I}})\bs \mcl{J}^0(P_{\ul{\gl^{'}/I}})\big),\\&   & I(x'')=K\in \big(\mcl{J}(P_{\ul{\gl}''})\bs \mcl{J}^0(P_{\ul{\gl}''})\big) \}\\
	&=\{(x',x'')\in \grppp\oplus \grpppp\mid& \ol{x}'\in \gp\grpppI,x''\in \gp\grpppp\}\\
	&= \{(x',x'')\in \grppp\oplus \grpppp \mid& x'\in \R e'_I+\gp\grppp,x''\in \gp\grpppp\}\\
	&=\R e_I+\gp\grpp.
}
This proves (a).

We prove (b). Consider the sum in the RHS of Equation~\ref{Eq:EighteenEighteen}. This sum is exactly given by the action of $(\autgp)_I$ on $X^{\ul{\gl}}_I$ for all $I\in \mcl{J}^0(P_{\ul{\gl}})$. This proves (b).

We prove (c). To prove this we use the fact the residue field $\R/\gp\R$ is a field with exactly $q$ elements.
Let $s:\R/\gp\R \lra \R$ be a section for the projection map $\R\lra \R/\gp\R$. Assume $s(0)=0,s(1)=1$. Given a transitive $(\autgp)_I$-orbit $O\subseteq \gp\grpp$ we have $q$ transitive $(\autgp)_I$-orbits \equ{s(a)e_I+O\subseteq \R e_I+\gp\grpp\text{ for every }a\in \R/\gp\R.} 
This is because: For any $x,y\in \gp\grpp,g\in (\autgp)_I$ we have $gx=y \Llra g(s(a)e_I+x)=s(a)e_I+y$. Also $s(a_1)e_I+O_1\cap s(a_2)e_I+O_2=\es$ unless $s(a_1)=s(a_2),O_1=O_2$ which is if and only if $a_1= a_2,O_1=O_2$. We also observe that the orbit decomposition of $X^{\ul{\gl}}_I$ is given as:
\equ{\us{a\in \R/\gp\R,O\in (\autgp)_I\bs \gp\grpp}{\bigsqcup} s(a)e_I+O=\R e_I+\gp\grpp=X^{\ul{\gl}}_I.}
This proves (c).

We prove (d) and (e). Note that if $\gl_k>1$, we have $\gp\hgrpp=\hgrmpp$. If $\gl_k=1$, we have
$\gp\hgrpp=\gp\grpp=\grmpp$. Now we use the proof of Theorem~\ref{theorem:OrbitPolytoHeightZeroPoly}. This proves (d) and (e).

Hence Theorem~\ref{theorem:FirstSummand} follows.
\end{proof}
\begin{theorem}
\label{theorem:SecondSummandFirstPart}
Let $\ul{\gl}=(\gl_1>\gl_2>\cdots>\gl_k)\in \Gl_0$ be a nonempty partition with distinct parts. Let $\ul{\gn}=(\gl_1-\gl_k>\gl_2-\gl_k>\cdots>\gl_{k-1}-\gl_k)$.
Let $I\in \mcl{J}(P_{\ul{\gl}}),J\in \mcl{J}(P_{\ul{\gl'/I}}),K\in \mcl{J}(P_{\ul{\gl}''})$.
Define $X^{\ul{\gl}}_{I,J,K}=\{(x',x'')\in \grpp\mid I(\ol{x}')=J,I(x'')=K\}$  and
$\ga^{\ul{\gl}}_{I,J,K}=\mid(\grppp)_{J\cup K}\oplus\big((\grpppp)^*_K+(\grpppp)_J\big)\mid$. Let $I_1\in \mcl{J}(P_{\ul{\gn}}),J_1\in \mcl{J}(P_{\ul{\gn'/I_1}}),K_1\in \mcl{J}(P_{\ul{\gn}''})$.
Define $X^{\ul{\gn}}_{I_1,J_1,K_1}=\{(x',x'')\in \grnpp\mid I(\ol{x}')=J_1,I(x'')=K_1\}$  and
$\ga^{\ul{\gn}}_{I_1,J_1,K_1}=\mid(\grnppp)_{J_1\cup K_1}\oplus\big((\grnpppp)^*_{K_1}+(\grnpppp)_{J_1}\big)\mid$. For any partition $\ul{\gd}$, let $\mcl{J}^0(P_{\ul{\gd}})$ be the set of height zero ideals in $\mcl{J}(P_{\ul{\gd}})$.
Then 
 \equan{TwentyThreeTwentyThree}{
 	\us{\os{I\in \mcl{J}^0(P_{\ul{\gl}})}{(0,\gl_k)\nin\max(I)}}{\sum}& \bigg(\us{(J,K)\in \mcl{J}^0(P_{\ul{\gl^{'}/I}})\times \mcl{J}(P_{\ul{\gl}''})}{\sum} \frac{\mid X^{\ul{\gl}}_{I,J,K}\mid}{\ga^{\ul{\gl}}_{I,J,K}}\bigg)\\&=\us{I_1\in \mcl{J}^0(P_{\ul{\gn}})}{\sum}\bigg( \us{(J_1,K_1)\in \mcl{J}^0(P_{\ul{\gn^{'}/I_1}})\times \mcl{J}(P_{\ul{\gn}''})}{\sum}\frac{\mid X^{\ul{\gn}}_{I_1,J_1,K_1}\mid}{\ga^{\ul{\gn}}_{I_1,J_1,K_1}}\bigg).}
\end{theorem}
\begin{proof}

The number of ideals $I\in \mcl{J}(P_{\ul{\gl}})$ of height zero, that is, $I\in \mcl{J}^0(P_{\ul{\gl}})$ is $(\gl_1-\gl_2+1)(\gl_2-\gl_3+1)\ldots(\gl_{k-1}-\gl_k+1)$ if $k>1$ and one if $k=1$. These ideals are in bijection with the ideals $I_1\in \mcl{J}(P_{\ul{\gn}})$ where $\ul{\gn}=\big((\gl_1-\gl_k)>(\gl_2-\gl_k)>\cdots>(\gl_{k-1}-\gl_k)\big)$ with the bijection being 
\equ{\max(I)\lra \max(I_1)=\{(v,\gl-\gl_k)\mid (v,\gl)\in \max(I),(v,\gl)\neq (0,\gl_k)\}.}
The inverse bijection is given as follows.
\equan{19Point5}{\max(I_1)\lra &\max(I)=\{(v,\gn+\gl_k)\mid (v,\gn)\in \max(I_1)\}\text{ if } I_1\text{ is of height zero,}\\
	\max(I_1)\lra &\max(I)=\{(v,\gn+\gl_k)\mid (v,\gn)\in \max(I_1)\}\cup\{(0,\gl_k)\}\\&\text{ if } I_1\text{ is of positive height.}}
Let $I\in \mcl{J}^0(P_{\ul{\gl}})$ with $\max(I)=\{(v_1,t_1),(v_2,t_2),\cdots,(v_s=0,t_s)\mid t_1>t_2>\cdots>t_s\}$. We have $v_1>v_2>\cdots>v_{s-1}>v_s=0$ and $t_1-v_1>t_2-v_2>\cdots>t_{s-1}-v_{s-1}>t_s-v_s=t_s\geq \gl_k$.

\equan{TwentyTwenty}{\ul{\gl}'&=(t_1>t_2>\cdots>t_{s-1}>t_s),\\ 
	\ul{\gl^{'}/I} &=\big((v_1+t_2-v_2)>(v_2+t_3-v_3)>\cdots>(v_{s-2}+t_{s-1}-v_{s-1})>(v_{s-1}+t_s)\big).} 

Consider the case $t_s>\gl_k$. Under the above bijection the ideal $I\lra I_1\in \mcl{J}(P_{\ul{\gn}})$ where
$\max(I_1)=\{(v_1,t_1-\gl_k),(v_2,t_2-\gl_k),\cdots,(v_{s-1},t_{s-1}-\gl_k),(v_s=0,t_s-\gl_k)\}$. The ideal $I_1$ is of height zero and we have

\equan{TwentyOneTwentyOne}{\ul{\gn}'&=(t_1-\gl_k>t_2-\gl_k>\cdots>t_{s-1}-\gl_k>t_s-\gl_k),\\
	\ul{\gn^{'}/I_1}&=\big((v_1+t_2-v_2-\gl_k)>(v_2+t_3-v_3-\gl_k)>\cdots\\&>(v_{s-2}+t_{s-1}-v_{s-1}-\gl_k)>(v_{s-1}+t_s-\gl_k)\big).} 

Consider the case $t_s=\gl_k$.
Here $\max(I_1)=\{(v_1,t_1-\gl_k),(v_2,t_2-\gl_k),\cdots,(v_{s-1},t_{s-1}-\gl_k)\}$ is a positive height ideal and we have
\equan{TwentyTwoTwentyTwo}{\ul{\gn}'&=(t_1-\gl_k>t_2-\gl_k>\cdots>t_{s-1}-\gl_k),\\
	\ul{\gn^{'}/I_1}&=\big((v_1+t_2-v_2-\gl_k)>(v_2+t_3-v_3-\gl_k)>\cdots\\&>(v_{s-2}+t_{s-1}-v_{s-1}-\gl_k)>v_{s-1}=v_{s-1}+t_s-\gl_k\big).} 
Note that the last part of $\ul{\gl^{'}/I}$ which is $v_{s-1}+t_s$ is always greater than $\gl_k$. 
\begin{center}
	*****************************************************************
\end{center}

Suppose $I\in \mcl{J}^0(P_{\ul{\gl}})$ is a height zero ideal, such that $t_s>\gl_k$, which occurs if and only if $I_1\in \mcl{J}^0(P_{\ul{\gn}})$ is a height zero ideal, then $\gl_k$ is a part of $\ul{\gl}''$ and it is the last part of $\ul{\gl}''$. So the set of ideals of height zero, that is, the set $\mcl{J}^0(P_{\ul{\gl}''})$ is in bijection with the set $\mcl{J}(P_{\ul{\gn}''})$ of all ideals.

Let $J\in \mcl{J}^0(P_{\ul{\gl^{'}/I}}),K\in \mcl{J}(P_{\ul{\gl}''})$
and $\max(J)=\{(w_1,r_1),(w_2,r_2),\cdots,$ $(w_l=0,r_l)\}$. Then $r_l-w_l=r_l\geq v_{s-1}+t_s>\gl_k$ and therefore for $1\leq i\leq l$ we have $r_i-w_i\geq r_l-w_l=r_l>\gl_k$. Let $J_1\in \mcl{J}(P_{\ul{\gn^{'}/I_1}})$ be such that $\max(J_1)=\{(w_1,r_1-\gl_k),(w_2,r_2-\gl_k),\cdots,(w_l=0,r_l-\gl_k)\}$. Then $J_1\in \mcl{J}^0(P_{\ul{\gn^{'}/I_1}})$. Here we observe in this way, that there is a bijection between the sets $\mcl{J}^0(P_{\ul{\gl^{'}/I}})$ and $\mcl{J}^0(P_{\ul{\gn^{'}/I_1}})$.

Suppose $K$ is a height zero ideal and $(0,\gl_k)\nin \max(K)$. Let $K_1\in \mcl{J}(P_{\ul{\gn}''})$ be the height zero ideal associated to $K$ under the bijection $\mcl{J}^0(P_{\ul{\gl}''})\llra \mcl{J}(P_{\ul{\gn}''})$. 

We have in this scenario, \equa{\mid X^{\ul{\gl}}_{I,J,K}\mid&=\frac{\mid \grppp\mid}{\mid \grpppI\mid}\mid (\grpppI)^*_J\mid\mid (\grpppp)^*_K\mid\\
	&=q^{\gl_k}\frac{\mid \grnppp\mid}{\mid \grnpppI\mid}q^{(s-1)\gl_k}\mid (\grnpppI)^*_{J_1}\mid q^{(k-s)\gl_k}\mid (\grnpppp)^*_{K_1}\mid\\&=q^{k\gl_k}\mid X^{\ul{\gn}}_{I_1,J_1,K_1}\mid.}
Also \equa{\ga^{\ul{\gl}}_{I,J,K}&=\mid (\grppp)_{J\cup K}\mid\mid (\grpppp)^*_K+(\grpppp)_J\mid \\ &=
	q^{s\gl_k}\mid (\grnppp)_{J_1\cup K_1}\mid q^{(k-s)\gl_k} \mid (\grnpppp)^*_{K_1}+(\grnpppp)_{J_1}\mid\\&= q^{k\gl_k} \ga^{\ul{\gn}}_{I_1,J_1,K_1}.}
So we have 
\equ{\frac{\mid X^{\ul{\gl}}_{I,J,K}\mid}{\ga^{\ul{\gl}}_{I,J,K}}=\frac{\mid X^{\ul{\gn}}_{I_1,J_1,K_1}\mid}{\ga^{\ul{\gn}}_{I_1,J_1,K_1}}.}

Now suppose $(0,\gl_k)\in \max(K)$. Let $K_1\in \mcl{J}(P_{\ul{\gn}''})$ be the ideal of positive height associated to $K$ under the bijection $\mcl{J}^0(P_{\ul{\gl}''})\llra \mcl{J}(P_{\ul{\gn}''})$. Let $K'\in \mcl{J}(P_{\ul{\gl}})$ be such that $\max(K')=\max(K)\bs\{(0,\gl_k)\}$. Then $K'$ is an ideal of positive height. Any positive height ideal $K_2\in \mcl{J}(P_{\ul{\gl}})$ satisfies the condition, $K'\subseteq K_2\sbnq K$ for some unique $K\in \mcl{J}(P_{\ul{\gl}})$ such that $(0,\gl_k)\in \max(K)$.
The ideal $K$ is obtained from $K_2$ as follows $K=K_2\cup [\langle \{(0,\gl_k)\}\rangle]_{\ul{\gl}''}$. Now it is clear that for this ideal $K$ we have $K'\subseteq K_2\sbnq K$ and $K$ is uniquely determined from $K_2$.

Now we observe that for the ideal $K_2$ we have $\max(K_2)\bs [J]_{\ul{\gl}''}=\max(K)\bs [J]_{\ul{\gl}''}=\max(K')\bs [J]_{\ul{\gl}''}$.
Moreover we have the ideals $J\cup K_2=J\cup K=J\cup K'$ in the big fundamental poset $P$ as defined in~\ref{Eq:FP}. Hence we have 
\equa{\ga^{\ul{\gl}}_{I,J,K}&=\ga^{\ul{\gl}}_{I,J,K'}=\ga^{\ul{\gl}}_{I,J,K_2}\\
	\ga^{\ul{\gl}}_{I,J,K}&=\mid (\grppp)_{J\cup K}\mid\mid (\grpppp)^*_K+(\grpppp)_J\mid \\ &=
	q^{s\gl_k}\mid (\grnppp)_{J_1\cup K_1}\mid q^{(k-s)\gl_k} \mid (\grnpppp)^*_{K_1}+(\grnpppp)_{J_1}\mid\\&= q^{k\gl_k} \ga^{\ul{\gn}}_{I_1,J_1,K_1}.}
We also observe that 
\equa{\mid\us{K'\subseteq K_2\sbnq K}{\sqcup} (\grpppp)^*_{K_2}\mid&=q^{(k-s)\gl_k}\mid (\grnpppp)^*_{K_1}\mid\\
	\mid\us{K'\subseteq K_2\sbnq K}{\sqcup} X^{\ul{\gl}}_{I,J,K_2}\mid &=q^{k\gl_k}\mid X^{\ul{\gn}}_{I_1,J_1,K_1}\mid.}
This proves that 
\equ{\us{K'\subseteq K_2\subseteq K}{\sum}\frac{\mid X^{\ul{\gl}}_{I,J,K_2}\mid}{\ga^{\ul{\gl}}_{I,J,K_2}}=\frac{\mid \us{K'\subseteq K_2\subseteq K}{\sqcup} X^{\ul{\gl}}_{I,J,K_2}\mid}{\ga^{\ul{\gl}}_{I,J,K}}=\frac{\mid X^{\ul{\gn}}_{I_1,J_1,K_1}\mid}{\ga^{\ul{\gn}}_{I_1,J_1,K_1}}.}

Hence we have proved the identity in Equation~\ref{Eq:TwentyThreeTwentyThree}. This proves Theorem~\ref{theorem:SecondSummandFirstPart}.
\end{proof}
\begin{theorem}
\label{theorem:ThirdSummandFirstPart}
Let $\ul{\gl}=(\gl_1>\gl_2>\cdots>\gl_k)\in \Gl_0$ be a nonempty partition with distinct parts. Let $\ul{\gn}=(\gl_1-\gl_k>\gl_2-\gl_k>\cdots>\gl_{k-1}-\gl_k)$.
Let $I\in \mcl{J}(P_{\ul{\gl}}),J\in \mcl{J}(P_{\ul{\gl'/I}}),K\in \mcl{J}(P_{\ul{\gl}''})$.
Define $X^{\ul{\gl}}_{I,J,K}=\{(x',x'')\in \grpp\mid I(\ol{x}')=J,I(x'')=K\}$  and
$\ga^{\ul{\gl}}_{I,J,K}=\mid(\grppp)_{J\cup K}\oplus\big((\grpppp)^*_K+(\grpppp)_J\big)\mid$. Let $I_1\in \mcl{J}(P_{\ul{\gn}}),J_1\in \mcl{J}(P_{\ul{\gn'/I_1}}),K_1\in \mcl{J}(P_{\ul{\gn}''})$.
Define $X^{\ul{\gn}}_{I_1,J_1,K_1}=\{(x',x'')\in \grnpp\mid I(\ol{x}')=J_1,I(x'')=K_1\}$  and
$\ga^{\ul{\gn}}_{I_1,J_1,K_1}=\mid(\grnppp)_{J_1\cup K_1}\oplus\big((\grnpppp)^*_{K_1}+(\grnpppp)_{J_1}\big)\mid$. For any partition $\ul{\gd}$, let $\mcl{J}^0(P_{\ul{\gd}})$ be the set of height zero ideals in $\mcl{J}(P_{\ul{\gd}})$.
Then 
\equan{TwentySixTwentySix}{
	\us{\os{I\in \mcl{J}^0(P_{\ul{\gl}})}{(0,\gl_k)\in\max(I)}}{\sum}& \bigg(\us{J\in \mcl{J}^0(P_{\ul{\gl^{'}/I}}), K\in \mcl{J}(P_{\ul{\gl}''})}{\sum} \frac{\mid X^{\ul{\gl}}_{I,J,K}\mid}{\ga^{\ul{\gl}}_{I,J,K}}\bigg)\\&=\us{I_1\in \mcl{J}(P_{\ul{\gn}})\bs \mcl{J}^0(P_{\ul{\gn}})}{\sum}\bigg( \us{J_1\in \mcl{J}^0(P_{\ul{\gn^{'}/I_1}}),K_1\in \mcl{J}(P_{\ul{\gn}''})}{\sum}\frac{\mid X^{\ul{\gn}}_{I_1,J_1,K_1}\mid}{\ga^{\ul{\gn}}_{I_1,J_1,K_1}}\bigg).}

\end{theorem}
\begin{proof}
Here we observe that 
$I\in \mcl{J}^0(P_{\ul{\gl}})$ is a height zero ideal, such that $t_s=\gl_k$ occurs if and only if $I_1\in \mcl{J}(P_{\ul{\gn}})\bs \mcl{J}^0(P_{\ul{\gn}})$ is an ideal of positive height. The set $\mcl{J}^0(P_{\ul{\gl^{'}/I}})$ of ideals of height zero is in bijection with the set $\mcl{J}^0(P_{\ul{\gn^{'}/I_1}})$ of height zero ideals. Here $\gl_k$ is a part of $\ul{\gl}'$ and it is the last part of $\ul{\gl}'$. So the last part of $\ul{\gl}''$ is greater than $\gl_k$.  Let $\ul{\gl}'''$ be the partition obtained by adding $\gl_k$ to the partition $\ul{\gl}''$ at the end.	The set  $\mcl{J}^0(P_{\ul{\gl}'''})$ of ideals of height zero is is in bijection with the set $\mcl{J}(P_{\ul{\gn}''})$ of all ideals. 

Let $J\in \mcl{J}^0(P_{\ul{\gl^{'}/I}}),K\in \mcl{J}(P_{\ul{\gl}''})$
and $\max(J)=\{(w_1,r_1),(w_2,r_2),\cdots,$ $(w_l=0,r_l)\}$. Then $r_l-w_l=r_l\geq v_{s-1}+t_s>\gl_k$ and therefore for $1\leq i\leq l$ we have $r_i-w_i\geq r_l-w_l=r_l>\gl_k$. Let $J_1\in \mcl{J}(P_{\ul{\gn^{'}/I_1}})$ be such that $\max(J_1)=\{(w_1,r_1-\gl_k),(w_2,r_2-\gl_k),\cdots,(w_l=0,r_l-\gl_k)\}$. Then $J_1\in \mcl{J}^0(P_{\ul{\gn^{'}/I_1}})$. Here we observe in this way, that there is a bijection between the sets $\mcl{J}^0(P_{\ul{\gl^{'}/I}})$ and $\mcl{J}^0(P_{\ul{\gn^{'}/I_1}})$.

Let $K'\in \mcl{J}^0(P_{\ul{\gl}'''})$ be such that $(0,\gl_k)\nin \max(K')$. Then $K=K'\cap P_{\ul{\gl}''}\in \mcl{J}^0(P_{\ul{\gl}''})$. Let $\max(K')=\max(K)=\{(x_1,a_1),\cdots,(x_m=0,a_m)\}$. Let $K_1\in \mcl{J}^0(P_{\ul{\gn}''})$ be such that $\max(K_1)=\{(x_1,a_1-\gl_k),(x_2,a_2-\gl_k),\cdots,(x_m=0,a_m-\gl_k)\}$. Then it is clear that 

\equ{\frac{\mid X^{\ul{\gl}}_{I,J,K}\mid}{\ga^{\ul{\gl}}_{I,J,K}}=\frac{\mid X^{\ul{\gn}}_{I_1,J_1,K_1}\mid}{\ga^{\ul{\gn}}_{I_1,J_1,K_1}}.}

Let $K'\in \mcl{J}^0(P_{\ul{\gl}'''})$ be such that $(0,\gl_k)\in \max(K')$. Let $K\in \mcl{J}(P_{\ul{\gl}''})$ be such that $\max(K)=\max(K')\bs\{(0,\gl_k)\}$. Let $K_1\in \mcl{J}(P_{\ul{\gn}''})\bs\mcl{J}^0(P_{\ul{\gn}''})$ be the ideal of positive height.
Now given any ideal $K_2\in \mcl{J}(P_{\ul{\gl}''})$ of positive height, there is a unique ideal $K'\in \mcl{J}^0(P_{\ul{\gl}'''})$ be such that $(0,\gl_k)\in \max(K')$ satisfying $K\subseteq K_2\subseteq K'\bs P_{(\gl_k)}$ where $K\in \mcl{J}(P_{\ul{\gl}''})$ is such that $\max(K)=\max(K')\bs\{(0,\gl_k)\}$.
  Then it is clear that 

\equ{\us{K\subseteq K_2\subseteq K'\bs P_{(\gl_k)}}{\sum}\frac{\mid X^{\ul{\gl}}_{I,J,K_2}\mid}{\ga^{\ul{\gl}}_{I,J,K_2}}=\frac{\mid \us{K\subseteq K_2\subseteq K'\bs P_{(\gl_k)}}{\sqcup} X^{\ul{\gl}}_{I,J,K_2}\mid}{\ga^{\ul{\gl}}_{I,J,K}}=\frac{\mid X^{\ul{\gn}}_{I_1,J_1,K_1}\mid}{\ga^{\ul{\gn}}_{I_1,J_1,K_1}}.}

This proves the identity in Equation~\ref{Eq:TwentySixTwentySix}. Hence the theorem follows.
\end{proof}
\begin{theorem}
	\label{theorem:SecondSummandSecondPart}
	Let $\ul{\gl}=(\gl_1>\gl_2>\cdots>\gl_k)\in \Gl_0$ be a nonempty partition with distinct parts. Let $\ul{\gn}=(\gl_1-\gl_k>\gl_2-\gl_k>\cdots>\gl_{k-1}-\gl_k)$.
	Let $I\in \mcl{J}(P_{\ul{\gl}}),J\in \mcl{J}(P_{\ul{\gl'/I}}),K\in \mcl{J}(P_{\ul{\gl}''})$.
	Define $X^{\ul{\gl}}_{I,J,K}=\{(x',x'')\in \grpp\mid I(\ol{x}')=J,I(x'')=K\}$  and
	$\ga^{\ul{\gl}}_{I,J,K}=\mid(\grppp)_{J\cup K}\oplus\big((\grpppp)^*_K+(\grpppp)_J\big)\mid$. Let $I_1\in \mcl{J}(P_{\ul{\gn}}),J_1\in \mcl{J}(P_{\ul{\gn'/I_1}}),K_1\in \mcl{J}(P_{\ul{\gn}''})$.
	Define $X^{\ul{\gn}}_{I_1,J_1,K_1}=\{(x',x'')\in \grnpp\mid I(\ol{x}')=J_1,I(x'')=K_1\}$  and
	$\ga^{\ul{\gn}}_{I_1,J_1,K_1}=\mid(\grnppp)_{J_1\cup K_1}\oplus\big((\grnpppp)^*_{K_1}+(\grnpppp)_{J_1}\big)\mid$. For any partition $\ul{\gd}$, let $\mcl{J}^0(P_{\ul{\gd}})$ be the set of height zero ideals in $\mcl{J}(P_{\ul{\gd}})$.
	Then 
\equan{TwentyFourTwentyFour}{
	\us{\os{I\in \mcl{J}^0(P_{\ul{\gl}})}{(0,\gl_k)\nin\max(I)}}{\sum}& \bigg(\us{J\in \mcl{J}(P_{\ul{\gl^{'}/I}})\bs\mcl{J}^0(P_{\ul{\gl^{'}/I}}), K\in \mcl{J}^0(P_{\ul{\gl}''})}{\sum} \frac{\mid X^{\ul{\gl}}_{I,J,K}\mid}{\ga^{\ul{\gl}}_{I,J,K}}\bigg)\\&=\us{I_1\in \mcl{J}^0(P_{\ul{\gn}})}{\sum}\bigg( \us{J_1\in \mcl{J}(P_{\ul{\gn^{'}/I_1}})\bs\mcl{J}^0(P_{\ul{\gn^{'}/I_1}}),K_1\in \mcl{J}(P_{\ul{\gn}''})}{\sum}\frac{\mid X^{\ul{\gn}}_{I_1,J_1,K_1}\mid}{\ga^{\ul{\gn}}_{I_1,J_1,K_1}}\bigg).}
	
\end{theorem}
\begin{proof}
Suppose $I\in \mcl{J}^0(P_{\ul{\gl}})$ is a height zero ideal, such that $t_s>\gl_k$, which occurs if and only if $I_1\in \mcl{J}^0(P_{\ul{\gn}})$ is a height zero ideal, then $\gl_k$ is a part of $\ul{\gl}''$ and it is the last part of $\ul{\gl}''$. So the set of ideals of height zero, that is, the set $\mcl{J}^0(P_{\ul{\gl}''})$ is in bijection with the set $\mcl{J}(P_{\ul{\gn}''})$ of all ideals. Let $\ul{\gl}'''$ be the partition obtained by adding the part $\gl_k$ to $\ul{\gl^{'}/I}$. Then the set of all ideals $J'\in \mcl{J}(P_{\ul{\gl}'''})$ of height zero such that $(0,\gl_k)\in \max(J')$ is in bijection with $\mcl{J}(P_{\ul{\gn^{'}/I_1}})\bs\mcl{J}^0(P_{\ul{\gn^{'}/I_1}})$.  Now the rest of the proof follows using standard arguments. This proves the theorem.
\end{proof}
\begin{theorem}
	\label{theorem:ThirdSummandSecondPart}
	Let $\ul{\gl}=(\gl_1>\gl_2>\cdots>\gl_k)\in \Gl_0$ be a nonempty partition with distinct parts. Let $\ul{\gn}=(\gl_1-\gl_k>\gl_2-\gl_k>\cdots>\gl_{k-1}-\gl_k)$.
	Let $I\in \mcl{J}(P_{\ul{\gl}}),J\in \mcl{J}(P_{\ul{\gl'/I}}),K\in \mcl{J}(P_{\ul{\gl}''})$.
	Define $X^{\ul{\gl}}_{I,J,K}=\{(x',x'')\in \grpp\mid I(\ol{x}')=J,I(x'')=K\}$  and
	$\ga^{\ul{\gl}}_{I,J,K}=\mid(\grppp)_{J\cup K}\oplus\big((\grpppp)^*_K+(\grpppp)_J\big)\mid$. Let $I_1\in \mcl{J}(P_{\ul{\gn}}),J_1\in \mcl{J}(P_{\ul{\gn'/I_1}}),K_1\in \mcl{J}(P_{\ul{\gn}''})$.
	Define $X^{\ul{\gn}}_{I_1,J_1,K_1}=\{(x',x'')\in \grnpp\mid I(\ol{x}')=J_1,I(x'')=K_1\}$  and
	$\ga^{\ul{\gn}}_{I_1,J_1,K_1}=\mid(\grnppp)_{J_1\cup K_1}\oplus\big((\grnpppp)^*_{K_1}+(\grnpppp)_{J_1}\big)\mid$. For any partition $\ul{\gd}$, let $\mcl{J}^0(P_{\ul{\gd}})$ be the set of height zero ideals in $\mcl{J}(P_{\ul{\gd}})$.
	Then 
\equan{TwentySevenTwentySeven}{
	\us{\os{I\in \mcl{J}^0(P_{\ul{\gl}})}{(0,\gl_k)\in\max(I)}}{\sum}& \bigg(\us{J\in \mcl{J}(P_{\ul{\gl^{'}/I}})\bs\mcl{J}^0(P_{\ul{\gl^{'}/I}}), K\in \mcl{J}^0(P_{\ul{\gl}''})}{\sum} \frac{\mid X^{\ul{\gl}}_{I,J,K}\mid}{\ga^{\ul{\gl}}_{I,J,K}}\bigg)\\&=\us{I_1\in \mcl{J}(P_{\ul{\gn}})\bs\mcl{J}^0(P_{\ul{\gn}})}{\sum}\bigg( \us{J_1\in \mcl{J}(P_{\ul{\gn^{'}/I_1}})\bs\mcl{J}^0(P_{\ul{\gn^{'}/I_1}}),K_1\in \mcl{J}^0(P_{\ul{\gn}''})}{\sum}\frac{\mid X^{\ul{\gn}}_{I_1,J_1,K_1}\mid}{\ga^{\ul{\gn}}_{I_1,J_1,K_1}}\bigg).}
	
\end{theorem}
\begin{proof}
Here we observe that 
$I\in \mcl{J}^0(P_{\ul{\gl}})$ is a height zero ideal, such that $t_s=\gl_k$ occurs if and only if $I_1\in \mcl{J}(P_{\ul{\gn}})\bs \mcl{J}^0(P_{\ul{\gn}})$ is an ideal of positive height. Here $\gl_k$ is a part of $\ul{\gl}'$ and it is the last part of $\ul{\gl}'$. So the last part of $\ul{\gl}''$ is greater than $\gl_k$.  The set $\mcl{J}^0(P_{\ul{\gl}''})$ of ideals of height zero is in bijection with the set $\mcl{J}^0(P_{\ul{\gn}''})$ of height zero ideals.
Let $\ul{\gl}'''$ be the partition obtained by adding $\gl_k$ to the partition $\ul{\gl}''$ at the end.	The set  $\mcl{J}^0(P_{\ul{\gl}'''})$ of ideals $J'$ of height zero such that $(0,\gl_k)\in \max(J')$ is in bijection with the set $\mcl{J}(P_{\ul{\gn}''})$ of all ideals. Now the rest of the proof follows using standard arguments. This proves the theorem.
\end{proof}
\begin{proof}[Proof of Theorem~\ref{theorem:HeightZeroPartitiontoHeightZeroSmallerPartitions}]
We prove (1),(2) simultaneously now. First we have using Equation~\ref{Eq:Alternative}
\equ{n^0_{\ul{\gl}}(q)=\us{I\in \mcl{J}^0(P_{\ul{\gl}}) }{\sum}\mid (\autgp)_I\bs\grpp\mid= \us{I\in \mcl{J}^0(P_{\ul{\gl}}) }{\sum}\bigg( \us{J\in \mcl{J}(P_{\ul{\gl^{'}/I}}),K\in \mcl{J}(P_{\ul{\gl}''})}{\sum}\frac{\mid X^{\ul{\gl}}_{I,J,K}\mid}{\ga^{\ul{\gl}}_{I,J,K}}\bigg).}
The above sum splits as two sums, one of which is given in the RHS of Equation~\ref{Eq:EighteenEighteen} and the other one is given below.  
\equan{NineteenNineteen}{ \us{I\in \mcl{J}^0(P_{\ul{\gl}}) }{\sum}\bigg( \us{(J,K)\in \mcl{J}^0(P_{\ul{\gl^{'}/I}})\times \mcl{J}(P_{\ul{\gl}''}) \bigcup \mcl{J}(P_{\ul{\gl^{'}/I}})\times \mcl{J}^0(P_{\ul{\gl}''})}{\sum}\frac{\mid X^{\ul{\gl}}_{I,J,K}\mid}{\ga^{\ul{\gl}}_{I,J,K}}\bigg).}

Theorem~\ref{theorem:FirstSummand} then gives the first summand in Theorem~\ref{theorem:HeightZeroPartitiontoHeightZeroSmallerPartitions} which is $qn^0_{\ul{\gm}}(q)$ if $\gl_k>1$ and $qn_{\ul{\gm}}(q)$ if $\gl_k=1$. 
We need to prove that the remaining sum in~\ref{Eq:NineteenNineteen} is  $n^0_{\ul{\gn}}(q)+n^0_{\ul{\gd}}(q)$ if $\gl_{k-1}-\gl_k>1$ and 
$n^0_{\ul{\gn}}(q)+n_{\ul{\gd}}(q)$ if $\gl_{k-1}-\gl_k=1$.

Now we break the subsum given in~\ref{Eq:NineteenNineteen} into two subsums as:

\equa{ \us{I\in \mcl{J}^0(P_{\ul{\gl}}) }{\sum}&\bigg( \us{(J,K)\in \mcl{J}^0(P_{\ul{\gl^{'}/I}})\times \mcl{J}(P_{\ul{\gl}''}) \bigcup \mcl{J}(P_{\ul{\gl^{'}/I}})\times \mcl{J}^0(P_{\ul{\gl}''})}{\sum}\frac{\mid X^{\ul{\gl}}_{I,J,K}\mid}{\ga^{\ul{\gl}}_{I,J,K}}\bigg)\\
	&=\us{\os{I\in \mcl{J}^0(P_{\ul{\gl}})}{(0,\gl_k)\nin\max(I)} }{\sum}\bigg( \us{(J,K)\in \mcl{J}^0(P_{\ul{\gl^{'}/I}})\times \mcl{J}(P_{\ul{\gl}''}) \bigcup \mcl{J}(P_{\ul{\gl^{'}/I}})\times \mcl{J}^0(P_{\ul{\gl}''})}{\sum}\frac{\mid X^{\ul{\gl}}_{I,J,K}\mid}{\ga^{\ul{\gl}}_{I,J,K}}\bigg)\\
	&+\us{\os{I\in \mcl{J}^0(P_{\ul{\gl}})}{(0,\gl_k)\in\max(I)} }{\sum}\bigg( \us{(J,K)\in \mcl{J}^0(P_{\ul{\gl^{'}/I}})\times \mcl{J}(P_{\ul{\gl}''}) \bigcup \mcl{J}(P_{\ul{\gl^{'}/I}})\times \mcl{J}^0(P_{\ul{\gl}''})}{\sum}\frac{\mid X^{\ul{\gl}}_{I,J,K}\mid}{\ga^{\ul{\gl}}_{I,J,K}}\bigg).}

To complete the proof of Theorem~\ref{theorem:HeightZeroPartitiontoHeightZeroSmallerPartitions}, we need the identities given by Equation
~\ref{Eq:TwentyFiveTwentyFive} and Equation~\ref{Eq:TwentyEightTwentyEight}. 

\equan{TwentyFiveTwentyFive}{\us{\os{I\in \mcl{J}^0(P_{\ul{\gl}})}{(0,\gl_k)\nin\max(I)}}{\sum}& \bigg(\us{(J,K)\in \mcl{J}^0(P_{\ul{\gl^{'}/I}})\times \mcl{J}(P_{\ul{\gl}''}) \bigcup \mcl{J}(P_{\ul{\gl^{'}/I}})\times \mcl{J}^0(P_{\ul{\gl}''})}{\sum} \frac{\mid X^{\ul{\gl}}_{I,J,K}\mid}{\ga^{\ul{\gl}}_{I,J,K}}\bigg)\\&=\us{I_1\in \mcl{J}^0(P_{\ul{\gn}})}{\sum}\bigg( \us{J_1\in \mcl{J}(P_{\ul{\gn^{'}/I_1}}),K_1 \in \mcl{J}(P_{\ul{\gn}''})}{\sum}\frac{\mid X^{\ul{\gn}}_{I_1,J_1,K_1}\mid}{\ga^{\ul{\gn}}_{I_1,J_1,K_1}}\bigg)=n^0_{\ul{\gn}}(q).} 

\equan{TwentyEightTwentyEight}{\us{\os{I\in \mcl{J}^0(P_{\ul{\gl}})}{(0,\gl_k)\in\max(I)}}{\sum}& \bigg(\us{(J,K)\in \mcl{J}^0(P_{\ul{\gl^{'}/I}})\times \mcl{J}(P_{\ul{\gl}''}) \bigcup \mcl{J}(P_{\ul{\gl^{'}/I}})\times \mcl{J}^0(P_{\ul{\gl}''})}{\sum} \frac{\mid X^{\ul{\gl}}_{I,J,K}\mid}{\ga^{\ul{\gl}}_{I,J,K}}\bigg)\\&=\us{I_1\in \mcl{J}(P_{\ul{\gn}})\bs \mcl{J}^0(P_{\ul{\gn}})}{\sum}\bigg( \us{(J_1,K_1)\in \mcl{J}^0(P_{\ul{\gn^{'}/I_1}})\times \mcl{J}(P_{\ul{\gn}''}) \bigcup \mcl{J}(P_{\ul{\gn^{'}/I_1}})\times \mcl{J}^0(P_{\ul{\gn}''})}{\sum}\frac{\mid X^{\ul{\gn}}_{I_1,J_1,K_1}\mid}{\ga^{\ul{\gn}}_{I_1,J_1,K_1}}\bigg)\\&=\mid \autngp\bs(\gp\grnpp\times \hgrnpp)\mid.} 
Observe that the last part of $\ul{\gn}$ is $\gl_{k-1}-\gl_k$. We have already observed that if $\gl_{k-1}-\gl_k>1$ then $\mid \autngp\bs(\gp\grnpp\times \hgrnpp)\mid=n^0_{\ul{\gd}}(q)$ and if $\gl_{k-1}-\gl_k=1$ then $\mid \autngp\bs(\gp\grnpp\times \hgrnpp)\mid=n_{\ul{\gd}}(q)$ where $\gp\grnpp\cong \mcl{A}_{\ul{\gd}}$.

Equation~\ref{Eq:TwentyFiveTwentyFive} is obtained by combining Equation~\ref{Eq:TwentyThreeTwentyThree} and Equation~\ref{Eq:TwentyFourTwentyFour}. Equation~\ref{Eq:TwentyEightTwentyEight} is obtained by combining Equation~\ref{Eq:TwentySixTwentySix} and Equation~\ref{Eq:TwentySevenTwentySeven}.

We established the identities in Equations~\ref{Eq:TwentyThreeTwentyThree},~\ref{Eq:TwentyFourTwentyFour},~\ref{Eq:TwentySixTwentySix},~\ref{Eq:TwentySevenTwentySeven} in Theorems~\ref{theorem:SecondSummandFirstPart},\linebreak ~\ref{theorem:SecondSummandSecondPart},~\ref{theorem:ThirdSummandFirstPart}~\ref{theorem:ThirdSummandSecondPart} respectively. 

Theorem~\ref{theorem:HeightZeroPartitiontoHeightZeroSmallerPartitions}(3) is easy.
Therefore we have completed the proof of Theorem~\ref{theorem:HeightZeroPartitiontoHeightZeroSmallerPartitions}.

\end{proof}
\begin{remark}
Theorem~\ref{theorem:HeightZeroPartitiontoHeightZeroSmallerPartitions} can be generalized to partitions which need not have distinct parts.

Also many more  nice identities such as the following have been observed.
Let $\ul{\gl}=(\gl_1>\gl_2>\cdots>\gl_k)$ be a partition with distinct parts.
\begin{itemize}
	\item Suppose $k\geq 2,\gl_k=1,\gl_1-\gl_2>1$. Let $\ul{\gm}=(\gl_1-1>\gl_2>\cdots>\gl_k),\ul{\gn}=(\gl_2>\cdots>\gl_k)$.
	Then we have \equ{n^0_{\ul{\gl}}(q)=qn^0_{\ul{\gm}}(q)+2n^0_{\ul{\gn}}(q).}
	\item Suppose $k\geq 2,\gl_k=1,\gl_1-\gl_2=1$. Let $\ul{\gm}=(\gl_2^2>\gl_3>\cdots>\gl_k),\ul{\gn}=(\gl_2>\gl_3>\cdots>\gl_k),\ul{\gd}=
	\big((\gl_2-1)^2\geq \gl_3>\cdots>\gl_k\big)$. 
	Then we have \equ{n^0_{\ul{\gl}}(q)=qn^0_{\ul{\gm}}(q)+2n^0_{\ul{\gn}}(q)+n^0_{\ul{\gd}}(q).}
\end{itemize}
\end{remark}
\begin{theorem}
	\label{theorem:PositivityConjecture}
	Let $\ugl\in \Gl_0$ be a partition and $\R$ be a discrete valuation ring with maximal ideal generated by a uniformizing element $\gp$ having finite residue field ${\bf k}=\RR {} \cong \mbb{F}_q$. Let $\grpp=\us{i=1}{\os{k}{\oplus}} (\RR {\gl_i})^{\gr_i}$ be its associated finite $\R$-module and $\autgp$ be its automorphism group. Let $n_{\ul{\gl}}(q)$ be the number of orbits of pairs in $\grpp \times \grpp$ for the diagonal action of $\autgp$ on $\grpp \times \grpp$. Then $n_{\ul{\gl}}(q)$ is a polynomial of degree $\gl_1$ with non-negative integer coefficients.
\end{theorem}
\begin{proof}
This theorem follows from Theorems~[\ref{theorem:OrbitPolytoHeightZeroPoly},~\ref{theorem:HeightZeroPolyDistinctParts},~\ref{theorem:HeightZeroPartitiontoHeightZeroSmallerPartitions}] and the fact that $n^0_{(1)}(q)=q,n_{\es}(q)=n^0_{\es}(q)=1$. 
\end{proof}

\end{document}